\documentclass[11pt]{article}

\usepackage{amsmath,amssymb,amsfonts,enumitem,amsthm,accents,pstricks,pst-node,pstricks-add,mathrsfs,xy}
\xyoption{all}

\usepackage{color}
\usepackage{makeidx}
\usepackage{graphicx}
\usepackage{setspace}
\usepackage{extarrows}
\usepackage{tocloft}

\usepackage[left=1.1in, right=1.1in, top=1in, bottom=1in]{geometry}
\usepackage{lipsum}
\usepackage{titlesec}
\numberwithin{equation}{section}
\newtheorem{theorem}{Theorem}[section]
\newtheorem{lemma}[theorem]{Lemma}
\newtheorem{conjlemma}[theorem]{Conjectural Lemma}
\newtheorem{corollary}[theorem]{Corollary}
\newtheorem{proposition}[theorem]{Proposition}
\newtheorem{conjecture}[theorem]{Conjecture}

\theoremstyle{definition}
\newtheorem{definition}[theorem]{Definition}
\newtheorem{remark}[theorem]{Remark}
\newtheorem{example}[theorem]{Example}

\newtheorem{claim}[theorem]{Claim}

\def\ov#1{\overline{#1}}
\def\tn#1{\textnormal{#1}}
\def\mf#1{\mathfrak{#1}}
\def\wt#1{\widetilde{#1}}
\def\wh#1{\widehat{#1}}
\def\scz{\scriptsize}
\def\vr{\varrho}
\def\ll{\left\langle}
\def\rr{\right\rangle}
\def\mc{\mathcal}
\def\lra{\longrightarrow}
\def\dbar{\bar\partial}

\newcommand\uvec[1]{\underaccent{\vec}{#1}}
\newcommand{\ucev}[1]{\reflectbox{\ensuremath{\uvec{\reflectbox{\ensuremath{#1}}}}}}

\def\bEq#1{\begin{equation}\label{#1}}
\def\eEq{\end{equation}}

\def\bsEq{\begin{equation*}}
\def\esEq{\end{equation*}}

\def\bDf#1{\begin{definition}\label{#1}}
\def\eDf{\end{definition}}

\def\bTh#1{\begin{theorem}\label{#1}}
\def\eTh{\end{theorem}}

\def\bCn#1{\begin{conjecture}\label{#1}}
\def\eCn{\end{conjecture}}

\def\bClm#1{\begin{claim}\label{#1}}
\def\eClm{\end{claim}}

\def\bLm#1{\begin{lemma}\label{#1}}
\def\eLm{\end{lemma}}

\def\bCLm#1{\begin{conjlemma}\label{#1}}
\def\eCLm{\end{conjlemma}}

\def\bRm#1{\begin{remark}\label{#1}}
\def\eRm{\end{remark}}

\def\bEx#1{\begin{example}\label{#1}}
\def\eEx{\end{example}}

\def\bPr#1{\begin{proposition}\label{#1}}
\def\ePr{\end{proposition}}

\def\bCr#1{\begin{corollary}\label{#1}}
\def\eCr{\end{corollary}}

\def\bFg#1{\begin{figure}\label{#1}}
\def\eFg{\end{figure}}

\def\bPf{\begin{proof}}
\def\ePf{\end{proof}}

\def\bIt{\begin{itemize}[leftmargin=*]}
\def\eIt{\end{itemize}}

\def\bEn{\begin{enumerate}[label=$(\arabic*)$,leftmargin=*]}
\def\eEn{\end{enumerate}}

\def\bEnalph{\begin{enumerate}[label=$(\alph*)$,leftmargin=*]}
\def\eEnalph{\end{enumerate}}

\def\AK{\tn{AK}}
\def\coker{\tn{coker}}

\def\nd{\tn{d}}

\def\Def{\tn{Def}}
\def\ev{\tn{ev}}
\def\GW{\tn{GW}}
\def\id{\tn{id}}
\def\ker{\tn{ker}}

\def\Obs{\tn{Obs}}
\def\ord{\tn{ord}}

\def\st{\tn{st}}
\def\Symp{\tn{Symp}}

\def\cB{\mc{B}}
\def\cC{\mc{C}}

\def\cM{\mc{M}}
\def\cO{\mc{O}}
\def\cE{\mc{E}}

\def\cN{\mc{N}}
\def\cR{\mc{R}}

\def\cL{\mc{L}}

\def\cU{\mc{U}}
\def\cC{\mc{C}}
\def\cH{\mc{H}}

\def\cT{\mc{T}}
\def\E{\mathbb E}
\def\R{\mathbb R}
\def\C{\mathbb C}
\def\Z{\mathbb Z}
\def\Q{\mathbb Q}
\def\P{\mathbb P}
\def\H{\mathbb H}
\def\T{\mathbb T}

\def\N{\mathbb N}
\def\D{\mathbb D}
\def\L{\mathbb L}
\def\K{\mathbb K}
\def\CK{\mathbb{C}\mathbb{K}}
\def\V{\mathbb V}
\def\mfi{\mf{i}}
\def\mfj{\mf{j}}

\def\mfd{\mf{d}}
\def\mfs{\mf{s}}
\def\mft{\mf{t}}
\def\mfz{\mf{z}}

\def\la{\lambda}

\def\ep{\epsilon}
\def\De{\Delta}
\def\de{\delta}
\def\om{\omega}
\def\Om{\Omega}
\def\si{\sigma}
\def\Si{\Sigma}
\def\al{\alpha}

\def\ze{\zeta}
\def\na{\nabla}
\def\eset{\emptyset}

\makeindex
\begin{document}
\title{Deformation theory of log pseudo-holomorphic curves\\
and logarithmic Ruan-Tian perturbations}
\author{Mohammad Farajzadeh-Tehrani}
\date{\today}
\maketitle

\begin{abstract}
In a previous paper \cite{FT1}, for any logarithmic symplectic pair $(X,D)$ of a symplectic manifold $X$ and a simple normal crossings symplectic divisor $D$, we introduced the notion of log pseudo-holomorphic curve and proved a compactness theorem for the moduli spaces of stable log curves. 
In this paper, we introduce a natural set up for studying the deformation theory of log (and relative) curves. 
As a result, we obtain a logarithmic analog of the space of Ruan-Tian perturbations for these moduli spaces. 
For a generic compatible pair of an almost complex structure and a log perturbation term, we prove that the subspace of simple maps in each stratum is cut transversely. Such perturbations enable a geometric construction of Gromov-Witten type invariants for certain semi-positive pairs $(X,D)$ in arbitrary genera. In future works, we will use local perturbations and a gluing theorem to construct log Gromov-Witten  invariants of arbitrary such pair~$(X,D)$.
\end{abstract}
\begin{small}
\tableofcontents
\end{small}

%------------------------------------------------------------------------------------------------------
%---------------------Introuduction
%------------------------------------------------------------------------------------------------------
\section{Introduction}\label{Intro_s}
Studying pairs $(X,D)$ of a smooth complex projective variety $X$ and a  normal crossings (or \textbf{NC}) divisor $D$ has a long history in complex algebraic geometry.
In the symplectic category, M. McLean, A. Zinger, and the author recently introduced topological notations of NC symplectic divisor and configuration in arbitrary dimension; see \cite{FMZ1,FMZ3}.
We also constructed a class of almost K\"ahler structures that is suitable for defining and studying $J$-holomorphic curves relative to an NC symplectic divisor; see Definitions~\ref{SCD_dfn} and \ref{stdJ_dfn}.  
The theory of $J$-holomorphic curves has been a key tool in the study of symplectic manifolds ever since its inception by Gromov in the 1980s. 
Given a closed symplectic manifold $(X,\om)$, an $\om$-compatible or tame almost complex structure $J$ on $X$, $g,k\!\in\!\N$, and $A\!\in\!H_2(X,\Z)$ in the stable range (either $A\!\neq\!0$ or $2g\!+\!k\!\geq\! 3$), the objects of interest are the moduli spaces $\cM_{g,k}(X,A)$ of equivalence classes of genus $g$ degree $A$ $k$-marked  $J$-holomorphic maps 
into~$X$. 
Similarly, in the presence of a Lagrangian $L\!\subset\!X$, the objects of interests are the moduli spaces of marked bordered $J$-holomorphic curves with boundary on~$L$. 
In the presence of an NC symplectic divisor $D$, we study moduli spaces of $J$-holomorphic curves that intersect $D$ at finitely many points with pre-determined tangency orders.  
The particular choice of the symplectic structure $\om$ is often not important and might be hidden. 
In particular, we say $D\!\subset\!X$ is an SNC symplectic divisor whenever there is a symplectic form $\om$ on $X$ with respect to which $D$ is an SNC symplectic divisor.
Construction of such moduli spaces and the related enumerative/algebraic invariants involves:
(i) compactifying moduli spaces of such curves, (ii) putting some sort of oriented smooth structure on the moduli space (transversality and orientation problems, and the gluing analysis), and (iii) calculating/analyzing/using the resulting invariants or algebraic structures.\\

\noindent 
In a previous paper \cite{FT1}, we introduced a natural and explicit way of compactifying moduli spaces of pseudo-holomorphic curves (called log compactification) relative to a simple normal crossings (or \textbf{SNC}) symplectic divisor. 
Even if $D$ is smooth, our log compactification is somewhat different and smaller than the well-known relative compactification in \cite{IP1,LR,Li}. As the naming suggest, we expect it to be closely related to the log compactification of Gross-Siebert \cite{GS} and Abramovich-Chen \cite{AC}  in the algebraic setting.
In this paper, we introduce a natural set up for studying the deformation theory of log curves based on the logarithmic tangent bundle $TX(-\log D)$ associated with any symplectic \textbf{logarithmic pair} $(X,D)$, and address the transversality problem to some extent. 

%---------------------------------------------------------
\subsection{Log pseudo-holomorphic curves}
\noindent
Let 
$$
[N]\!\equiv\!\{1,\ldots,N\} \qquad \forall~N\!\in\!\N.
$$
Let $D\!\equiv\!\bigcup_{i\in [N]}D_i \!\subset\!X$ be an SNC symplectic divisor and $J$ be an $\om$-tame almost complex structure on $X$ such that 
$$
J(TD_i)\!=\!TD_i\qquad \forall i\!\in\![N].
$$ 
For every $J$-holomorphic map $u\colon\!\Si\!\lra\!X$ (with smooth domain) representing the homology class $A\!\in\!H_2(X,\Z)$, either $\tn{Im}(u)\!\subset\!D$ or there is a finite set of points $\{z_1,\ldots,z_k\}\!\subset\!\Si$ and a corresponding set $\mfs$ of vectors with non-negative integer coefficients
\bEq{stuple_e}
\mfs=\big(s_a=(s_{ai})_{i\in [N]}\big)_{a\in [k]}\!\in\!(\N^N)^k, \quad\tn{satisfying}\quad \sum_{a\in [k]} s_{ai}=A\cdot D_i\qquad \forall~i\!\in\![N],
\eEq
such that $u^{-1}(D)\!\subset\!\{z_1,\ldots,z_k\}\!\subset\!\Si$ and $u$ has a tangency of order $s_{ai}$ with $D_i$ at $z_a$ and nowhere else. Thus, $\mfs$ classifies the intersection type of $k$-marked $J$-holomorphic maps that are not mapped into $D$.  Let 
$$
\cM_{g,\mfs}(X,D,A)\subset \cM_{g,k}(X,A)
$$ 
denote the subspace\footnote{A marked point $z_a$ with $s_a\!=\!0\in \Z^N$ is a classical marked point. If $D$ is smooth, we can arrange the points $z_1,\ldots,z_k$ such that
$$
\mfs=(s_1,\ldots,s_{k_1},0,\ldots,0)\in \N^k, \qquad s_i\neq 0 \quad \forall~i\!\in\! [k_1].
$$ 
Then, the common practice in the relative theory is to denote $\cM_{g,\mfs}(X,D,A)$
by 
$$
\cM_{g,k_1,\ell}(X,D,A)\subset  \cM_{g,k}(X,A), \quad \ell=k-k_1.
$$
Separating marked points into different types is notationally cumbersome and not useful.} of all such $J$-holomorphic curves of type $\mfs$. Constructing a nice compactification\footnote{i.e., a naturally defined compact space that contains $\cM_{g,\mfs}(X,D,A)$ as a subset.} of $\cM_{g,\mfs}(X,D,A)$ and the related enumerative/algebraic invariants in the sense of (i)-(iii) above has been a challenging question that people have been working on from various perspectives for at least two decades. Technically speaking, the main goal is:\vskip.1in

\noindent
($\star$)~\textit{to construct a natural geometric compactification $\ov\cM_{g,\mfs}(X,D,A)$ of $\cM_{g,\mfs}(X,D,A)$ so that the definition of the contact vector $\mfs$ naturally extends to every element of $\ov\cM_{g,\mfs}(X,D,A)$, and $\ov\cM_{g,\mfs}(X,D,A)$ is (virtually) smooth enough to admit a natural class of cobordant  Kuranishi structures of the expected real dimension}
\bEq{log-dim_e}
2\Big(c_1^{TX}(A)+(n-3)(1-g)+k-A\cdot D\Big).
\eEq
\vskip.1in
\noindent
We refer to \cite{FF,MW} for the technical terms in ($\star$). For each $a\!\in\![k]$, let 
\bEq{Ia_e}
I_a\!=\!\{i\!\in\![N]~~~\tn{s.t}~~~s_{ai}>0\}\subset [N]
\eEq
be the subset of indices where $u$ intersects $D_i$ at the $a$-th marked point $z_a$ non-trivially. Let 
$$
D_I=\bigcap_{i\in I} D_i \qquad \forall~I\!\subset\![N].
$$
In ($\star$), we furthermore expect the natural stabilization  and evaluation  maps 
\bEq{ev-st_e}
\tn{st}\colon\!\ov\cM_{g,\mfs}(X,D,A)\!\lra\!\ov\cM_{g,k} \quad \tn{and}\quad
\tn{ev}\equiv(\tn{ev}_a)_{a\in [k]} \colon\!\ov\cM_{g,\mfs}(X,D,A)\!\lra\!X^{\mfs}\!\equiv\!\prod_{a\in [k]} D_{I_a}
\eEq
to be continuous (or smooth). 
Similarly to the classical case, if $\ov\cM_{g,\mfs}(X,D,A)$ has a ``nice" orbifold structure of the expected real dimension (\ref{log-dim_e}), Gromov-Witten (or \textbf{GW}) invariants of the logarithmic pair $(X,D)$ with primary insertions can be  defined by intersecting the image of $\ov\cM_{g,\mfs}(X,D,A)$ under $\st\!\times\!\ev$ with appropriate cycles and counting (with sign and $\Q$-weights)  the number of intersection points. With $J$ carefully chosen from an appropriate space of almost K\"ahler structures, the resulting rational numbers will be independent of $J$, and only depend on the genus $g$, homology class $A$, tangency data $\mfs$, cohomology classes of the aforementioned cycles in $\ov\cM_{g,k}$ and $X^{\mfs}$, and the deformation equivalence class $[X,D,\om]$ of $(X,D,\om)$.\\

\noindent
The case where $D$ is smooth was treated at the turn of the millennium. The so called ``relative'' theory of \cite{Li} in algebraic geometry, and \cite{IP1,LR} in the symplectic category\footnote{The construction is not complete in any of these two papers; see \cite{FZ1}.} address ($\star$) when $D$ is smooth.  More recently, Gross-Siebert and Abramovich-Chen \cite{AC,Ch,GS} used sophisticated  techniques from logarithmic complex algebraic geometry to fully address ($\star$) in the algebraic case. In the algebraic construction, $\ov\cM_{g,\mfs}^{\log}(X,D,A)$ is required to have a richer structure. Roughly speaking, it should be a coarse moduli space for the functor  which assigns to a log scheme $B$ the set of all families of ``good" log curves in $X$ with base $B$.  Ultimately, completing all the steps needed for ($\star$) will allow the formulation of a symplectic analogue of logarithmic GW invariants in \cite{AC,GS} as well-defined invariants of the deformation equivalence class of logarithmic symplectic pairs $(X,D)$. In that regard, there has been a series of papers by B. Parker (see \cite{BP1,BP2} and the references therein) and a paper by Ionel \cite{I} that aim to address ($\star$) in the symplectic category; see  \cite{FT2} for some comments on these works. Our approach in \cite{FT1} and here is substantially different from these approaches. We avoid changing the target (as in \cite{I}) or putting extra sheaves on it (as in \cite{BP1,BP2}).\\

\noindent
In \cite{FT1}, we introduced a geometric notion of \textbf{log $J$-holomorphic curve} relative to an arbitrary SNC symplectic divisor and proved the following compactness result for the moduli space $\ov\cM^{\log}_{g,\mfs}(X,D,A)$ of  stable log $J$-holomorphic curves.

\begin{theorem}[{\cite[Thm~1.1]{FT1}}]
Let $(X,\om)$ be a closed symplectic manifold and $D\!\subset\!X$ be an SNC symplectic divisor. For  ``suitable'' choice of $J$, the Gromov sequential convergence topology on $\ov\cM_{g,k}(X,A)$ lifts to a compact metrizable sequential convergence topology
on $\ov\cM^{\log}_{g,\mfs}(X,D,A)$ such that the forgetful map 
\bEq{MiotaMap_e}
\iota \colon \ov\cM^{\log}_{g,\mfs}(X,D,A)\lra \ov\cM_{g,k}(X,A)
\eEq
is a continuous local embedding. If $g\!=\!0$, then $\iota$ is globally an  embedding.
\eTh
\noindent
Intuitively, if $g\!>\!0$, the map $\iota$ behaves like an immersion.
If $D$ is smooth, there exists a surjective map from the relative moduli space $\ov\cM^{\tn{rel}}_{g,\mfs}(X,D,A)$ in \cite{Li,LR,IP1} to the log moduli space $\ov\cM^{\log}_{g,\mfs}(X,D,A)$; see \cite[Sec~3.5]{FT1}.\\

\noindent 
For each $i\!\in\![N]$, let $\cN_XD_i$ denote the normal bundle of $D_i$ in $X$. An element of $\ov\cM^{\log}_{g,\mfs}(X,D,A)$ is the equivalence\footnote{Two tuples $f_{\log}\!=\!\big((u_v,\Si_v,\mfj_v,\ze_v\!=\!(\ze_{v,i})_{i\in I_v}\big)_{v\in \V},\vec{z}\big)$ and $f'_{\log}\!=\!\big((u'_v,\Si'_v,\mfj'_v,,\ze'_v=(\ze'_{v,i})_{i\in I_v}\big)_{v\in \V},\vec{z}\,'\big)$ are equivalent if there exists a holomorphic reparametrization $h \colon\! \Si\!\lra\!\Si'$ such that $u'\circ h\!=\!u$, $h(\vec{z})=\vec{z}\,'$, and $h_{v}^*\ze'_{h(v),i}\!=\!c_{v,i}\ze_{v,i}$ for all $v\!\in\!\V$ and some $c_{v,i}\!\in\!\C^*$. Here $v\!\lra\!h(v)$ is the induced map on the vertices of the dual graph.} class of a tuple 
$$
f_{\log}=\big(u_v,\Si_v,\mfj_v,\vec{z}_v,\ze_v=(\ze_{v,i})_{i\in I_v}\big)_{v\in \V}
$$
such that 
\bEn
\item for each $v\!\in\!\V$, $I_v\!\subset\![N]$ is the maximal subset where the image of $u_v$ is contained in $D_{I_v}$; 
\item for each $v\!\in\!\V$ and $i\!\in\!I_v$, $\ze_{v,i}$ is a meromorphic section of $u_v^*\cN_XD_i$; 
\item forgetting $\ze_v$, $f=(u_v,\Si_v,\mfj_v,\vec{z}_v\big)_{v\in \V}$ defines a classical nodal map in $\ov\cM_{g,k}(X,A)$; 
\item and, the conditions listed below are satisfied. 
\eEn
First, for each $x\!\in\!\Si_v$, the pair $(u_v,\ze_v)$ gives rise to a well-defined \textbf{tangency order} vector in $\Z^N$ with $D$ at $x$, which we denote by $\ord_x(u_v,\ze_v)$. 
Then, the conditions in (4) are: 
\bEnalph
\item\label{ord_l} every point in $\Si$ with a non-trivial tangency order vector is either a marked point or a nodal point, 
\item the tangency order vector at the marked point $z_a$ is the pre-determined vector $s_a\!\in\!\Z^N$ in~$\mfs$;
\item\label{CNode_l} the tangency order vectors at the nodes are the opposite of each other;
\item\label{CT_l} there exists a vector-valued function $s\colon\!\V\!\lra\!\R^N$ such that $s_v\!=\!s(v)\!\in\! \R_{+}^{I_v}\!\times \{0\}^{[N]-I_v}$ for all $v\!\in\!\V$, and $s_{v'}\!-\!s_{v}$ is a positive multiple of the contact order vector of any nodal point on $\Si_v$ connected to $\Si_{v'}$, for all $v,v'\!\in\!\V$;
\item\label{GObs_it} and, certain group element $\tn{ob}_{\Gamma}(f_{\log})\!\in\!\mc{G}(\Gamma)$  associated to $f_{\log}$ is equal to $1$.
\eEnalph
Conditions~\ref{ord_l}-\ref{CT_l} are combinatorial conditions on the dual graph of $f_{\tn{log}}$. A tuple $f_{\tn{log}}$ satisfying all the conditions except the last one\footnote{Or sometimes, the last two conditions.} is called a \textbf{pre-log curve}.
The last two conditions ensure that each boundary stratum of  $\ov\cM^{\log}_{g,\mfs}(X,D,A)$ has positive  expected complex co-dimension, and every nodal log curve is virtually smoothable, respectively. The moduli space $\ov\cM^{\log}_{g,\mfs}(X,D,A)$ can be described without mentioning the meromorphic sections $\ze_{v,i}$. Whenever such a section exists, it is unique up to multiplication by a non-zero constant; see Remark~\ref{Noze_e}. In (\ref{MiotaMap_e}), different lifts of a stable map in $\ov\cM_{g,k}(X,A)$ to $\ov\cM_{g,\mfs}^{\log}(X,D,A)$  are characterized by different choices of tangency order vectors at the nodal points, satisfying Condition~\ref{CNode_l} and (\ref{stuple_e}) on each smooth component.
See Section~\ref{LogJnu_ss} for the details. 

%----------------
\subsection{Deformation theory}
\textit{Our first goal in this paper is to introduce a natural set up for studying the deformation theory of log pseudo-holomorphic curves}. This setup can also be used to address the transversality issue in other applications of moduli spaces of curves relative to a smooth or SNC divisor, such as in the construction of the relative Fukaya category \cite{DF}. 
In the latter, it would also be useful for properly addressing the orientation problem. The key to this setup is the observation that the deformation theory of log $J$-holomorphic curves relies on the linearization of CR equation (\ref{Jnu-holo_e}) as an operator acting on the set of sections of the \textbf{log tangent bundle} $TX(-\log D)$ instead of $TX$. In the holomorphic case, the log tangent sheaf is the sheaf of holomorphic tangent vector fields in $TX$ whose restriction to each $D_i$ is tangent to $D_i$. The construction in the symplectic case is similar but depends on some auxiliary data. The deformation equivalence class of the complex vector bundle $TX(-\log D)$ only depends on the deformation equivalence class of $(X,D,\om)$. Furthermore, we have 
$$
c_1^{TX(-\log D)}= c_1^{TX} - \sum_{i\in [N]} \tn{PD}(D_i) \in H^2(X,\Z);
$$
see Section~\ref{LogBundle_ss} or \cite{FMZ2,FMZ4}.
Therefore, in analogy with the classical dimension formula 
\bEq{exp-dim_e}
d=\tn{exp-dim}_{\R}~\ov\cM_{g,k}(X,A)=2\Big(c_1^{TX}(A)+(n-3)(1-g)+k\Big),
\eEq
the expected dimension in (\ref{log-dim_e}) can be re-written as  
\bEq{dlog_e}
d^{\log}=\tn{exp-dim}_{\R}~\ov\cM^{\log}_{g,\mfs}(X,D,A)= 2\Big(c_1^{TX(-\log D)}(A)+(n-3)(1-g)+k\Big).
\eEq

\noindent
The analytical setup in the classical case follows the following steps. Given a smooth domain $(\Si,\mfj)$, let $\tn{Map}_{A}(\Si,X)$ denote the space of all smooth maps $u\colon\!\Si\!\lra\!X$ that represent the homology class~$A$. Let 
$$
\cE_A\lra\tn{Map}_{A}(\Si,X)
$$
be the infinite dimensional bundle whose fiber over $u$ is $\Gamma(\Si, \Om^{0,1}_{\Si,\mfj}\otimes_{\C} u^*TX)$. The Cauchy-Riemann (or \textbf{CR}) equation 
\bEq{Jnu-holo_e}
\dbar u\equiv \frac{1}{2}\big( \nd u+ J \nd u \circ \mfj)
\eEq 
can be seen as a section of this infinite dimensional bundle. To be precise, we need to consider a Sobolev completion of these spaces for the Implicit Function Theorem to apply, but, by elliptic regularity, every solution of $\dbar u \!=\!0$ will be smooth. The linearization of the $\dbar$-section at any $J$-holomorphic map $u$ is an $\R$-linear map
$$
\tn{D}_u \dbar\colon \Gamma(\Si,u^*TX)\lra \Gamma(\Si,\Om^{0,1}_{\Si,\mfj}\otimes_\C u^*TX)
$$ 
that is the sum of a $\C$-linear $\dbar$-operator and a compact operator. Therefore, it is a Fredholm operator and Riemann-Roch applies; i.e., it has finite dimensional kernel and co-kernel, and
\bEq{RR_e}
\dim_\R \Def(u)\!-\!\dim_\R \Obs(u)\!=\!2\big( \tn{deg}(u^*TX) \!+\! \dim_\C\! X (1\!-\!g)\big),
\eEq
where $\Def(u)\!=\!\ker(\tn{D}_u\dbar)$ and $ \Obs(u)\!=\!\coker(\tn{D}_u\dbar)$.
The first space corresponds to infinitesimal deformations of $u$ (over the fixed smooth marked domain) and the second one is the obstruction space for integrating elements of $\Def(u)$ to actual deformations. 
If $\Obs(u)\!\equiv\! 0$, by Implicit Function Theorem \cite[Thm~A.3.3]{MS2}, in a small neighborhood $B_\ep(u)$ of $u$ in $\tn{Map}_A(\Si,X)$ the set of $J$-holomorphic maps $V_u \equiv \dbar^{-1}(0)\cap B_\ep(u)
$ is a smooth manifold of real dimension (\ref{RR_e}), all the elements of $\Def(u)$ are smooth, and $T_u V_u\cong \Def(u)$;
see \cite[Thm~3.1.5]{MS2}. The manifold $V_u$ carries a natural orientation.\vskip.1in

\noindent
In the logarithmic case, given $\big(\Si,\mfj,\vec{z}\!=\!(z_1,\ldots,z_k)\big)$, $A$, and $\mfs$, we generalize this construction in a natural way. In Section~\ref{Generic_ss}, we construct a configuration space 
$$
\tn{Map}_{A,\mfs}((\Si,\vec{z}),(X,D))\subset \tn{Map}_{A}(\Si,X)
$$ 
whose elements are smooth maps that have contact type $\mfs$ with $D$ at $\vec{z}$ in a suitable sense. Let
$$
\cE_{A,\mfs}\lra\tn{Map}_{A,\mfs}\big((\Si,\vec{z}),(X,D)\big)
$$
be the infinite dimensional bundle whose fiber over $u$ is $\Gamma\big(\Si, \Om^{0,1}_{\Si,\mfj}\otimes_{\C} u^*TX(-\log D)\big)$.
The section $\dbar$ of $\cE_A$ restricts to a section $\dbar^{\log}$ of $\cE_{A,\mfs}$. Recall that there is a $\C$-linear homomorphism 
\bEq{IotaMap_e}
\iota\colon TX(-\log D)\lra TX
\eEq 
(covering $\tn{id}_X$) that is an isomorphism away from $D$. This homomorphism induces $\C$-linear maps 
$$
\aligned
\iota_1&\colon \Gamma(\Si,u^*TX(-\log D))\lra  \Gamma(\Si,u^*TX),\\
 \iota_2&\colon \Gamma(\Si,\Om^{0,1}_{\Si,\mfj}\otimes_\C u^*TX(-\log D))\lra  \Gamma(\Si,\Om^{0,1}_{\Si,\mfj}\otimes_\C u^*TX).
\endaligned
$$   
The linearization $\tn{D}_u^{\log}\dbar$ of $\dbar^{\log}$ is a Fredholm\footnote{After taking Sobolev completions of these spaces.} linear map
$$
\tn{D}^{\log}_u\dbar\colon\Gamma(\Si,u^*TX(-\log D))\lra \Gamma(\Si,\Om^{0,1}_{\Si,\mfj}\otimes_\C u^*TX(-\log D))
$$
such that $\iota_2\circ \tn{D}^{\log}_u\dbar = \tn{D}_u\dbar \circ \iota_1$ (on the subset of smooth sections). Furthermore, by Riemann-Roch and Implicit-Function Theorem, if $\tn{coker}(\tn{D}^{\log}_u\dbar)\!=\! 0$, the set of $J$-holomorphic maps of contact type $\mfs$ close to $u$ form an oriented smooth manifold of real dimension
$$
2\big( \tn{deg} (u^*TX(-\log D))\!+\! \dim_\C\! X (1\!-\!g)\big).
$$
Considering the deformations of the marked domain $(\Si,\mfj,\vec{z})$, we get the dimension formula (\ref{dlog_e}) and the deformation-obstruction long exact sequence 
$$
0\lra ~\tn{aut}(C) \stackrel{\de}{\lra} 
\Def_{\log}(u)\lra \Def_{\log}(f) \lra~\Def(C) \stackrel{\de}{\lra} 
\Obs_{\log}(u)  \lra  ~\Obs_{\log}(f) \lra ~0,
$$
where $f\!=\!\big(u,C\!=\!(\Si,\mfj,\vec{z})\big)$, 
$$
\aligned
\tn{aut}(C)\!=\!H^0_{\dbar}(T\Si(-\log\vec{z}))&,\qquad \Def(C)\!\cong\! H^1_{\dbar}(T\Si(-\log \vec{z})),\\
\Def_{\log}(u)\!=\!\ker(\tn{D}^{\log}_u\dbar)&, \qquad \Obs_{\log}(u)\!=\!\coker(\tn{D}^{\log}_u\dbar).
\endaligned
$$
If $\tn{Obs}_{\log}(f)\!=\!0$, then a small neighborhood $B(f)$ of $f$ in $\cM_{g,\mfs}(X,D,A)$ is a smooth orbifold of the expected dimension (\ref{dlog_e}). In Sections~\ref{NGSmooth_ss} and \ref{NodalMap_ss}, we will extend this setup to log maps with smooth domain and image in a stratum $D_I$, and to general nodal log maps, respectively. 

\subsection{Transversality}\label{Transversality_ss}
\noindent
For the general construction of Gromov-Witten type invariants of every arbitrary pair $(X,D)$, we need a gluing theorem that generalizes the known gluing theorem in the classical case. 
We also need to generalize the theory of Kuranishi structures to allow toric singularities. In the classical case, by restricting to the subset of simple maps, transversality can be achieved by global perturbations of the $\dbar$-equation and $J$. Then, in the case of \textit{semi-positive} symplectic manifolds, as worked out by McDuff-Salamon \cite{MS1} in genus $0$ and Ruan-Tian \cite{RT} in arbitrary genera, the classical analogue of the map $\tn{st}\times\tn{ev}$ in  (\ref{ev-st_e}) over each stratum of non-simple maps factors through a positive complex co-dimension space of simple maps, and thus can be ignored. Therefore, in semi-positive situations, gluing and virtual techniques are not needed, and GW invariants can be defined by a direct count of perturbed pseudo-holomorphic curves in the following sense. 

\bDf{positive_dfn}
A closed real $2n$-dimensional symplectic manifold $(X,\om)$  is called \textbf{semi-positive} if 
$$
c_1^{TX}(A)\!\geq\! 3-n ~~\Rightarrow~~ c_1^{TX}(A)\!\geq\! 0
$$ 
for all $A\!\in\! \pi_2(X)$ such that $\om(A)\!>\!0$. It is called \textbf{positive}\footnote{Other terms such as {monotone}, {strongly semi-positive}, {convex}, etc. have also been used in the literature.} if 
$$
c_1^{TX}(A)\!\geq\! 3-n ~~\Rightarrow~~ c_1^{TX}(A)\!>\! 0
$$  
for all $A\!\in\! \pi_2(X)$ such that $\om(A)\!>\!0$. We say $[X,\om]$ is semi-positive/positive if $(X,\om')$ for some $\om'$ deformation equivalent to $\om$ is semi-positive/positive.
\eDf

\noindent
In particular, every symplectic manifold of real dimension $6$ or less is semi-positive and every symplectic 4-manifold is positive.
Let $\cM^\star_{g,k}(X,A)\subset \cM_{g,k}(X,A)$ denote the subspace of \textbf{simple} (not multiple-cover) maps. 
By \cite[Thm~3.1.5]{MS2}, for generic $J$, $\cM^\star_{g,k}(X,A)$ is a naturally oriented smooth manifold of the expected real dimension (\ref{exp-dim_e}).
If $g\!=\!0$ and $(X,\om)$ is positive, for generic $J$, a deliberate dimension counting argument shows that the image in $\ov\cM_{0,k}\!\times\! X^k$ of the complement 
$$
\ov\cM_{0,k}(X,A)\setminus \cM_{0,k}^\star(X,A)
$$
under
$$
\tn{st}\times\tn{ev}\colon\!\ov\cM_{0,k}(X,A,\nu)\!\lra\!\ov\cM_{0,k}\times X^k
$$
is a set of real codimension at least $2$ and can be ignored. In other words, the inclusion $\cM_{0,k}^\star(X,A)\!\subset\! \ov\cM_{0,k}(X,A)$ gives rise to a \textbf{pseudocycle} whose homology class
$$
\GW^{X}_{0,k,A}\subset H_d( \ov\cM_{0,k}\times X^k,\Z)
$$
is independent of the choice of  $J$ and the particular choice of $\om$ in its deformation equivalence class; see \cite[Thm 6.6.1]{MS2} for more details.\\

\noindent
In the higher genus case, the same argument does not work and one needs to perturb the Cauchy-Riemann equation to take care of constant and multiple-cover maps.
\noindent
Given a smooth Riemann surface $\Si$ with complex structure $\mfj$ and a (sufficiently) smooth map $u\colon\!\Si\!\lra\! X$, the space of \textbf{perturbation terms} for the pair $(u,\Si)$ is the infinite dimensional vector space 
$$
\Gamma(\Sigma,\Om_{\Si,\mfj}^{0,1}\otimes_\C u^*TX)
$$
of smooth $u^*TX$-valued $(0,1)$-forms with respect to $\mfj$ on $\Si$ and $J$ on $TX$. 
Given a perturbation term $\nu$, we say $u$ is $(J,\nu)$-holomorphic if it satisfies the perturbed CR equation $\dbar u\!=\!\nu$.
If $2g\!+\!k\!\geq\!3$, Ruan-Tian defined a class of global perturbation terms $\nu$, where each $\nu$ is a section in 
\bEq{RTNuSpace_e}
\Gamma\big(\ov{\mf{U}}_{g,k}\times X, \pi_1^*\Om^{0,1}_{g,k}\otimes_\C \pi_2^*TX\big).
\eEq
Here, $\ov{\mf{U}}_{g,k}$ is a universal family over a ``regular"  covering 
$\ov{\mf{M}}_{g,k}$ of $\ov\cM_{g,k}$. For each $\nu$, the moduli space of interest is the set
$$
\ov\cM_{g,k}(X,A,J,\nu)
$$ 
of the equivalence classes of stable degree $A$ genus $g$ $k$-marked nodal $(J,\nu)$-holomorphic maps $(\phi,u,C)$. Here $\phi\colon\! C\!\lra\! \ov{\mf{U}}_{g,k}$ is a holomorphic map\footnote{It inductively contracts bubble components with unstable domain.} from a genus $g$ $k$-marked nodal curve $C$ onto a fiber of the universal family $\ov{\mf{U}}_{g,k}$ and $\dbar u =(\phi,u)^*\nu$; see Section~\ref{LogJnuMaps_s}. By \cite[Crl 3.9]{RT}, the space $\ov\cM_{g,k}(X,A,J,\nu)$ is Hausdorff and compact with respect to a similarly defined Gromov convergence topology. By \cite[Thm 3.16]{RT}, for each $\om$-tame $J$ and generic  perturbation term $\nu$, the main stratum $\cM_{g,k}(X,A,J,\nu)$ consisting of maps with smooth domain\footnote{Since $2g\!+\!k\!\geq\!3$, every such map is automatically simple according to Definition~\ref{Simple_dfn0}.} is cut out transversely by the $\{\dbar-\nu\}$-section and thus it is  a smooth manifold of the expected dimension (\ref{exp-dim_e}). 
By \cite[Thm 3.11]{RT}, it has a canonical orientation. Furthermore, if $(X,\om)$ is semi-positive, by \cite[Prp 3.21]{RT}, for generic $(J,\nu)$, the image of  the complement of $\cM_{g,k}(X,A,J,\nu)$ in $\ov\cM_{g,k}(X,A,J,\nu)$ under $\st\times \ev$
is contained in images of maps from smooth even-dimensional manifolds of at least $2$ real dimension less than the main stratum. Thus, similarly to the positive case and after dividing by the degree of the regular covering used to define $\nu$,
the inclusion 
$$
\cM_{g,k}(X,A,J,\nu)\!\subset\! \ov\cM_{g,k}(X,A,J,\nu)
$$ 
gives rise to a GW homology class
$$
\GW^{X}_{g,k,A}\subset H_d( \ov\cM_{g,k}\times X^k,\Q)
$$
independent of the choice of the admissible almost complex structure $J$, the perturbation $\nu$, or the particular choice of $\om$ in its deformation equivalence class. 
We conclude that, in the semi-positive situations, the resulting GW invariants are \textit{enumerative} in the sense that they can be interpreted as a finite $\Q$-weighted count of $(J,\nu)$-holomorphic maps of fixed degree and genus meeting some prescribed cycles at the marked points.\\

\noindent
\textit{Our second goal is to introduce a logarithmic analogue of Ruan-Tian perturbations of the Cauchy-Riemann section for log moduli spaces, and use them to 
achieve transversality over the subspace of simple maps.}\vskip.1in

\noindent
The definition of the obstruction bundle $\cE_{A,\mfs}$ above and (\ref{RTNuSpace_e}) suggest that a logarithmic perturbation term should be an element of 
$$
\Gamma\big(\ov{\mf{U}}_{g,k}\times X, \pi_1^*\Om^{0,1}_{g,k}\otimes_\C \pi_2^*TX(-\log D)\big).
$$
For each $\nu_{\log}$ in this space, with $\iota$ as in (\ref{IotaMap_e}), $\nu=\iota(\nu_{\log})$ is the associated classical perturbation term in (\ref{RTNuSpace_e}).
The construction of log $(J,\nu)$-holomorphic moduli spaces works for arbitrary $\nu\!=\!\iota(\nu_{\log})$. However, for simplicity, we restrict to a subclass of such sections $\mc{H}_{g,k}(X,D)$, described below, whose elements have a standard form near each stratum $D_I$ of $D$.  With the perturbation term $\nu_{\log}$ as above and a similarly defined log linearization map $\tn{D}^{\log}_u\{\dbar-\nu\}$ in place of $\tn{D}_u\{\dbar-\nu\}$,  the proof of the main results follow the same general steps as in \cite{RT}. However, we face some difficulties in dealing with non-simple nodal maps. \\

\noindent
We consider a particular set of almost K\"ahler auxiliary data $\tn{AK}(X,D,\om)$ defined in \cite{FMZ1} and compatible logarithmic perturbations in the following sense. Since \cite[Thm~1.1]{FT1} also includes integrable $J$, similar results for perturbations compatible with integrable almost complex structures can be obtained. An element $(\cR,J)$ in $\tn{AK}(X,D,\om)$  consists of 
\bIt
\item a regularization $\cR$, that is a compatible set of symplectic identifications of neighborhoods of $\{D_I\}_{I\subset [N]}$ in their normal bundles with their neighborhoods in $X$ in the sense of Definition~\ref{omRegV_dfn}; 
\item and an $\om$-tame almost complex structure $J$ compatible with $\cR$ in the sense of Definition~\ref{stdJ_dfn}.
\eIt
The space $\tn{AK}(X,D,\om)$ might be empty for some choices of $\om$. Let\footnote{This space is denoted by $\tn{Symp}^+(X,D)$ in \cite{FMZ1,FMZ3}.}  $\tn{Symp}(X,D)$ be the space of all symplectic structures $\om$ on $X$ with respect to which $D$ is an SNC symplectic divisor. Let $\tn{AK}(X,D)$ be the space of tuples $(\om,\cR,J)$ where $\om\in \tn{Symp}(X,D)$ and $(\cR,J)\!\in\!\tn{AK}(X,D,\om)$. As a consequence of \cite[Thm~2.13]{FMZ1}, the projection map 
$$
\tn{AK}(X,D)\lra \tn{Symp}(X,D), \qquad (\om,\cR,J)\lra \om,
$$
is a weak homotopy equivalence. In particular, starting from any $\om$, we can deform it without changing its cohomology class to another $\om'$ such that $\tn{AK}(X,D,\om')\neq \eset$. Therefore, the subspace $\tn{AK}(X,D)_{[\om]}$ of tuples $(\om',\cR,J)$ such that $\om'$ is deformation equivalent to $\om$ in $\tn{Symp}(X,D)$ is path connected.  \\

\noindent
For any fixed tuple $(\om,\cR,J)$ and $g,k\!\in\!\N$, with $2g+k\!\geq\! 3$, we define a class of ``$(\cR,J)$-compatible" perturbation terms $\nu$ in Definition~\ref{stdJnu_dfn}. The compatibility condition requires $\nu$ to be of a standard form with respect to the regularization $\cR$ in the following sense. 
For every $I\!\subset\![N]$, the normal bundle $\cN_XD_I$ of $D_I$ in $X$ admits a decomposition into a direct sum of complex line bundles
\bEq{DeccNXDI_e}
\cN_XD_I\cong \bigoplus_{i\in I} \cN_{X}D_i|_{D_I}.
\eEq  
The regularization map $\Psi_I$ in $\cR$ gives a stratified identification of a neighborhood $D_I$ in $\cN_XD_I$ with a neighborhood of that in $X$ so that 
\bIt
\item $\Psi_I^*\om$ is of some standard form,
\item on the horizontal subspace $T^{\tn{hor}}\cN_XD_I\cong \pi_I^*TD_I$ defined via the connections in $\cR$, $\Psi_I^*J$ is the pull back of~$J|_{TD_I}$ via the projection map $\pi_I\colon\!\cN_XD_I\!\lra\!D_I$,
\item and in the vertical subspace $T^{\tn{ver}}\cN_XD_I\cong \pi_I^*\cN_XD_I$, $\Psi_I^*J$ is the direct sum complex structure on the right-hand side of (\ref{DeccNXDI_e}). 
\eIt
Then, similarly to the definition of $J$, the horizontal component of $\nu$ is required to be the pull back from $D_I$ of some $\nu_I$, and the vertical component of $\nu$ is required to be $(\C^*)^I$-equivariant with respect to component-wise action of $\C^*$ on the right-hand side of (\ref{DeccNXDI_e}). Let $\cH_{g,k}(X,D)$ be the set of such tuples $(\om,\cR,J,\nu)$. As a consequence of \cite[Thm~2.13]{FMZ1}, the projection map 
$$
\cH_{g,k}(X,D)\lra \tn{Symp}(X,D), \qquad (\om,\cR,J,\nu)\lra \om,
$$
is again a weak homotopy equivalence. This implies that any invariant of the deformation equivalence classes in $\cH_{g,k}(X,D)$ is an invariant of the symplectic deformation equivalence class of $(X,D,\om)$. The subspace $\cH_{g,k}(X,D)_{[\om]}$ of tuples $(\om',\cR,J,\nu)$ such that $\om'$ is deformation equivalent to $\om$ in $\tn{Symp}(X,D)$ is path connected.
Given $(\om,\cR,J,\nu)\!\in\!\cH_{g,k}(X,D)$, in Section~\ref{LogJnuMaps_s}, we construct the moduli space 
$$
\ov\cM^{\log}_{g,\mfs}(X,D,A,\nu)
$$ 
of equivalence classes of stable $k$-marked genus $g$ degree $A$ log $(J,\nu)$-holomorphic curves of contact type $\mfs$, similarly to \cite[Sec~3.2]{FT1}. \\

\noindent
The following result is the straightforward generalization of \cite[Thm~1.1]{FT1} to the $(J,\nu)$-holomorphic case. This time, the rescaling map normal to $D_i$ in the proof of \cite[Prp~4.10]{FT1} yields a meromorphic section with respect to the $\dbar$-operator $\tn{D}^{\cN_i}_u\{\dbar -\nu\} $ in (\ref{CRLNuNormal_ei}) instead of $u^*\dbar_{\cN_XD_i}=\tn{D}^{\cN_i}_u \dbar $ in \cite[Lmm~2.1]{FT1}.

\begin{theorem}\label{Compactness_th}
Suppose $X$ is closed, $D=\!\bigcup_{i\in [N]}D_i\!\subset\! X$ is an SNC symplectic divisor, $A\!\in\!H_2(X,\Z)$, $g,k\!\in\!\N$,  and $\mfs\!\in\!\big(\Z^N\big)^k$.  
For every $(\om,\cR,J,\nu)\!\in\!\cH_{g,k}(X,D)$, the Gromov sequential convergence topology on $\ov\cM_{g,k}(X,A,\nu)$ lifts to a compact metrizable sequential convergence topology
on $\ov\cM^{\log}_{g,\mfs}(X,D,A,\nu)$ such that the forgetful map 
\bEq{FogetLog_e}
 \ov\cM^{\log}_{g,\mfs}(X,D,A,\nu)\lra \ov\cM_{g,k}(X,A,\nu)
\eEq
is a continuous local embedding.
If $g\!=\!0$, then (\ref{FogetLog_e}) is a global embedding.
\eTh

\noindent
Furthermore, the forgetful and evaluation maps $\tn{st}$ and $\tn{ev}$ in (\ref{ev-st_e}) are continuous. Each moduli space $\ov\cM^{\log}_{g,\mfs}(X,D,A,\nu)$ is coarsely stratified by the subspaces
$$
\cM_{g,\mfs}(X,D,A,\nu)_{\Gamma}
$$
consisting of log $(J,\nu)$-holomorphic curves with the decorated dual graph $\Gamma$. The set of vertices $\V$ corresponds to the smooth components of $\Si$, the set of edges $\E$ corresponds to the nodes of $\Si$, and the set of roots $\L$ corresponds to the marked points. Each $v\!\in\!\V$ is decorated by the degree $A_v\!\in\!\H_2(X,\Z)$, genus $g_v\!\in\! \N$, and the index set $I_v\!\subset\![N]$.  An edge $e\!\in\!\E$ is decorated by the index set $I_e\!\subset\![N]$ if $u(q_{e})\!\in\!D_{I_e}\setminus \partial D_{I_e}$. Each edge can be oriented in two ways $\uvec{e}$ and $\ucev{e}$ such that the node $q_e$ is obtained by identifying two nodal points $q_{\uvec{e}}$ and $q_{\scz\ucev{e}}$. Each oriented edge $\uvec{e}$ is decorated by the contact order vector $s_{\uvec{e}}$ satisfying Condition~\ref{CNode_l} of Page~\pageref{CNode_l} and (\ref{stuple_e}) on each smooth component. The set of such $\Gamma$ is finite. Similarly to \cite[Sec 5]{FT1}, associated to every admissible decorated dual graph $\Gamma$ we get a $\Z$-linear map
$$
\vr=\vr_{\Gamma}\colon \D=\Z^\E \oplus \bigoplus_{v\in \V} \Z^{I_v}\lra \T=\bigoplus_{e\in \E} \Z^{I_e};
$$  
see  (\ref{DtoT_e}). Condition~\ref{CT_l} in Page~\pageref{CT_l} means that either $\Gamma$ is the trivial one-vertex graph ($\D=\T=0$, corresponding to the virtually main stratum $\cM_{g,\mfs}(X,D,A,\nu))$ or $\tn{ker}(\vr)$ has an element in the positive quadrant. The complex torus in Condition~\ref{GObs_it} has Lie algebra $\tn{coker}(\vr)\otimes \C$.

\bDf{Simple_dfn}
A $(J,\nu)$-holomorphic log map is called \textbf{simple} if the underlying $(J,\nu)$-holomorphic map is simple.
\eDf

\noindent
For each decorated dual graph $\Gamma$, let 
$$
\cM^{\star}_{g,\mfs}(X,D,A,\nu)_{\Gamma}\subset \cM_{g,\mfs}(X,D,A,\nu)_{\Gamma}
$$
denote the open subspace consisting of simple maps. 

\bTh{TransGamma_thm} 
Suppose $X$ is a closed symplectic manifold, $D=\bigcup_{i\in [N]} D_i$ is an SNC symplectic divisor,  $A\!\in\!H_2(X,\Z)$, $g,k\!\in\!\N$,  and $\mfs\!\in\!(\Z^N)^k$. 
For each admissible decorated dual graph $\Gamma$, the following statements hold.
\bEn
\item If $2g\!+\!k\!\geq\!3$, for any given choice of universal family in (\ref{UFamily_e}), there exists a Baire set of second category $\cH^{\Gamma}_{g,k}(X,D)\!\subset\!\cH_{g,k}(X,D)$ such that for each $(\om,\cR,J,\nu)\!\in\!\cH^{\Gamma}_{g,k}(X,D)$, the subspace of simple maps 
$$
\cM^\star_{g,\mfs}(X,D,A,\nu)_\Gamma\!\subset\!\cM_{g,\mfs}(X,D,A,\nu)_\Gamma
$$ 
is a naturally oriented smooth manifold of the real dimension 
\bEq{dGamma_e}
2\big(c_1^{TX(-\log D)}(A) + (n-3)(1-g) + k - \dim \tn{ker}(\vr)\big);
\eEq 
the restriction of $\tn{st}\times\tn{ev}$ in (\ref{ev-st_e}) to $\cM^\star_{g,\mfs}(X,D,A,\nu)_{\Gamma}$ is smooth.
\item  If $g\!=\!0$, $\nu\!\equiv\!0$, and $A\!\neq\! 0$ or $k\!\geq\! 3$, the same statement holds for $J$ in a Baire set of second category $\tn{AK}^{\Gamma}(X,D)\!\subset\!\tn{AK}(X,D)$.
\eEn
\eTh

\noindent
In a future work, by considering local perturbations, in the sense of Kuranishi structures, and construing a gluing map (outlined in \cite[Sec~3.4]{FT1}), we will use the same techniques used in the proof of Theorem~\ref{TransGamma_thm}  to construct Kuranishi charts for $\ov\cM_{g,\mfs}(X,D,A)$.\\

\noindent
In the classical case (no $D$), $\vr$ is the trivial map $\Z^\E\lra 0$.
If $N\!=\!1$ ($D$ smooth), $\tn{dim}~\D-\tn{dim}~\T \geq 0$ and $\tn{dim}~\tn{ker}(\vr)$ is often very large.
If $N\!>\!1$, there are nodal configurations $\Gamma$ with arbitrary many nodes, $\tn{dim}~\D-\tn{dim}~\T < 0$, and 
$\dim_\R \K_\R(\Gamma)=1$; see \cite[Ex~3.12]{FT1}. As a consequence, if $D$ is not empty or smooth, we need a sharp dimension counting argument for dealing with non-simple nodal maps.

\subsection{Semi-positive pairs}\label{SPLP_ss}
\textit{Our final goal is to introduce a class of semi-positive pairs $(X,D)$ for which one can use the perturbed moduli spaces above to construct log Gromov-Witten invariants in arbitrary genera without constructing a virtual fundamental class.}\\

\noindent
For $(X,D\!=\!\bigcup_{i\in [N]} D_i,\om)$ as before, we say $D$ is \textbf{Nef}\footnote{The terminology is inspired by the corresponding notion in algebraic geometry but is different.} if $A\cdot D_i\!\geq\! 0$ for all $A\!\in\! \pi_2(X)$ such that $\om(A)\!>\!0$.
Let 
$$
\ell_{I,A}=\tn{min}(A\cdot \bigcup_{j\in [N]-I} D_j,2)\qquad \forall~A\!\in\!H_2(X,\Z),\quad I\subset [N].
$$

\bDf{NCsemipos_dfn}
Let $(X,\om)$ be a closed $2n$-dimensional symplectic manifold and $D$ be a Nef SNC symplectic divisor in $(X,\om)$.  We say $(X,D,\om)$  is  \textbf{semi-positive} if 
\bEq{SP1_e}
c_1^{TX(- \log D)}(A)\!\geq 3\!-\!n\!+\!I\!-\!\ell_{I,A}   ~~\Rightarrow~~ c_1^{TX(-\log D)}(A)\!\geq\! 0
\eEq
for all $A\!\in\! \pi_2(D_I)$ such that $\om(A)\!>\!0$. We say $(X,D,\om)$  is  \textbf{positive} if  
\bEq{PExtra_e}
c_1^{TX(- \log D)}(A)\!\geq 3\!-\!n\!+\!I\!-\!\ell_{I,A} ~~\Rightarrow~~ c_1^{TX(-\log D)}(A)\!>\! 0
\eEq
for all $A\!\in\! \pi_2(D_I)$ such that $\om(A)\!>\!0$. 
\eDf

\noindent
We say $[X,D,\om]$ is semi-positive/positive if $(X,D,\om')$ for some $\om'$ deformation equivalent to $\om$ is semi-positive/positive. One may remove the Nef condition at the expense of altering the left-hand side. \\

\noindent
In Section~\ref{RMCMap_ss}, we show that under the semi-positivity condition (\ref{SP1_e}) (resp. positivity condition (\ref{PExtra_e})), multiple-cover log $J$-holomorphic spheres happen in dimensions less than or equal (resp. less) than somewhere injective maps. As we explain below, unlike the classical case, this is not sufficient for finding a suitable upper bound for the dimension of the image of non-simple maps in $\ov\cM_{g,k}\times X^{\mfs}$ when $D$ is not smooth.\\

\noindent
For every decorated dual graph $\Gamma=\Gamma(\V,\E,\L)$ with $|\V|\!\geq\!2$, let
$$
\cM^{\tn{ns}}_{g,\mfs}(X,D, A,\nu)_\Gamma=\cM_{g,\mfs}(X,D,A,\nu)_\Gamma-\cM^{\star}_{g,\mfs}(X,D,A,\nu)_\Gamma
$$
be the subset of \textbf{non-simple} (or \textbf{multiply-covered}) maps. An arbitrary stable log  $(J,\nu)$-map fails to be simple if either it contains a non-trivial bubble component that is a multiple-cover, or if it contains two non-trivial bubbles with the same image. 

\bPr{MC_pr} If $D$ is smooth,
\bEn
\item under the semi-positivity condition of Definition~\ref{NCsemipos_dfn}, for $(\om',\cR,J,\nu)$ in a Baire set of second category $\cH^{\Gamma,\tn{ns}}_{g,k}(X,D)_{[\om]}\!\subset\!\cH_{g,k}(X,D)_{[\om]}$,
the image of  $\cM^{\tn{ns}}_{g,\mfs}(X,A,J,\nu)_\Gamma$ under $\st\times \ev$
lies in the image of smooth maps from finitely many smooth even-dimensional manifolds of at least $2$ real dimension less than the dimension of the main stratum; 
\item similarly, under the same condition, if $g\!=\!0$, $k\leq2$, and $\nu\!\equiv\!0$, for $(\om',\cR,J)$ in a Baire set of second category $\tn{AK}^{\Gamma,\tn{ns}}(X,D)_{[\om]}\!\subset\!\tn{AK}(X,D)_{[\om]}$, the image of  $\cM^{\tn{ns}}_{0,\mfs}(X,D,A)_\Gamma$ under $\st\times \ev$ lies in the image of smooth maps from finitely many smooth even-dimensional manifolds of at least $2$ real dimension less than the expected dimension.
\eEn
\ePr

\noindent
For $D$ smooth, the semi-positivity/positivity conditions of Definition~\ref{NCsemipos_dfn} are essentially the same as \cite[Dfn~4.7]{FZ1}. There is some confusion in \cite{IP1, IP2} regarding the proper semi-positivity requirements in the relative case; see \cite[Rmk~4.9]{FZ1}.
Their work does not include a detailed proof of the relative analogue of Proposition~\ref{MC_pr}. The claim is that the result  follows from the classical result of Ruan-Tian by looking at the image of the moduli space in $\ov\cM_{g,k}(X,A,\nu)$. \\

\noindent
The first statement below follows from Theorem~\ref{TransGamma_thm}.(1) and Proposition~\ref{MC_pr}.(1). The second follows from a standard family version of these results. The third, follows by lifting two perturbations $\nu_1$ and $\nu_2$ obtained from two regular families $\mf{U}_{g,k}^{(1)}\lra \mf{M}_{g,k}^{(1)}$ and $\mf{U}_{g,k}^{(2)}\lra \mf{M}_{g,k}^{(2)}$ to their fiber product regular family and connecting them by a regular path, as in \cite[p. 34-35]{Z2}.

\bCr{main_cr}
Suppose $(X,\om)$ is a closed symplectic manifold, $D$ is a smooth symplectic divisor,  $A\!\in\!H_2(X,\Z)$, $g,k\!\in\!\N$ with $2g\!+\!k\!\geq\!3$,  and $\mfs\!\in\!\N^k$. If $[X,D,\om]$ is  semi-positive, for any choice of a regular universal family as in (\ref{UFamily_e}), there exists a Baire subset  
$$
\cH^{\tn{reg}}_{g,k}(X,D)_{[\om]}\!\in\!\cH_{g,k}(X,D)_{[\om]}
$$ 
of the second category such that for every $(\om',\cR,J,\nu)$ in this set
\bEn
\item the map 
$$
\st\times \ev\colon\!\cM_{g,\mfs}(X,D,A,\nu)\!\lra \ov\cM_{g,k}\!\times\!X^\mfs 
$$
defines a pseudo-cycle of real dimension
$$ 
d^{\log}=2\Big(c_1^{TX(-\log D)}(A)+(n-3)(1-g)+k\Big);
$$
\item the integral homology class $\wt{\GW}^{X,D}_{g,\mfs,A}$ in $\ov\cM_{g,k} \times X^\mfs$ determined by this pseudo-cycle is independent of the choice of $(\om',\cR,J,\nu)$; 
\item furthermore, the rational class
$$
\GW^{X,D}_{g,\mfs,A}\equiv \frac{1}{\tn{deg}~p}\wt{\GW}^{X,D}_{g,\mfs,A}\in H_{2d^{\log}}(\ov\cM_{g,k}\!\times\!X^\mfs, \Q),
$$
where $\tn{deg}~p$ is the degree of the regular covering used to define $\nu$, is an invariant of the deformation equivalence class of 
$\om\!\in\!\tn{Symp}(X,D)$.
\eEn
\eCr

\noindent
Similarly, the first statement below follows from Theorem~\ref{TransGamma_thm}.(2), Proposition~\ref{MC_pr}.(2), and Lemma~\ref{MCLemma}. The second one follows from a family version of these results.  

\bCr{main2_cr}
Suppose $(X,\om)$ is a closed symplectic manifold, $D$ is a smooth symplectic divisor,  $A\!\in\!H_2(X,\Z)$, $k\leq 2$,  and $\mfs\!\in\!\N^k$. If $[X,D,\om]$ is positive, there exists a Baire subset  
$$
\tn{AK}^{\tn{reg}}(X,D)_{[\om]}\!\in\!\tn{AK}(X,D)_{[\om]}
$$ 
of the second category such that for every $(\om',\cR,J)$ in this set
\bEn
\item the map 
$$
\st\times \ev\colon\!\cM^\star_{0,\mfs}(X,D,A)\!\lra \ov\cM_{0,k}\!\times\!X^\mfs 
$$
defines a pseudo-cycle of real dimension
$$ 
2\Big(c_1^{TX(-\log D)}(A)+n-3+k\Big);
$$
\item and, the integral homology class $\GW^{X,D}_{0,\mfs,A}$ in $\ov\cM_{0,k} \!\times\! X^\mfs$ determined by this pseudo-cycle is independent of the choice of $(\om',\cR,J)$.
\eEn
\eCr

\noindent
For the reason stated at the last paragraph of Section~\ref{Transversality_ss}, Definition~\ref{NCsemipos_dfn} is not strong enough for proving Proposition~\ref{MC_pr} (and consequently Corollaries~\ref{main_cr} and \ref{main2_cr}) for arbitrary SNC case ($N\!>\!1$). More precisely, in the proof of the classical version of  Proposition~\ref{MC_pr}, we get a stratification
$$
\cM^{\tn{ns}}_{g,k}(X, A,\nu)_\Gamma= \bigcup_{\gamma} \cM^{\gamma}_{g,k}(X, A,\nu)_\Gamma
$$ 
where $\gamma=(\Gamma,\Gamma'')$ and $\Gamma''$ characterizes the topological type of the underlying simple curves. For each $\gamma$, we get a fiberation 
$$
\pi_\gamma\colon \cM^{\gamma}_{g,k}(X, A,\nu)_\Gamma\lra \cM_{g,k}(X, A,\nu)_{\Gamma''}
$$
where the image is the space of underlying simple curves, and $\tn{st}\times\tn{ev}$ factors through $\pi_\gamma$. In the logarithmic/relative case, we first consider a pre-log space $\cM_{g,\mfs}^{\tn{plog}}(X,D,A,\nu)_\Gamma$ for which similar fiberations
$$
\pi_\gamma\colon \cM^{\tn{plog},\gamma}_{g,\mfs}(X, D,A,\nu)_\Gamma\lra \cM_{g,\mfs}(X,D, A,\nu)_{\Gamma''}
$$ 
can be constructed. Then 
$$
 \cM^{\gamma}_{g,\mfs}(X, D,A,\nu)_\Gamma=\tn{ob}_{\Gamma}^{-1}(1),
$$
where 
$$
\tn{ob}_{\Gamma}\colon \cM^{\tn{plog},\gamma}_{g,\mfs}(X, D,A,\nu)_\Gamma\lra \mc{G}(\Gamma)
$$
is the obstruction map in Condition~\ref{GObs_it} of Page~\pageref{GObs_it}. The point is that $\tn{ob}_{\Gamma}$ does not factor through $\pi_\gamma$. If $D$ is smooth, we can ignore $\tn{ob}_\Gamma$ because $\tn{dim}~\D-\tn{dim}~\T \geq 0$. Otherwise, as Example~\ref{MCissue_ex} shows, $\cM^{\gamma}_{g,\mfs}(X, D,A,\nu)_\Gamma$ can be larger than the main stratum. For this reason, we need a stronger semi-positivity condition that yields a bound on the dimension of $\cM^{\tn{plog},\gamma}_{g,\mfs}(X, D,A,\nu)_\Gamma$.\\

\noindent
Given $[X,D,\om]$, for each $I\!\subset\![N]$ and $A\!\in\!\pi_2(X)$, let $\de_{I,A}$ denote the minimum number of (geometric) intersection points of a degree $A$ $J$-holomorphic sphere\footnote{It is possible to write down a weaker definition of $\de_{I,A}$ without mentioning $J$-holomorphic spheres. Also, we define this number to be zero if such a $J$-holomorphic sphere does not exist.} in $D_I$ with $\bigcup_{i\in [N]-I}D_i$, for all $(\om',J)\!\in\!\tn{AK}(X,D)_{[\om]}$. 

\bDf{NCsemipos_dfnN}
Let $(X,\om)$ be a closed $2n$-dimensional symplectic manifold and $D$ be a Nef SNC symplectic divisor in $(X,\om)$.  We say $(X,D,\om)$  is  \textbf{strongly-semi-positive} if  
\bEq{SP2_eN}
c_1^{TX(- \log D)}(A)\!\geq 3\!-\!n\!+|I|\!-\!\ell_{I,A}   ~~\Rightarrow~~ c_1^{TX(-\log D)}(A)\!\geq\! \tn{max}\{0,2-\de_{I,A}\}
\eEq
for all $A\!\in\! \pi_2(D_I)$ such that $\om(A)\!>\!0$ (except possibly $(I,\de_{I,A})\neq (\eset,0)$).
We say $(X,D)$ is {strongly-positive} if it is strongly-semi-positive and positive.
\eDf

\noindent
In many examples, Condition~\ref{SP2_eN} and \ref{SP1_e} are equivalent.

\bEx{ProjSpace_ex}
A transverse union of $d$ hyperplanes in $\P^n$ is semi-positive ($=$ strongly-semi-positive) whenever $d\!\notin \![n+2,2n+1]$ and it is positive ($=$ strongly-positive) whenever $d\!\notin \![n+1,2n+1]$. An interesting case is when $X\!=\!\P^n$ and $D$ is the degree $n+1$ toric boundary divisor $\partial \P^n$. In this case $TX(-\log D)$ is the trivial complex vector bundle generated by (the pre-image in $TX(-\log D)$ of $n\!+\!1$ vector fields 
$$
x_0\partial x_0, \ldots, x_{n}\partial x_{n}, \qquad\tn{satisfying}\quad \sum_{i=0}^{n} x_i\partial x_i=0.
$$
In other words, $(\P^n,\partial \P^n)$ is a log Calabi-Yau pair.
If furthermore $n\!=\!3$, then 
$$
\tn{exp-dim}_{\C}~\ov\cM^{\log}_{g,\mfs}(X,D,A)=k
$$
is independent of $g$ and $A$. Note that $k$ is at least $2$. Since the quintic Calabi-Yau 3-fold can be degenerated to an SNC configuration whose components are  blowups of $(\P^3,\partial \P^3)$, the GW numbers arising from $\ov\cM^{\log}_{g,\mfs}(\P^3,\partial \P^3,A)$ should be related to the GW invariants of the quintic CY 3-fold.
\eEx

\noindent
The claim is that under the stronger conditions of Definition~\ref{NCsemipos_dfnN}, Proposition~\ref{MC_pr} (and consequently Corollaries~\ref{main_cr} and \ref{main2_cr}) holds for all $N$.

\bCn{Claim}
Suppose $X$ is closed, $D=\bigcup_{i\in [N]} D_i\subset X$ is an SNC symplectic divisor,  $A\!\in\!H_2(X,\Z)$, $g,k\!\in\!\N$ with $2g\!+\!k\!\geq\!3$,  and $\mfs\!\in\!(\N^N)^k$. If $[X,D,\om]$ is  strongly-semi-positive (resp. strongly-positive), then Proposition~\ref{MC_pr}, and consequently Corollaries~\ref{main_cr} and \ref{main2_cr}, hold.
\eCn

\noindent
We lay the foundation and explain the difficulties for proving Conjecture~\ref{Claim} in Section~\ref{MCMap_ss}. We provide several examples that illustrate the issues. We also explain the consequence of (\ref{SP2_eN}).
We plan to address this conjecture in a future work. There are, however, some special but interesting cases, such as when $X$ is toric and $D$ is its boundary divisor (see Remark~\ref{de2_rmk}), or when $D_i$ are $(0,A)$-hollow (Donaldosn divisors of sufficiently high degree) in the sense of \cite[Dfn~1]{FZ2}, where Conjecture~\ref{Claim} can be confirmed with easier arguments. \\
%\bPr{Toric_prp}
%Suppose $\de_{I,A}\!\geq\!2$ for all $A\!\in\! \pi_2(D_I)$ such that $\om(A)\!>\!0$ (except possibly $(I,\de_{I,A})\neq (\eset,0)$). Then Conjecture~\ref{Claim} is correct.
%\ePr

\noindent
\textbf{Acknowledgments.} I would like to thank A. Zinger, D. Pomerleano, and the  referee of \cite{FT1} for their helpful comments.
%------------------------------------------------------------------------------------------------------
%---------------------SNC structures
%------------------------------------------------------------------------------------------------------
\section{SNC divisors and the associated structures; review}\label{SC_s}
\noindent
In Section~\ref{DivDef_ss}, we recall the notions of simple normal crossings (SNC) symplectic divisor and symplectic regularizations for such objects introduced in \cite{FMZ1}. Regularizations allow us to define a suitable space of almost K\"ahler structures in Sections~\ref{Aux_ss} and perturbations in Section~\ref{Perturbation_ss}. In Section~\ref{LogBundle_ss}, we review the notion of logarithmic tangent bundle associated to SNC symplectic divisors  introduced in \cite{FMZ2}.
Readers familiar with the definitions, notations, and results of \cite{FMZ1,FMZ2} may skip this section. We refer to \cite{FMZ2} for a relatively short review of these concepts.

%---------------------------------------------------
\subsection{SNC divisors and regularizations}\label{DivDef_ss}
\noindent
Let $X$ be a (smooth) manifold. 
For any submanifold $D\!\subset\!X$, let
$$
\cN_XD\equiv \frac{TX|_D}{TD}\lra D
$$
denote the normal bundle of~$D$ in~$X$.

\begin{definition}[{\cite[Definition~2.1]{FMZ1}}]\label{SCD_dfn}
An \textbf{SNC symplectic divisor} in a  symplectic manifold $(X,\om)$ 
is a finite transverse union $D\!\equiv\!\bigcup_{i \in [N]} D_i$ 
of smooth symplectic divisors $\{D_i\}_{i\in [N]}$ such that for every $I\!\subset\![N]$
the submanifold
$$
D_I\!\equiv\!\bigcap_{i\in I} D_i\subset X
$$
is symplectic and its symplectic and intersection orientations are the same. 
\eDf

\noindent
Let $\tn{Symp}(X,D)$ (this is denoted by $\tn{Symp}^+(X,D)$  in \cite{FMZ1}) denote the space of symplectic structures $\om$ on $X$ such that $D$ is an SNC divisor in $X$ with respect to $\om$. The particular choice of $\om$ in its deformation equivalence class is not important in the construction of (relative) GW invariants. \\

\noindent
By the transversality assumption, the homomorphisms 
\bEq{DecNI_e}
\cN_XD_I\lra \bigoplus_{i\in I}\cN_XD_i\big|_{D_I} \qquad \forall~I\!\subset\![N],
\eEq
induced by the inclusions $TD_I\!\subset\!TD_i|_{D_I}$ are isomorphisms. These vector bundle isomorphisms are not symplectic unless $\{D_i\}_{i\in [N]}$ intersect orthogonally.
For $I'\!\subset\!I\!\subset\![N]$, define
$$
\cN_{I;I'}=\bigoplus_{i\in I-I'}\!\!\cN_XD_i|_{D_I}\,;
$$
under the decomposition (\ref{DecNI_e}), $\cN_{I;I'}$ is isomorphic to the normal bundle of $D_I$ in $D_{I'}$.
We denote by 
$$
\pi_{I;I'} \colon\! \cN_{I;I'}\lra D_I, \quad
\pi_I\colon\!\cN_X D_I\lra D_I,
$$
the natural projection maps.\\

\noindent
A \textbf{system of regularizations for} $\{D_i\}_{i\in [N]}$ in~$X$ is
a collection of smooth embeddings
$$
\Psi_I\!: \cN'_X D_I \lra X \qquad \forall~I \!\subset\! [N],
$$
from open neighborhoods $\cN'_XD_I\!\subset\!\cN_X D_I$ of~$D_I$
so that $\Psi_I|_{D_I}\!=\!\id_{D_I}$,
$\nd\Psi_I$ induces the identity map on~$\cN_XD_I$, and
$$
\Psi_I\big(\cN_{I;I'}\!\cap\!\tn{Dom}(\Psi_I)\big)=D_{I'}\!\cap\!\tn{Im}(\Psi_I)
\qquad\forall~I'\!\subset\!I\!\subset\![N]\,.
$$
This implies that $\nd\Psi_I$ induces an isomorphism
\bEq{fDIIdfn_e}
\mf{D}\Psi_{I;I'}\!:  \pi_{I;I'}^*\cN_{I;I-I'}\big|_{\cN_{I;I'}\cap\tn{Dom}(\Psi_I)}
\lra \cN_XD_{I'}\big|_{D_{I'}\cap\tn{Im}(\Psi_I)};
\eEq
see \cite[Section~2.2]{FMZ1}.
In the $I\!=\!I'$ case, this derivative is the identity map.

\bDf{RegV_dfn}
A \textbf{regularization} for $D$ in $X$ is a
system of regularizations for $\{D_i\}_{i\in [N]}$ in~$X$ as above
such~that
$$\tn{Dom}(\Psi_I)=\mf{D}\Psi_{I;I'}^{-1}\big(\tn{Dom}(\Psi_{I'})\big), 
\quad  \Psi_I\!=\!\Psi_{I'}\!\circ\!\mf{D}\Psi_{I;I'}|_{\tn{Dom}(\Psi_I)}
\qquad \forall~I'\!\subset\!I\!\subset\![N]\,.
$$
\eDf
\vskip.1in 

\noindent
For each $i\!\in\![N]$, an \textbf{$\om$-compatible Hermitian structure} on $\cN_XD_i$ is a triple $(\mf{i}_{i},\rho_{i},\na^{(i)})$,
where $\mf{i}_{i}$ is an $\om$-compatible (fiber-wise) complex structure on~$\cN_XD_i$, 
$\rho_{i}$ is a Hermitian metric with real part
$$
\rho_{\R}(\cdot,\cdot)=\om|_{\cN_XD_i}(\cdot,\mf{i}_{i}\cdot),
$$ 
and $\na^{(i)}$ is a Hermitian connection compatible with $(\mf{i}_{i},\rho_{i})$.
For each $i\!\in\![N]$, the space of $\om$-compatible Hermitian structures on $\cN_XD_i$ is non-empty and contractible. 
Each triple $(\mf{i}_{i},\rho_{i},\na^{(i)})$ as above
determines a 1-form $\al_{i}$ on $\cN_XD_i\!-\!D_i$ 
whose restriction to each fiber $\cN_XD_i|_x\!-\{x\}\!\cong\!\C^*$ 
is the 1-form $\nd\theta$ with respect to the polar coordinates $(r,\theta)$ on~$\C$. 
We also denote by~$\rho_{i}$ the square-of-the-norm function on~$\cN_XD_i$.
For each $i\!\in\!I\!\subset\![N]$, we denote the tuple induced by $(\mf{i}_{i},\rho_{i},\na^{(i)},\al_{i})$ on $\cN_XD_i|_{D_I}$ by $(\mf{i}_{I;i},\rho_{I;i},\na^{(I;i)},\al_{I;i})$.

\begin{definition}[{\cite[Definition~2.9]{FMZ1}}]\label{omRegV_dfn}
If $D\equiv \bigcup_{i\in [N]} D_i$ is an SNC symplectic divisor in $(X,\om)$, an \textbf{$\om$-regularization} for $D$ in $X$ consists of a choice of Hermitian structure
$(\mf{i}_{i},\rho_{i},\na^{(i)})$ on $\cN_X D_i$
for all $i\!\in\![N]$ together with
a regularization for~$D$ in~$X$ as in Definition~\ref{RegV_dfn} so that
\bEq{omStd_e}
\Psi_I^*\om\!=\!\pi_I^*(\om|_{TD_I})+\frac{1}{2} \sum_{i\in I} \nd(\rho_{I;i} \al_{I;i})
\qquad \forall~I\!\subset\! [N],
\eEq
and (\ref{fDIIdfn_e}) is an isomorphism of split Hermitian vector bundles
for all $I'\!\subset\!I\!\subset\![N]$. 
\eDf

\noindent
If $N\!=\!1$, i.e.~$D$  is a smooth divisor, an $\om$-regularization is a single map
\bEq{PsiMaps_e}
\Psi\!: \cN'_X D \lra X
\eEq
as in \cite[Lmm~3.14]{MS2} without any further compatibility condition.
We define the space of auxiliary data $\tn{Aux}(X,D)$ to be the space of pairs 
$(\om,\cR)$, where $\om\!\in\!\Symp(X,D)$ and $\cR$ 
is an $\om$-regularization of~$D$ in~$X$. 
Let 
$$
\Pi\!:\Symp(X,D)\lra H^2(M;\R)
$$ 
be the map sending $\om$ to its de Rham equivalence class~$[\om]$.
The following is a weaker version of the main result of \cite{FMZ1} 
for SNC symplectic divisors.

\begin{theorem}[{\cite[Theorem 2.13]{FMZ1}}]\label{SymptoAux_Thx}
For $(X,D)$ as above, the projection~maps 
$$
\pi\!:\tn{Aux}(X,D)\lra \Symp(X,D),\qquad
\pi\big|_{\Pi^{-1}(\al)}\!: \{\Pi\!\circ\!\pi\}^{-1}(\al) \lra \Pi^{-1}(\al), 
~~ \al\!\in\!H^2_{\tn{dR}}(M),$$
are weak homotopy equivalences.
\eTh

\noindent
Given $\om$-regularizations $\cR$ and $\cR'$, we say $\cR'$ is a \textbf{shrinking} of $\cR$ if the Hermitian data in $\cR$ and $\cR'$ are the same and  
$$
\tn{Dom}(\Psi'_I)\!\subset\tn{Dom}(\Psi_I) \quad \tn{and}\quad \Psi'_I=\Psi_I|_{\tn{Dom}(\Psi'_I)}\qquad \forall~I\!\subset\![N].
$$ 
Two $\om$-regularizations $\cR$ and $\cR'$ are said to be \textbf{equivalent} if they have a common shrinking. The latter defines an equivalence relation among $\om$-regularizations.\\

\noindent
If $D$ is an SNC symplectic divisor in $(X,\om)$, then, for each $I\!\subset[N]$, 
\bEq{parDI_e}
\partial D_I= \bigcup_{I\subsetneq J} D_J
\eEq
is an SNC symplectic divisor in $(D_I,\om|_{TD_I})$ and an $\om$-regularization $\cR$ for $D$ in $X$ restricts to an $\om|_{TD_I}$-regularization $\cR_I$ for $\partial D_I$ in $D_I$. For $I\!=\!\eset$, the convention is $D_\eset\!=\!X$ and $\partial X \!=\!D$.

%---------------------------------------------------
\subsection{Almost K\"ahler auxiliary data}\label{Aux_ss}
\noindent
In order to define relative GW invariants of $(X,D,\om)$, we need to consider an almost complex structure $J$ on $X$ that is both $\om$-tame and $D$-compatible. The biggest set of such almost complex structures that one may consider is the set of $\om$-tame (or compatible) $J$ such that $J(TD_i)\!=\!TD_i$ for all $i\!\in\![N]$ and 
\bEq{intInnormal_e}
N_J(u,v)\in TD_i \qquad \forall~i\!\in\![N],~x\!\in\!D_i,~u,v\!\in\!T_xX,
\eEq
where $N_J$ is Nijenhueis tensor of $J$. Condition~(\ref{intInnormal_e}) is needed to ensure that certain operators are complex linear (see \cite[(2.20)]{FT1}), or equally, certain sequence of almost complex structures on $\cN_XD_i$ converges to a standard one (see  \cite[Lmm~4.5]{FT1}). In this paper, however, similarly to \cite{LR} and \cite{EGH}, we restrict to a special class of almost complex structures arising from regularizations in the following sense.

\bDf{stdJ_dfn}
Suppose $D\!=\! \bigcup_{i\in [N]} D_i$ is an SNC symplectic divisor in $(X,\om)$, $\cR$ is an $\om$-regularization for $D$ in $X$ as in Definition~\ref{omRegV_dfn}, and $J$ is an $\om$-tame almost complex structure on $X$ such that  $J(TD_i)\!=\!TD_i$ for all $i\!\in\![N]$. We say $J$ is $\cR$-compatible, if 
\bEq{stdJ_e}
\Psi_{I}^*J=  \pi_{I}^*(J_I)\oplus \pi_I^* \bigoplus_{i\in I}  \mfi_{I;i}\qquad \forall~I\subset [N],
\eEq
where $J_I\!\equiv\! J\big|_{TD_{I}}$, $\mfi_{I;i}$ is the complex structure on $\cN_{X}D_i|_{D_I}$ pre-determined in $\cR$, and the righthand side of (\ref{stdJ_e}) is the direct sum complex structure corresponding to the decomposition
\bEq{Tdecomp_e}
T\cN_XD_I\cong \pi_{I}^*TD_{I}\oplus \pi_{I}^*\cN_XD_{I},
\eEq
given by the connections in $\cR$.
\eDf

\noindent
Given any $(\om,\cR)\!\in\!\tn{Aux}(X,D)$, we denote the space of $\om$-compatible almost complex structures compatible with a shrinking of $\cR$ by $\AK(X,D,\om)_\cR$. This space is non-empty and contractible. We denote the space of compatible tuples $(\om,\cR,J)$ by $\AK(X,D)$. It follows from Theorem~\ref{SymptoAux_Thx} that the projection map
\bEq{AKtoSymp_e}
\AK(X,D)\!\lra\!\Symp(X,D) 
\eEq
is a weak-homotopy equivalence. This implies that any invariant of the deformation equivalence classes in $\AK(X,D)$ is an invariant of the symplectic deformation equivalence classes in $\tn{Symp}(X,D)$.
%---------------------------------------------------
\subsection{Logarithmic tangent bundle}\label{LogBundle_ss}
\noindent
In this section, we review the notion of logarithmic tangent bundle in complex geometry and the analogous notion for SNC symplectic divisors introduced in \cite{FMZ2}. A detailed description of this construction will appear in \cite{FMZ4}. We will show in Section~\ref{Deformation_s} that the linearization of Cauchy-Riemann operator for log maps is a lift of the classical linearization map to the logarithmic tangent bundle.\\

\noindent
Let $X$ be a smooth holomorphic manifold and $D\!\subset\!X$ be a normal crossings divisor; i.e. locally around every point $p\!\in\!X$ there are holomorphic coordinates $(x_1,\ldots,x_n)$, with $n\!=\!\dim_\C X$, such that 
$$
D\!\equiv \!( x_1\cdots x_k\!=\!0) \subset X\quad\tn{for~some}~~~k\leq n.
$$
In such coordinates, the sheaf $\cT X$ of holomorphic sections of the complex tangent bundle $TX$ is generated by 
$$
\partial_{x_1}, \cdots, \partial_{x_n}
$$
and the \textbf{log tangent sheaf} $\cT X(-\log D)$ is the sub-sheaf  generated by 
$$
x_1\partial_{x_1}, \ldots, x_k \partial_{x_k},
~~~\partial_{x_{k+1}}, \ldots, \partial_{x_n}.
$$
It is dual to the sheaf $\Omega^{1}_X(\log D)$ of meromorphic $1$-forms with at most simple poles along $D_i$.
Since $\cT X(-\log D)$ is locally free, it is the sheaf of holomorphic sections of a holomorphic vector bundle $TX(-\log D)$. The inclusion $\cT X(-\log D)\!\subset\!\cT X$ gives rise to a holomorphic homomorphism
$$
\iota\colon T X(-\log D)\lra T X
$$
which is an isomorphism away from $D$. \\

\noindent
Given an SNC symplectic divisor  $D\!=\! \bigcup_{i\in [N]} D_i$ in $(X,\om)$, an $\om$-regularization $\cR$ for $D$ in $X$ as in Definition~\ref{omRegV_dfn}  gives rise to a real rank-$2n$ vector bundle 
$$
TX(-\log D)_{\cR} \lra X
$$
satisfying
$$
 TX(-\log D)_{\cR}|_{D_I}=  TD_I(-\log \partial D_I)\oplus D_I\!\times\! \C^I \qquad~\forall~I\!\subset\![N],
$$
where $\partial D_I\!\subset\!D_I$ is the SNC divisor in (\ref{parDI_e}).
As in the holomorphic case, there exists a canonical homomorphism 
\bEq{LogToTX_e}
\iota\colon TX(-\log D)_{\cR}\lra TX
\eEq
which is an isomorphism away from $D$.
Furthermore, any $\cR$-compatible almost complex structure $J$ on $TX$ gives rise to a similarly denoted complex structure on $TX(-\log D)_{\cR}$  such that the homomorphism (\ref{LogToTX_e}) is complex linear. 
The deformation equivalence class $TX(-\log D)$ of the complex vector bundle $(TX(-\log D)_{\cR},J)$, which we call the \textbf{log tangent bundle of} $(X,D)$, only depends on the deformation equivalence class of  $(X,D,\om)$. Furthermore, 
$$
c\big(TX(-\log D)\big) = c(TX)/ \prod_{i\in [N]} \big(1+ \tn{PD}(D_i)\big);
$$
see \cite{FMZ2, FMZ4}. In particular,
$$
c_1^{TX(-\log D)} = c_1^{TX}- \sum_{i\in [N]} \tn{PD}(D_i).
$$

\noindent
For a smooth divisor $D\!\subset\!X$, with notation as in (\ref{PsiMaps_e}), we have
$$
\aligned
&TX(-\log D)_{\cR}= \big(\Psi^{-1\,*} (\pi^* TD \oplus \cN'_XD\!\times\!\C) \sqcup 
TX|_{X-D}\big)\big/\!\!\sim,\\
&\Psi^{-1\,*} (\pi^* TD \oplus \cN_X'D\!\times\!\C)
  \ni\big(\Psi(v),u \oplus c\big)\sim\big(\Psi(v),\nd_v\Psi(u+cv)\big)\in T(X\!-\!D),
  \endaligned
$$
where in the last equation, via the isomorphism $T\cN_XD\cong \pi^*TD\oplus \pi^*\cN_XD$ given by the connections in $\cR$, we think of 
$u+cv$ as a tangent vector in $T_v\cN_XD$, for all $v\!\in\!\cN_XD$. In the general SNC case, with notation as in Section~\ref{DivDef_ss}, we have
\bEq{TXlog_e}
TX(-\log D)_{\cR}= \bigg(\bigsqcup_{I\subset [N]}
\Psi_I^{-1\,*} \big(  \pi_{I}^*TD_{I} |_{\cN^\circ_XD_I} \oplus \cN^\circ_XD_I\times \C^I \big)\bigg)\Big/\!\!\sim
\eEq
where 
$$
\cN^\circ_XD_I = \pi_{I}^{-1}\big(D_I\!-\!\partial D_I\big)\cap \tn{Dom}(\Psi_I)
$$
and the identification maps that give rise to the equivalence relations $\sim$ on the overlaps are given by
$$
\aligned
&\pi_{I}^*TD_{I} |_{\cN^\circ_XD_I} \oplus \cN^\circ_XD_I\times \C^I \ni \big((p,v); \xi\oplus (c_i)_{i\in I}\big) \lra\\
 &\qquad \qquad\qquad\qquad\big((q,\tilde{v}); \wt{\xi}\oplus (c_i)_{i\in I'} \big) \in \pi_{I'}^*TD_{I'} |_{\cN^\circ_XD_{I'}} \oplus \cN^\circ_XD_{I'}\times \C^{I'}\qquad \forall~I'\!\subset\!I\!\subset\![N],\\
 & p\!\in\! D_I\!-\!\partial D_I,~v=(v_i)_{i\in I}\!\in\!\cN^\circ_XD_I|_{p}, ~q= \Psi_{I,I'}(p, (v_i)_{i\in I-I'})\!\in\!D_{I'}, \\
 & \wt{v}=\mf{D}\Psi_{I,I'} \big((p, (v_i)_{i\in I-I'}); (v_j)_{j\in I'}) \big)\!\in\!\cN_XD_{I'}|_{q},~\wt{\xi}=d\Psi_{I,I'}|_{(v_i)_{i\in I-I'}}\big(\xi+\sum_{i\in I-I'} c_i v_i\big) \in T_{q} D_{I'}.
\endaligned
$$
on the overlap; see \cite[Sec~2.2]{FMZ1}.
For each $\cR$-compatible $J$, on the local chart 
\bEq{Chart_e}
\pi_{I}^*TD_{I} |_{\cN^\circ_XD_I} \oplus \cN^\circ_XD_I\times \C^I\lra  \cN^\circ_XD_I,
\eEq
the complex structure is given by the pull back of $J|_{TD_I}$ on the first summand and the trivial complex structure on the second summand. Furthermore, the $\C$-linear homomorphism (\ref{LogToTX_e}) is given by 
\bEq{LogTrans_e}
\pi_{I}^*TD_{I} |_{\cN^\circ_XD_I} \oplus \cN^\circ_XD_I\times \C^I \ni \big((p,v); \xi\oplus (c_i)_{i\in I}\big) \lra d\Psi_I|_{v} \big(\xi \oplus \sum_{i\in I} c_i v_i\big) \in T_{\Psi_I(v)}X,
\eEq
where via the identifications (\ref{Tdecomp_e}) and (\ref{DecNI_e}) for each 
$$
v\cong(v_i)_{i\in I}\!\in\!\cN_XD_I|_{p} \cong \bigoplus_{i\in I} \cN_XD_i|_{p}
$$ 
we think of $\xi \oplus \sum c_i v_i$ as a tangent vector in $T_v\cN_XD_I$.  \\

\noindent
Let $h$ be a Hermitian metric  on $TX(-\log D)_{\cR}$.  For example, at the cost of shrinking $\cR$, one can construct $h$ so that on the chart (\ref{Chart_e}) it is the direct sum of the standard Hermitian metric on $\cN^\circ_XD_I\!\times\C$ and the pull back of some Hermitian metric from $D_I$ on the first summand. For every $I\!\subset\![N]$,  via the inclusion 
$$
TD_I(-\log \partial D_I)\subset TX(-\log D)|_{D_I}
$$ 
and the identification 
$$
TD_I(-\log \partial D_I)|_{D_I\!-\partial D_I}= T(D_I\!-\!\partial D_I),
$$
$h$ induces a complete Hermitian metric on $D_I\!-\!\partial D_I$. Let $\tn{exp}_I\colon T(D_I\!-\!\partial D_I)\lra D_I\!-\!\partial D_I$ be the exponential map of $h_I$.
Define 
\bEq{Log-exp_e}
\aligned 
&\tn{exp}^{\log}\colon TX(-\log D)_{\cR}\lra X,\quad (\xi\oplus c)|_p \lra \tn{exp}_I(p,\xi) \in D_I\!-\!\partial D_I,\\
 & \forall~I\!\subset\![N],~p\!\in\!D_I-\partial D_I,~\xi\oplus c \in T_pX(-\log D)_{\cR}\cong T_pD_I \oplus \C^I.
 \endaligned
\eEq
This map is smooth and it is a logarithmic version of the exponential map in the classical sense. In fact, for any 
$$
p\!\in\!D_I-\partial D_I\quad \tn{and} \quad v=(v_i)_{i\in I}\in \cN^\circ_XD_I|_{p},
$$ 
via the identification (\ref{LogTrans_e}), the map $\tn{exp}^{\log}$ is approximately (to the first order in $|v|$) equal to 
\bEq{exp-exp_e}
 (v, \xi\oplus c) \lra \Big(\tn{exp}_I(\xi), \tn{Pal}\big((\tn{e}^{c_i}v_i)_{i\in I}\big)\Big) \in \cN^\circ_XD_I,
\eEq
where $\tn{Pal}\big((\tn{e}^{c_i}v_i)_{i\in I}\big)$ is the parallel translation of the vector $(\tn{e}^{c_i}v_i)_{i\in I}\!\in\!\cN_XD_I$ along the path $\tn{exp}_I(t\xi)|_{t\in [0,1]}$ in $D_I$. Putting $v\!=\!0$, (\ref{exp-exp_e}) becomes $\tn{exp}_I$ at $p$. This logarithmic exponential map will be used in constructing a Banach neighborhood of a $(J,\nu)$-holomorphic map in the space of all smooth maps of the same contact type.

%------------------------------------------------------------------------------------------------------
%---------------------------------------Construction
%------------------------------------------------------------------------------------------------------
\section{Moduli spaces of log $(J,\nu)$-curves}\label{LogJnuMaps_s}
\noindent
In this section, following the description of \cite[Sec~2]{RT} and \cite[Sec~2.1]{Z2}, we define a suitable space of perturbation terms $\nu$ over a ``regular covering'' of the Deligne-Mumford space $\ov\cM_{g,k}$ for any symplectic logarithmic pair $(X,D)$.
Then, following and generalizing the definition of log $J$-holomorphic curves in \cite{FT1}, we introduce the notion of log $(J,\nu)$-holomorphic curve. 

%--------------------------------------------------------
\subsection{Regular coverings}\label{RegCover_ss}
\noindent
For $g,k\!\in\!\N$ with $2g\!+\!k\!\geq\!3$, let 
\bEq{UCurve_e}
\pi\colon \ov\cU_{g,k}\equiv \ov\cM_{g,k+1} \lra \ov\cM_{g,k}
\eEq
be the universal curve, where $\pi$ is defined by forgetting the last marked point. For a marked curve $C\!=\![\Si,\mfj, z_1,\ldots,z_k]$ in $\ov\cM_{g,k}$, if  the automorphism group $\tn{Aut}(C)$ is non-trivial, then $\pi^{-1}(C)= \Si/\tn{Aut}(C)$ instead of $\Si$. Therefore, unless $g\!=\!0$, (\ref{UCurve_e}) is not a universal family and we can not directly define the perturbation term $\nu$ over $\ov\cU_{g,k}$. One can resolve this problem by taking appropriate finite covers of (\ref{UCurve_e}).\\

\noindent
Denote by $\cT_{g,k}$ the Teichm\"uler space of genus $g$ Riemann surfaces with $k$ marked points (punctures) and by $\mc{G}_{g,k}$ the corresponding mapping class group. We have
$$
\cM_{g,k}=\cT_{g,k}/ \mc{G}_{g,k}.
$$
Assume $g\!=\!g_1\!+\!g_2$ and $k\!=\!k_1\!+\!k_2$ with $2g_i\!+\!k_i\!\geq\!3$ for $i\!=\!1,2$. For any decomposition $S_1\cup S_2$ of $[k]$ with $|S_i|=k_i$, there exists a canonical immersion 
\bEq{B1_e}
\iota=\iota_{S_1,S_2}\colon \ov\cM_{g_1,k_1+1}\times \ov\cM_{g_2,k_2+1}\lra \partial \ov\cM_{g,k}
\eEq
which assigns to a pair of marked  curves 
$$
\big(C_i\!=\![\Si_i,\mfj_i,z_{i,1},\ldots,z_{i,k_i+1}]\big)_{i=1,2},
$$ 
the marked  curve
$$
\aligned
&C=[\Si,\mfj,z_1,\ldots,z_k], \quad \Si\!=\! \Si_1\!\sqcup\! \Si_2/ z_{1,k_1+1}\!\sim \!z_{2,k_2+1}\\
&\{z_1,\ldots,z_k\}\!=\! \{z_{1,1},\ldots,z_{1,k_1}\}\cup  \{z_{2,1},\ldots,z_{2,k_2}\},
\endaligned
$$ 
so that the remaining marked points are renumbered by $\{1,\ldots,k\}$ according to the decomposition $S_1\!\cup\!S_2$. There is also another natural immersion 
\bEq{B2_e}
\de\colon \ov\cM_{g-1,k+2}\lra \partial \ov\cM_{g,k}
\eEq
which is obtained by gluing  together the last two marked points.

\begin{definition}[{\cite[Dfn~2.1]{Z2}}]\label{UnivFamily_dfn}
Let $g,k\!\in\!\N$ with $2g\!+\!k\!\geq\!3$, and 
\bEq{Cover_e}
p\colon \ov{\mf{M}}_{g,k}\lra\ov\cM_{g,k}
\eEq
be a finite branched cover in the orbifold category.
A \textbf{universal family over} $ \ov{\mf{M}}_{g,k}$ is a tuple 
\bEq{UFamily_e}
\bigg(\pi\colon\ov{\mf{U}}_{g,k}\lra \ov{\mf{M}}_{g,k}, \mf{z}_1,\ldots,\mf{z}_k\bigg)
\eEq
where $\ov{\mf{U}}_{g,k}$ is a complex projective variety and $\pi$ is a projective morphism with disjoint sections $\mf{z}_1,\ldots,\mf{z}_k$ such that for each $c\!\in\!\ov{\mf{M}}_{g,k}$  the tuple 
$$
C\!=\!\big((\Si,\mfj)\!=\!\pi^{-1}(c), \vec{\mf{z}}(c)=(\mf{z}_1(c),\ldots,\mf{z}_k(c))\big)
$$ 
is a stable $k$-marked genus $g$ curve  whose equivalence class is $[C]\!=\!p(c)$.
\eDf

\bDf{RefFamily_dfn}
Let $g,k\!\in\!\N$ with $2g\!+\!k\!\geq\!3$. A cover (\ref{Cover_e}) is \textbf{regular} if
\bEn
\item it admits a universal family,
\item each topological component of $p^{-1}\big(\cM_{g,k}\big)$ is the quotient of $\cT_{g,k}$ by a subgroup of $\mc{G}_{g,k}$,
\item\label{sepnode_l} for each boundary divisor (\ref{B1_e}) we have 
$$
\big(\ov\cM_{g_1,k_1+1}\times \ov\cM_{g_2,k_2+1}\big)\times_{(\iota,p)}\ov{\mf{M}}_{g,k} \approx \ov{\mf{M}}_{g_1,k_1+1}\times \ov{\mf{M}}_{g_2,k_2+1},
$$ 
for some regular covers $\ov{\mf{M}}_{g_i,k_i+1}$ of  $\ov\cM_{g_i,k_i+1}$, and
\item\label{selfnode_l} for the boundary divisor (\ref{B2_e}) we have 
$$
\ov\cM_{g-1,k+2}\times_{(\de,p)} \ov{\mf{M}}_{g,k} \approx \ov{\mf{M}}_{g-1,k+2},
$$ 
for some regular cover $\ov{\mf{M}}_{g-1,k+2}$ of $\ov\cM_{g-1,k+2}$.
\eEn
\eDf
\noindent
The last two conditions are inductively well-defined. This definition is a modified version of \cite[Dfn~2.2]{Z2}. In \cite[Dfn~2.2]{Z2}, the last condition is missing; furthermore, $\ov{\mf{M}}_{g_i,k_i+1}$ and $\ov{\mf{M}}_{g-1,k+2}$ are only required to be ``some" cover of $\ov\cM_{g_i,k_i+1}$ and $\ov\cM_{g-1,k+2}$, respectively.
The existence of such regular covers is a consequence of \cite[Prp~2.2, Thm~2.3, Thm~3.9]{BoPi}; see also moduli space of curves with \textit{level $n$ structures} in \cite[p~285]{Mu}. In the genus $0$ case, for each $k\!\geq\!3$, the moduli space $\ov\cM_{0,k}$ itself is smooth and the universal curve (\ref{UCurve_e}) is already a universal family. The regular covers are only branched over the boundaries of the moduli space. Furthermore, the total space of a universal family as in (\ref{UFamily_e}) over a regular cover only has singularities of the form
$$
\{(x,y,t)\in \C^3\colon~xy=t^m\} \lra \C, \qquad (x,y,t)\lra t
$$
at the nodal points of the fibers of $\pi$. In the approach of \cite{RT}, for dealing with such singularities they consider embeddings of a universal family into $\P^N$ for  sufficiently large $N$. In this article, following \cite{LR,Z2}, we consider perturbations supported away from the nodes.

%--------------------------------------------------------
\subsection{Logarithmic Ruan-Tian perturbations}\label{Perturbation_ss}
\noindent
Let $g,k\!\in\!\N$ with $2g\!+\!k\!\geq\!3$ and fix a regular covering (\ref{Cover_e}) and a universal family (\ref{UFamily_e}). Denote by 
$$
\ov{\mf{U}}_{g,k}^\star\subset \ov{\mf{U}}_{g,k}
$$
the complement of the nodes of the fibers of the projection map $\pi$ in (\ref{UFamily_e}). Denote by 
$$
T_{g,k}= \tn{Ker}~\nd (\pi|_{\ov{\mf{U}}_{g,k}^\star}) \lra \ov{\mf{U}}_{g,k}^\star
$$
the vertical tangent bundle. The latter is a complex line bundle; we denote the complex structure by $\mfj_{\mf{U}}$. 
Then
$$
\Om^{0,1}_{g,k}:= (T_{g,k},-\mfj_{\mf{U}})^* \lra \ov{\mf{U}}_{g,k}^\star
$$
is the complex line bundle of vertical $(0,1)$-forms. It is possible to extend this construction to the nodal points by allowing simple poles and dual residues, or by embedding $\ov{\mf{U}}_{g,k}$ into some $\P^M$ as in \cite{RT}. \\

\noindent
Let $(X,\om)$ be a  symplectic manifold and $J$ be an $\om$-tame almost complex structure on $X$.  The classical space of perturbations considered in \cite{RT} (following the modification in \cite{Z2}) is the infinite dimensional linear space
\bEq{cHgk_e}
\cH_{g,k}(X,J)\!=\!\big\{ \nu\! \in\! \Gamma\big(\ov{\mf{U}}_{g,k}^\star\times X, \pi_1^*\Om^{0,1}_{g,k}\otimes_\C \pi_2^*TX\big)~~\tn{s.t.}~~\tn{supp}(\nu)\!\subset(\ov{\mf{U}}_{g,k}^\star-\bigcup_{a\in [k]} \tn{Im}(\mf{z}_a))\!\times\! X \big\},
\eEq
where $\pi_1,\pi_2$ are projection maps from $\ov{\mf{U}}_{g,k}^\star\times X$ onto the first and second components, respectively, and $\tn{supp}(\nu)$ is the closure of the complement of the vanishing locus of $\nu$ in the compact space $\ov{\mf{U}}_{g,k}\times X$. Let $\cH_{g,k}(X,\om)$ denote the space of tuples $(J,\nu)$ where $J$ is $\om$-tame and $\nu\!\in\!\cH_{g,k}(X,J)$. Note that given $\nu$ and a boundary component as in Definition~\ref{RefFamily_dfn}.\ref{sepnode_l} (resp.  Definition~\ref{RefFamily_dfn}.\ref{selfnode_l}), the restriction of $\nu$ to $ \ov{\mf{M}}_{g_1,k_1+1}$ gives a perturbation term in $\cH_{g_1,k_1}(X,J)$ (resp. $\cH_{g-1,k+2}(X,J)$).

\bDf{JnuModuli_dfn}
Suppose $g,k\!\in\!\N$ with $2g\!+\!k\!\geq\!3$, $\ov{\mf{U}}_{g,k}$ is a universal family as in (\ref{UFamily_e}), $(X,\om)$ is a symplectic manifold, $A\!\in\!H_2(X,\Z)$, and $(J,\nu)\!\in\!\cH_{g,k}(X,\om)$. A \textbf{degree $A$ genus $g$ $k$-marked $(J,\nu)$-map} is a tuple 
\bEq{JnuMap_e}
f=\Big( \phi ,u, \big(\Si,\mfj, (z_a)_{a\in [k]}\big) \Big)
\eEq
where $(\Si,\mfj,(z_a)_{a\in [k]})$ is a nodal genus $g$  $k$-marked complex curve, $\phi\colon\!\Si\! \lra\! \ov{\mf{U}}_{g,k}$ is a holomorphic map onto a fiber of $\ov{\mf{U}}_{g,k}$ preserving the marked points,  and $u\colon\!\Si\!\lra\!X$ represents the homology class $A$ and satisfies 
$$
\dbar u =(\phi,u)^*\nu.
$$
\eDf
\noindent
Two $k$-marked $(J,\nu)$-holomorphic maps 
$$
\Big( \phi_1 ,u_1, \big(\Si_1, \mfj_1,(z_{1,a})_{a\in [k]}\big) \Big)\quad\tn{and}\quad\Big( \phi_2 ,u_2, \big(\Si_2,\mfj_2, (z_{2,a})_{a\in [k]}\big) \Big)
$$
are \textbf{equivalent} if there exists a holomorphic identification 
$$
h\colon \Si_1\!\lra\!\Si_2, \qquad h(z_{1,a})=z_{2,a} \quad \forall~a\!\in\![k],
$$
such that $(\phi_1,u_1)\!=\!(\phi_2,u_2)\circ h$. A $(J,\nu)$-holomorphic map is \textbf{stable} if it has a finite automorphism group. 
For any fixed $J$, denote by 
$$
\ov\cM_{g,k}(X,A,\nu)
$$
the space of equivalence classes of $k$-marked genus $g$ degree $A$ stable $(J,\nu)$-holomorphic maps. The subspace of maps with smooth domain is denoted by $\cM_{g,k}(X,A,\nu)$.
A \textbf{contracted} component of $\Si$ in (\ref{JnuMap_e}) is a smooth component whose image under the map $\phi$ is just a point. 
A map (\ref{JnuMap_e}) is stable  if and only if the degree of the restriction of $u$ to every contracted component of $\Si$ containing only one or two special (nodal or marked) points is not zero.
If (\ref{JnuMap_e}) is stable, every connected cluster of contracted components is a tree of spheres, with a total of at most $2$ special\footnote{either a marked point or a nodal point connecting the cluster to an irreducible component of $\Si$ outside the cluster.} points, at least one of which is a nodal point. 
For generic $\nu$, the only components of $\Si$ contributing non-trivially to the automorphism group of (\ref{JnuMap_e}) are the contracted components.

\bDf{Simple_dfn0}
A $k$-marked $(J,\nu)$-holomorphic map $f$ as in (\ref{JnuMap_e}) is called \textbf{simple} if  the restriction $u_v$ of $u$ to each irreducible component $\Si_v$ of $\Si$ contracted by $\phi$ is not multiply-covered (or equally it is somewhere injective) whenever\footnote{this is automatically satisfied if $f$ is stable.} $u_v$ is not constant, and the images of two such components in $X$ are distinct. 
\eDf

\noindent
In the case where $g\!=\!0$ and $\nu\!\equiv \!0$ (i.e. when we are dealing with $J$-holomorphic maps and the construction does not involve any regular cover, etc.), by a simple map we mean a nodal map such that the restriction $u_v$ of $u$ to each irreducible component $\Si_v\!\cong\!\P^1$ is not multiply covered, and the images of two such components in $X$ are distinct. \\

\noindent
If $J$ is an $\om$-tame almost complex structure on $X$, let $\na$ be the Levi-Civita connection of the metric $\ll u, v\rr=\frac{1}{2}(\om(u,Jv)+\om(v,Ju))$ and 
\bEq{CLConn_e}
\wt\na_{v} \ze\!=\!\na_v \ze - \frac{1}{2} J(\na_v J) \ze\!=\! \frac{1}{2}\big(\na_v \ze-J \na_v (J\ze)\big)\qquad \forall~v\!\in\!TX,~\ze\!\in\!\Gamma(X,TX)
\eEq
be the associated Hermitian connection. The torsion $T$ of the modified $\C$-linear connection
\bEq{CLConn_e2}
\wh\na_{v} \ze\!=\!\wt\na_v \ze - A(\ze)v, \quad A(\ze)= \frac{1}{4} \big(\nabla_{J\ze}J+ J\nabla_{\ze}J\big) \qquad \forall~v\!\in\!TX,~\ze\!\in\!\Gamma(X,TX)
\eEq
is related to the Nijenhueis tensor normalized as in \cite[p.18]{MS2} by
$$
T_{\wh\na}(v,w)=-\frac{1}{4}N_J(v,w)\qquad\forall~v,w\!\in\!TX.
$$
If $J$ is $\om$-compatible, $\wt\na$ coincides with $\wh\na$. See \cite[Ch~3.1 and Appendix C]{MS2} for details. 
By \cite[(3.1.6)]{MS2}, for a $(J,\nu)$-holomorphic map $(\phi,u)$, the linearization of  $\dbar-\nu$ has the form
\bEq{CRLNu_e}
\aligned
&\tn{D}_u \{\dbar-\nu\}\colon \Gamma(\Si,u^*TX)\lra \Gamma(\Si,\Om^{0,1}_{\Si,\mfj}\otimes_\C u^*TX),\\
&\tn{D}_u\{\dbar-\nu\} (\ze)\!=\! (\wh\nabla \ze)^{(0,1)}+ \frac{1}{4}N_J(\ze,d u)- \wt\na_{\ze} \nu+B(\nu) \ze, 
\endaligned
 \eEq 
 where 
 $$
 B(\nu)=\frac{1}{4}\big(J\nabla_{\nu}J + \nabla_{J\nu} J\big)
 $$
and 
$$
(\wh\nabla\ze)^{(0,1)} = \frac{1}{2}(\wh\nabla \ze + J \wh\nabla \ze \circ \mfj )
$$
is the $(0,1)$-part of $\wh\nabla\ze$, see \cite[Ch~3.1]{MS2}. The last term in (\ref{CRLNu_e}) is zero if $J$ is $\om$-compatible.\\

\noindent
In the relative Gromov-Witten theory (i.e. when $D$ is smooth), the most general almost complex structure that one may consider is an  $\om$-tame or compatible almost complex structures $J$ on $X$ that preserves $TD$ and satisfies the Nijenhueis condition (\ref{intInnormal_e}) along $D$. For such an $\om$-tame almost complex structure $J$, the most general perturbation term $\nu$ in $\cH_{g,k}(X,J)$ that one may consider is one satisfying
\bEq{RelNu_e}
\aligned
 &\nu|_{\ov{\mf{U}}_{g,k}\times D} \in \cH_{g,k}(D,J|_{TD})\quad\tn{and}\\
  & \frac{1}{2}\big(J\nabla_\nu J+\nabla_{J\nu} J \big) w -\big(\wt\nabla_w\nu + J\wt\nabla_{Jw}\nu\big)\in  \Om^{0,1}_{g,k}\otimes_\C T_xD,\quad \forall x\!\in\!D,~w\!\in\!T_xX;
  \endaligned
\eEq
see \cite{IP1,FZ1,FZ2}. The first parenthesis in the second line of (\ref{RelNu_e}) is zero if $J$ is $\om$-compatible. 
If the image of $u$ lies in $D$, by the first condition in (\ref{RelNu_e}), $\tn{D}_u \{\dbar-\nu\}$ induces an  operator 
\bEq{CRLNuNormal_e}
\tn{D}^\cN_u \{\dbar-\nu\}\colon \Gamma(\Si,u^*\cN_XD)\lra \Gamma(\Si,\Om^{0,1}_{\Si,\mfj}\otimes_\C u^*\cN_XD).
 \eEq 
 The second condition in (\ref{RelNu_e}) together with (\ref{intInnormal_e}) imply that (\ref{CRLNuNormal_e}) is $\C$-linear (i.e. it is a $\dbar$-operator); see \cite{IP1} and \cite[(2.14)]{FZ2}. Furthermore, 
\bEq{SimplifiedNDu}
\tn{D}^{\cN}_u\{\dbar-\nu\} \big([\ze]\big)\!=\! \bigg[(\wh\nabla \ze)^{(0,1)}- \frac{1}{2}\big(\wt\na_{\ze} \nu -J \wt\na_{J\ze}\nu\big)\bigg]^{\cN}\qquad \forall~[\ze]^{\cN}\in \Gamma(\Si,u^*\cN_XD),
\eEq
where $[\ze]^{\cN}$ denotes the image in $\Gamma(\Si,u^*\cN_XD)$ of a section $\ze\!\in\!\Gamma(\Si,u^*TX)$. If $\nu\!\equiv\!0$,  $\tn{D}^{\cN}_u\dbar$ is simply $u^*\dbar_{\cN_XD}$ where
$$
\dbar_{\cN_XD}\colon\Gamma(D,\cN_XD)\lra \Gamma(D,\Om^{0,1}_{D,J}\otimes_\C\cN_XD)
$$
is the $\dbar$-operator associated to $J$ in \cite[(2.7)]{FT1}.\\

\noindent
Generalizing (\ref{RelNu_e}) to the SNC case, the most general pairs $(J,\nu)$ that one may consider are those satisfying (\ref{intInnormal_e}) and (\ref{RelNu_e}) along each $D_i$. In other words, we need $\nu$ to satisfy
\bEq{RelNu_e2}
\aligned
 &\nu|_{\ov{\mf{U}}_{g,k}\times D_I} \in \cH_{g,k}(D_I,J_I)\quad\forall~I\!\subset\![N]\quad\tn{and}\\
  & \frac{1}{2}\big(J\nabla_\nu J+\nabla_{J\nu} J \big) w- \big(\wt\nabla_w\nu + J\wt\nabla_{Jw}\nu\big) \in  \Om^{0,1}_{g,k}\otimes_\C T_xD_I,\quad \forall x\!\in\!D_I,~w\!\in\!T_xX.
  \endaligned
\eEq
For $(J,\nu)$ satisfying the first condition of (\ref{RelNu_e2}), if $(\phi,u)$ is a $(J,\nu)$-holomorphic map as in (\ref{JnuMap_e}) with smooth domain, for every $i\!\in\![N]$, either $\tn{Im}(u)\!\subset\!D_i$ or $u$ intersects $D_i$ at finitely many points with positive tangency orders; see \cite[p.10]{FZ2}.  
Therefore, there exists a maximal subset $I\!\subset\![N]$ (called \textbf{depth} of $u$) such that $\tn{Im}(u)\!\subset\!D_I$ and $u$ intersects every $D_i$, with $i\!\notin\!I$, non-negatively. If the image of $u$ lies in $D_I$, $\tn{D}_u \{\dbar-\nu\}$ induces $\C$-linear $\dbar$-operators
\bEq{CRLNuNormal_ei}
\tn{D}^{\cN_i}_u \{\dbar-\nu\}\colon \Gamma(\Si,u^*\cN_XD_i)\lra \Gamma(\Si,\Om^{0,1}_{\Si,\mfj}\otimes_\C u^*\cN_XD_i)\qquad \forall~i\!\in\!I.
 \eEq 
Meromorphic section defined with respect to these $\dbar$-operators will be used to define log $(J,\nu)$-holomorphic curves. \\

\noindent
Working with the conditions in (\ref{RelNu_e2}) is hard.  \textit{In the following, instead of imposing the conditions (\ref{RelNu_e2}) on $\nu$, we define a class of logarithmic perturbation terms $\nu_{\log}$, associated to each of which we get a classical perturbation term $\nu$ satisfying  (\ref{RelNu_e2}).}\\

\noindent
Let $D$ be an SNC symplectic divisor in $(X,\om)$, $\cR$ be an $\om$-regularization, and $J$ be an $\om$-tame and $\cR$-compatible almost complex structure $J$ on $X$. In Section~\ref{Generic_ss}, we will define a class of smooth maps $u\colon\!\Si\!\lra\!X$ containing representatives of $\cM_{g,\mfs}(X,D,A)$ for which $\dbar u$ lifts to a log CR section 
\bEq{dbarLog_e}
\dbar^{\log}u\in\! \Gamma(\Si,\Om^{0,1}_{\Si,\mfj}\otimes_\C u^*TX(-\log D)).
\eEq 
In comparison with (\ref{cHgk_e}), (\ref{dbarLog_e}) indicates that the
right set of perturbation terms for the log moduli spaces is (a subspace of)
\bEq{cHgk_e2}
\bigg\{ \nu_{\log}\! \in\! \Gamma\big(\ov{\mf{U}}_{g,k}^\star\times X, \pi_1^*\Om^{0,1}_{g,k}\otimes_\C \pi_2^*TX(-\log D)_{\cR}\big)~~\tn{s.t.}~~\tn{supp}(\nu_{\log})\!\subset\big(\ov{\mf{U}}_{g,k}^\star-\bigcup_{a\in [k]} \tn{Im}(\mf{z}_a)\big)\!\times\! X \bigg\}.
\eEq
Associated to each $\nu_{\log}$ we get a classical perturbation term 
\bEq{AsscNu_e}
\nu\!=\!\iota(\nu_{\log})\!\in\! \cH_{g,k}(X,J),
\eEq
where by abuse of notation $\iota$ denotes the $\C$-linear homomorphism  
$$
\pi_1^*(\Om^{0,1}_{g,k})\otimes_\C \pi_2^*TX(-\log D)_{\cR}\lra \pi_1^*(\Om^{0,1}_{g,k})\otimes_\C \pi_2^*TX
$$
induced by (\ref{LogToTX_e}). Conversely, we may think of log perturbations as those $\nu$ that lift to a section of 
$$
\pi_1^*\Om^{0,1}_{g,k}\otimes_\C \pi_2^*TX(-\log D)_{\cR}.
$$ 
\bLm{NUlogtoNU_lm}
For every $\nu_{\log}$ in (\ref{cHgk_e2}), the associated classical perturbation term $\nu=\iota(\nu_{\log})$ satisfies (\ref{RelNu_e2}).
\eLm
\bPf
Given a logarithmic perturbation term $\nu_{\log}$ as in (\ref{cHgk_e2}), for each $I\!\subset\![N]$, restricted to the neighborhood $\Psi_I(\cN'_XD_I)\!\subset\!X$ of $D_I$ and with respect to the decomposition
$$
\Psi_I^*TX(-\log D)_\cR\cong \pi_I^*TD_I(-\log \partial D_I) \oplus \cN'_XD_I\times \C^I
$$
we get a decomposition
\bEq{tildeNu_e}
\nu_{\log,I}\equiv \Psi_I^*{\nu_{\log}}= \wt{\nu}_{I,\log} \oplus \wt{\theta}_I
\eEq
where 
$$
\wt{\theta}_I=(\wt{\theta}_{I,i})_{i\in I}\in \Gamma\big(\ov{\mf{U}}_{g,k}^\star\times \cN'_XD_I, \pi_1^*\Om^{0,1}_{g,k}\otimes_\C\C^I\big)
$$
is a tuple of $(0,1)$-forms. From (\ref{AsscNu_e}), (\ref{LogTrans_e}), and the decomposition (\ref{Tdecomp_e}) we get 
\bEq{NufromNuLog_e}
\nu_{I}\equiv \Psi_I^*{\nu}= \wt{\nu}_I \oplus n_I, 
\eEq
where 
\bEq{C*inv_e2}
\aligned
&n_I|_{v_I}(w)= \big(\wt\theta_I|_{v_I}(w)\big)\cdot v_I\equiv (\wt\theta_{I;i}|_{v_I}(w) v_{I;i})_{i\in I}\!\in\!\cN_XD_I|_{x}\\
&\forall~ x\!\in\!D_I,\quad v_I\!=\!(v_{I;i})_{i\in I}\!\in\!\cN'_XD_I|_{x},\quad w\in T_{g,k}.
\endaligned
\eEq
Note that $n_I|_{D_I}\!=\!0$; thus, $\nu$ satisfies the first condition in (\ref{RelNu_e2}). Furthermore, by (\ref{omStd_e}), (\ref{stdJ_e}), (\ref{C*inv_e2}), and the first condition in (\ref{RelNu_e2})
$$
J\nabla_\nu J w,\nabla_{J\nu} J w \in  T_xD_I,\quad \forall x\!\in\!D_I,~w\!\in\!T_xX,
$$
and, with $\theta_I=\wt\theta_I|_{D_I}$,  
$$
[\wt\nabla_w\nu]^{\cN} =\theta_I|_{x}\cdot w_I \in \Om^{0,1}_{g,k}\otimes_\C \cN_XD_I|_{x},\quad \forall x\!\in\!D_I,~w= w^{\tn{hor}}\!\oplus\!w_I \in T_xX\cong T_xD_I \oplus \cN_XD_I|_{x}.
$$ 
Therefore, $\nu$ also satisfies the second condition in (\ref{RelNu_e2}) and, by (\ref{SimplifiedNDu}), (\ref{CRLNuNormal_ei}) has the form
\bEq{SimplifiedNDu2}
\aligned
&\tn{D}_u^{\cN_i}\{\dbar-\nu\} (\ze_{I,i})=\! (\wh\nabla \ze_{I,i})^{(0,1)}- (\phi,u)^* \theta_{I,i}\cdot  \ze_{I,i}\qquad \forall~i\!\in\!I\\
&\forall~\ze_I\!=\!(\ze_{I,i})_{i\in I}\in \Gamma(\Si,u^*\cN_XD)\cong \bigoplus_{i\in I} \Gamma(\Si,u^*\cN_XD_i).
\endaligned
\eEq
\ePf

\bRm{Lognu_rmk}
Similarly to \cite{FT1}, log $(J,\nu)$-holomorphic curves can be defined for arbitrary $\nu_{\log}$ in (\ref{cHgk_e2}).
However, as in \cite[Rmk~1.3]{FT1}, certain steps in the proof of the compactness and the rest of analytical work are complicated for arbitrary such $\nu_{\log}$ (or arbitrary $J$ satisfying (\ref{intInnormal_e})). 
To avoid these complications, almost complex structures and perturbations considered in this paper are rather special. In the definition below, in a neighborhood of each stratum $D_I$, via the identification map $\Psi_I$, $\wt{\nu}_I$ and $\wt{\theta}_I$ are required to be pull back via $\pi_I$ of similar terms along $D_I$. 
\eRm

\bDf{stdJnu_dfn} 
Given an SNC symplectic divisor $D\!\equiv\!\bigcup_{i\in [N]} D_i$ in $(X,\om)$, an $\om$-regularization $\cR$, an $\om$-tame and $\cR$-compatible almost complex structure $J$ on $X$ (see Definition~\ref{stdJ_dfn}), and a perturbation term $\nu_{\log}$ in (\ref{cHgk_e2}), we say $\nu_{\log}$ is $\cR$-\textbf{compatible} if
\bEq{stdnu_e}
\wt{\nu}_{I,\log}=\pi_{I}^*(\nu_{I,\log})\quad \tn{and}\quad \wt{\theta_I}=\pi_I^*(\theta_I)\qquad \forall I\subset [N],
\eEq
for some 
$$
\nu_{I,\log}\!\in\! \Gamma\big(\ov{\mf{U}}_{g,k}^\star\times D_I, \pi_1^*\Om^{0,1}_{g,k}\otimes_\C \pi_2^*TD_I(-\log \partial D_I)_{\cR_I}\big)  
$$
and
\bEq{thetaI_e}
\theta_I=({\theta}_{I,i})_{i\in I}\in \Gamma\big(\ov{\mf{U}}_{g,k}^\star\times D_I, \pi_1^*\Om^{0,1}_{g,k}\otimes_\C\C^I\big).
\eEq
\eDf

\noindent
The condition (\ref{stdnu_e}) implies that
\bEq{C*inv_e}
n_I|_{\al v}(w)=\al~n_I|_v(w)\!\in\!\cN_XD_I|_{x}\qquad \forall~\al\!\in\!(\C)^I,\quad x\!\in\!D_I,\quad v,\al v\!\in\!\cN'_XD_I|_{x},\quad w\in T_{g,k}.
\eEq
\noindent
In other words, in a neighborhood of each stratum $D_I$ in $X$, identified with a neighborhood of $D_I$ in $\cN_XD_I$ via the regularization map $\Psi_I$ in $\cR$, the horizontal component of $\nu$ is the pull back from $D_I$ of some perturbation term $\nu_I$ on $D_I$, and the vertical component of $\nu$ is $(\C^*)^I$-equivariant with respect to the component-wise multiplicative action of $(\C^*)^I$ on $\cN_XD_I$. As explained in \cite[p.~11]{FZ2}, under these assumptions, $\Psi_{I}^*\nu$ extends to a $(\C^*)^I$-equivariant perturbation term over the fiber product of $\P^1$-bundles 
\bEq{PP_e}
\prod_{i\in I} \P(\cN_XD_i|_{D_I}\oplus \C).
\eEq

\vskip.1in
\noindent 
Define $\cH_{g,k}(X,D,\om)$ (resp. $\cH_{g,k}(X,D)$) to be set of such tuples $(\cR,J,\nu)$ (resp. $(\om,\cR,J)$) where $\nu$ is the perturbation term associated\footnote{Since the map $\nu_{\log}\lra \nu$ is one-to-one, it is safe to use $\nu$ in place of $\nu_{\log}$ to keep the notation simple.} to $\nu_{\log}$ as in Definition~\ref{stdJnu_dfn}.  For any fixed $(\om,\cR,J)$ (resp. $(\om,\cR)$), define $\cH_{g,k}(X,D)_{\cR,J}$ (resp. $\cH_{g,k}(X,D)_{\cR}$) to be the set of perturbation terms $\nu\!\in\!\cH_{g,k}(X,J)$ (resp. $(J,\nu)\!\in\!\cH_{g,k}(X,\om)$) such that $\nu_{\log}$ is compatible with a shrinking $\cR'$ of $\cR$. Similarly to (\ref{AKtoSymp_e}), the natural projection map
$$
\cH_{g,k}(X,D)\lra \tn{Symp}(X,D)
$$
is a weak homotopy equivalence. \\

\noindent
The first component $(\wh\nabla \ze_{I,i})^{(0,1)}$ in (\ref{SimplifiedNDu2}) is a $\dbar$-operator on $u^*\cN_XD_i$ by itself and is independent of $\nu$. In what follows, we will denote it by 
$$
\dbar_{u^*\cN_XD_i}\,\ze_{I,i}\;.
$$
The $\dbar$-operators  $\dbar_{u^*\cN_XD_i}$ and $D_{u}^{\cN_i} \{\dbar\!-\!\nu\}$ define (usually different) holomorphic structures on the pull-back complex line bundle $u^*\cN_XD_i$. The latter is a deformation of the former via the $(0,1)$-form $(\phi,u)^* \theta_{I,i}$. The one defined by $\dbar_{u^*\cN_XD_i}$ will be used as a reference. The perturbation caused by $\theta_{I,i}$ will allow us to achieve transversality.  \\

\noindent
By (\ref{stdJ_e}), (\ref{stdnu_e}), and (\ref{C*inv_e2}), the $X$-component $u\colon\!\Si\!\lra\!X$ of any $(J,\nu)$-holomorphic map $(\phi,u)$, with image in a sufficiently small neighborhood of $D_I$ (more precisely, in $\tn{Im}(\Psi_I)$) is determined by its projection 
\bEq{Proju_e}
u_I\colon \Si\!\lra\!D_I
\eEq
and a set of sections 
\bEq{zeI_e}
\ze_I\equiv(\ze_{I,i})_{i\in I} \in \bigoplus_{i\in I} \Gamma(u_I^*\cN_XD_i)
\eEq
such that 
\bEq{localJnumap_e}
\dbar u_I=\nu_I\qquad \tn{and}\qquad \dbar \ze_I\equiv  \bigoplus_{i\in I} \dbar_{u_I^*\cN_XD_i}\ze_{I,i}= n_I(u_I,\ze_I)= ((\phi,u_I)^*\theta_I) \cdot \ze_I\;.
\eEq
The second equation in (\ref{localJnumap_e}) can also be written in the compact form
\bEq{localJnumap_e2}
D_{u_I}^{\cN_i}\{\dbar\!-\!\nu\}(\ze_{I,i})=0\qquad\forall~i\!\in\!I.
\eEq
By (\ref{localJnumap_e}) or equally (\ref{localJnumap_e2}),  for any $c_I\!\equiv\!(c_i)_{i\in I}\!\in\!\C^I$, if $(u_I,\ze_I)$ satisfies (\ref{localJnumap_e}) then 
$$
\big(u_I,c_I\cdot \ze_I=(c_i\ze_{I,i})_{i\in I}\big)
$$ 
satisfies (\ref{localJnumap_e}) as well. Therefore, in a neighborhood of each stratum $D_I$ in $X$, identified with a neighborhood of $D_I$ in $\cN_XD_I$ via the regularization map $\Psi_I$ in $\cR$, the set of $(J,\nu)$-holomorphic maps is invariant under the component-wise multiplicative action of $\C^I$ on $\cN_XD_I$.
This explains the motivation behind the extra assumption (\ref{stdnu_e}). In general, as a family $\{u_t\}_{t\lra 0}$ of $(J,\nu)$-holomorphic maps sinks into $D_I$, the corresponding sections $\ze_{t,I}$ in  (\ref{zeI_e}) converge to zero. Then the idea is that, by rescaling $\ze_{t,I}$ we get a $(J_I,\nu_I)$-holomorphic map $u_I$ with image in $D_I$ and a similarly denoted holomorphic (or meromorphic) section 
$$
\ze_I= \lim_{t\lra 0} \frac{\ze_{t,I}}{|\ze_{t,I}|}\in u_I^*\cN_XD_I
$$ 
that remembers the direction at which the maps have approached $D_I$. Conversely, in the gluing construction, given $(u_I,\ze_I)$, gluing is done by pushing $u_I$ out in the direction of $\ze_I$. These local observations explain the motivation behind Definition~\ref{PreLogMap_dfn}  in the next section.

%--------------------------------------------------------
\subsection{Construction of perturbed moduli spaces}\label{LogJnu_ss}
\noindent
Let $(X,D)$ and $(\om,\cR,J,\nu)$ be as in the previous section, with $\nu$ coming from an $\cR$-compatible logarithmic perturbation term $\nu_{\log}$ as in Definition~\ref{stdJnu_dfn}. Suppose 
\bEq{Proju_e2}
u_I\colon \Si\!\lra\!D_I
\eEq
is a $(J_I,\nu_I)$-holomorphic map with smooth domain not mapped into $\partial D_I$ (so its depth is $I$) and 
\bEq{zeI_e2}
\ze_I\equiv(\ze_{I,i})_{i\in I} \in \Om_{\tn{mero}}\big(\bigoplus_{i\in I} u_I^*\cN_XD_i\big)
\eEq
is a tuple of non-zero meromorphic sections with respect to the holomorphic structure defined by 
$$
\bigg(D_{u}^{\cN_i}\{\dbar\!-\!\nu\}\bigg)_{i\in I}=    \bigg(\dbar_{u^*\cN_XD_i}-(\phi,u)^*\theta_{I,i}\bigg)_{i\in I} .
$$
In other words, $u\!\cong\!(u_I,\ze_I)$ is as in (\ref{Proju_e}) and (\ref{zeI_e}), but $\ze_{I,i}$ are allowed to have poles. For each $x\!\in\!\Si$, the \textbf{contact order vector}
\bEq{OrderVec_e}
\ord_x(u_I,\ze_I)=\big(\ord^i_x(u_I,\ze_I)\big)_{i\in [N]}\!\in\!\Z^N
\eEq
is defined by 
\bEq{OrderVec_e2}
\ord^i_x(u_I,\ze_I)=\ord_x(u,D_i)\!\geq \!0\quad  \forall~i\!\in\![N]\!-\!I\qquad\tn{and}\qquad \ord^i_x(u_I,\ze_I)=\ord_x\ze_{I,i}\quad  \forall~i\!\in\!I.
\eEq
The first item in (\ref{OrderVec_e2}) is the order of tangency of $u$ to $D_i$ at $x$; this is zero if $u(x)\!\notin\!D_i$ and is positive otherwise. The second item in (\ref{OrderVec_e2}) is the order of zero/pole of $\ze_{I,i}$ at $x$. In particular, it only depends on the $\C^*$-equivalence class $[\ze_{I,i}]$ of the section $\ze_i$; changing each $\ze_{I;i}$ with a non-zero constant multiple  of that does not change $\ord^i_x(u_I,\ze_I)$. 

\bRm{otherdbar_rmk}
For each $i\!\in\!I$, in a local holomorphic trivialization of $u_I^*\cN_XD_i$ around any point $x\!\in\!\Si$, with respect to the holomorphic structure defined by $\dbar_{u^*\cN_{X}D_i}$, the second equation in (\ref{localJnumap_e}) has the form 
\bEq{Localdbar}
\dbar f_i=\theta_{i} f_i
\eEq
where $f_i$ is a complex valued function and $\theta_{i}$ is a $(0,1)$-form. Let $g_i$ be any $\C$-valued function such that $\dbar g_i\!=\!\theta_{i}$. Then every solution of (\ref{Localdbar}) is of the form 
$$
f_i=\tn{e}^{g_i}h_{i}
$$ 
where $\dbar h_i\!\equiv\!0$. Thus, $\ord_x f_i$ in (\ref{OrderVec_e2}) can also be defined to be the order of meromorphic function $h_i$ at $x$ which is independent of the choice of $g_i$ and hence $\nu$. This shows that tangency order remains the same when we deform $\nu$.
\eRm

\noindent
In order to define nodal $(J,\nu)$-holomorphic maps, we use decorated dual graphs of the following sort associated to every $k$-marked nodal domain $(\Si,z_1,\ldots,z_k)$ as in \cite[Sec~3.1]{FT1}. Let $\Gamma\!=\!\Gamma(\V,\E,\L)$ be a graph with the set of vertices $\V$, edges $\E$, and legs $\L$; the latter, also called flags or roots, are half edges that have a vertex at one end and are open at the other end. Let $\uvec{\E}$ be the set of edges with an orientation. Given $\uvec{e}\!\in\!\uvec{\E}$, let $\ucev{e}$ denote the same edge with the opposite orientation. For each $\uvec{e}\!\in\!\uvec{\E}$, let $v_1(\uvec{e})$ and $v_2(\uvec{e})$ in $\V$ denote the starting and ending points of the oriented edge, respectively. For $v,v'\!\in\!\V$, let $\E_{v,v'}$ denote the subset of edges between the two vertices and $\uvec{\E}_{v,v'}$ denote the subset of oriented edges from $v$ to $v'$. For every $v\!\in\!\V$, let $\uvec{\E}_{v}$ denote the subset of oriented edges starting from $v$. Such decorated graphs $\Gamma$ characterize different topological types of nodal marked surfaces $(\Si,\vec{z}\!=\!(z_1,\ldots,z_k))$ in the following way.
Each vertex $v\!\in\!\V$ corresponds to a smooth component $\Si_v$ of $\Si$ with genus $g_v$. Each edge $e\!\in\!\E$ corresponds to a node $q_e$ obtained by connecting $\Sigma_v$ and $\Sigma_{v'}$ at the points $q_{\uvec{e}}\!\in\!\Si_v$ and $q_{\scz\ucev{e}}\!\in\!\Si_{v'}$, where $e\!\in\!\E_{v,v'}$ and $\uvec{e}$ is an orientation on $e$ with $v_1(\uvec{e})\!=\!v$. The last condition uniquely specifies $\uvec{e}$ unless $e$ is a loop connecting $v$ to itself. Finally, each leg $l\!\in\!\L$ connected to the vertex $v(l)$ corresponds to a marked point $z_{a(l)}\!\in \Sigma_{v(l)}$ disjoint from the connecting nodes. Thus we have
\bEq{nodalcurve_e}
(\Si,\vec{z})\! =\! \coprod_{v\in \V}(\Si_v,\vec{z}_v\cup {q}_v)/\sim, \quad  q_{\uvec{e}}\!\sim\! q_{\scz\ucev{e}}\quad\forall~e\!\in\!\E,
\eEq
\vspace{-.1in}
\noindent
where 
$$
\vec{z}_v\!=\!\vec{z}\cap \Sigma_v\qquad\tn{and}\quad \quad {q}_v=\{q_{\uvec{e}}\colon \uvec{e}\!\in\!\uvec{\E}_v\}\qquad \forall~v\!\in\!\V.
$$
We treat $q_v$ as an un-ordered set of marked points on $\Si_v$. We say $\Gamma$ is the \textbf{decorated dual graph} of $(\Si,\vec{z})$. Initially, each vertex $v\!\in\!\V$ is decorated by the genus $g_v\!\in\!\N$ of the corresponding irreducible component. Further decorations will be introduced below.
A complex structure $\mfj$ on $\Sigma$ is a set of complex structures $(\mfj_v)_{v\in \V}$ on its components. By a (complex) marked nodal curve, we mean a  marked nodal real surface together with a complex structure $(\Si,\mfj,\vec{z})$. Given a  map $u\colon \Si\lra X$, each vertex $v\!\in\!\V$ receives an additional decoration that is the homology class $A_v\!\in\!H_2(X,\Z)$ represented by $u_v=u|_{\Si_v}$.
Figure~\ref{labled-graph_fg}-Left illustrates a decorated graph with $2$ flags and Figure~\ref{labled-graph_fg}-Right is the associated marked nodal domain with $(g_1,\ldots,g_5)\!=\!(0,0,2,1,0)$.\\
\begin{figure}
\begin{pspicture}(14,1.3)(17,3.8)
\psset{unit=.3cm}

\pscircle*(55,6){.25}\pscircle*(59,6){.25}
\pscircle*(55,9){.25}\pscircle*(61.5,9){.25}\pscircle*(59,9){.25}
\psline[linewidth=.05](59,6)(55,9)
\psline[linewidth=.05](59,6)(59,9)\psline[linewidth=.05](59,6)(61.5,9)
\psarc[linewidth=.05](53,7.5){2.5}{-36.9}{36.9}\psarc[linewidth=.05](57,7.5){2.5}{143.1}{216.0}
\psline[linewidth=.05](61.5,9)(60.5,10.5)\psline[linewidth=.05](61.5,9)(62.5,10.5)
\rput(60.5,11.2){1}\rput(62.5,11.2){2}
\rput(55,5){$(g_4,A_4)$}\rput(59.5,5){$(g_5,A_5)$}
\rput(54,10){$(g_1,A_1)$}\rput(63.8,8.5){$(g_3,A_3)$}
\rput(58.3,10){$(g_2,A_2)$}
\end{pspicture}

\begin{pspicture}(4,1.7)(7,2)
\psset{unit=.3cm}

\pscircle*(52.4,8.2){.25}\pscircle*(52.4,11.8){.25}
\pscircle(50,10){3}\psellipse[linestyle=dashed,dash=1pt](50,10)(3,1)
\pscircle(54,7){2}\psellipse[linestyle=dashed,dash=1pt](54,7)(2,.77)
\pscircle*(55,8.6){.25}\rput(56,9){$z_1$}
\pscircle*(55,5.4){.25}\rput(56,5.4){$z_2$}

\pscircle(45,10){2}\psellipse[linestyle=dashed,dash=1pt](45,10)(2,.77)
\pscircle*(47,10){.25}\pscircle*(43.4,8.8){.25}\pscircle*(43.4,11.2){.25}
\psarc(41.8,7.6){2}{270}{71}
\psarc(41.8,12.4){2}{-71}{90}
\psarc(41.8,10){4.4}{90}{270}\psarc(42.8,10){.6}{120}{240}
\psarc(40,8.7){1.5}{60}{120}
\psarc(40,11.2){1.5}{225}{315}

\psarc(54,13){2}{60}{300}
\psarc(58,13){2}{240}{120}
\psarc(56,16.464){2}{240}{300}
\psarc(56,9.546){2}{60}{120}
\psarc(54,11.7){1.5}{60}{120}
\psarc(54,14.2){1.5}{225}{315}
\psarc(58,11.7){1.5}{60}{120}
\psarc(58,14.2){1.5}{225}{315}

\end{pspicture}
\caption{A nodal curve in $\ov\cM_{4,2}$ and its dual graph.}
\label{labled-graph_fg}
\end{figure}
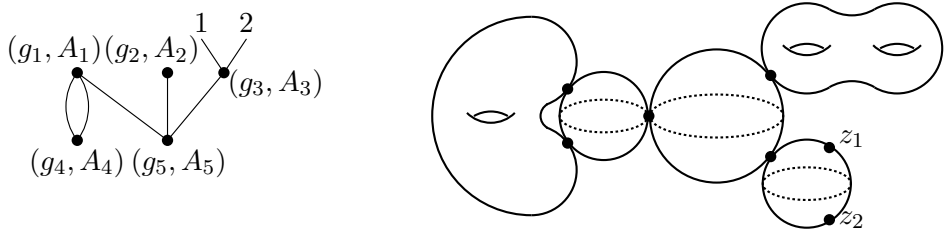 

\noindent
Assume $D\!=\!\bigcup_{i\in S} D_i\!\subset\!X$ is an SNC symplectic divisor, $(\om,\cR,J,\nu)\!\in\!\cH_{g,k}(X,D)$, and 
$$
u\!=\!(u_v)_{v\in \V}\colon\!(\Si,\mfj)\!\lra\!(X,J)
$$ 
is the $X$-component of a possibly nodal  $(J,\nu)$-holomorphic map $(\phi,u)$. In this situation, the dual graph of $(u,\Si)$ carries additional labelings
\bEq{LabelI_e}
I\colon \V,\E\lra \tn{subsets of}~[N], \qquad v\lra I_v\quad \forall v\!\in\!\V, \qquad e\lra I_e\quad \forall e\!\in\!\E,
\eEq
recording the minimal stratum $D_I$ that contains the image of $u_v$ and $u(q_e)$, respectively.

\bDf{PreLogMap_dfn} 
Suppose $D\!= \!\bigcup_{i\in [N]} D_i\!\subset\! X$ is an SNC symplectic divisor, $(\om,\cR,J,\nu)\!\in\!\cH_{g,k}(X,D)$, and 
$$
C\!\equiv\!(\Si,\mfj,\vec{z})=\bigg(\coprod_{v\in \V} C_v\equiv(\Si_v,\mfj_v,\vec{z}_v\cup q_v)\bigg)/\sim,\quad   q_{\uvec{e}}\!\sim\! q_{\scz\ucev{e}}\quad\forall~\uvec{e}\!\in\!\uvec{\E},
$$ 
is a connected nodal $k$-marked curve with smooth components $C_v$ and dual graph $\Gamma\!=\!\Gamma(\V,\E,\L)$ as in (\ref{nodalcurve_e}). A \textsf{pre-log $(J,\nu)$-holomorphic} map of contact type $\mfs\!=\!(s_a)_{a\in [k]}\!\in\!\big(\Z^N\big)^k$ from $C$ to $X$ is a collection 
\bEq{fplogSetUp_e}
f\equiv \Big(\phi,\big(f_{v}\!\equiv\! (u_v,\ze_v,C_v)\big)_{v\in \V}\Big)\;,
\eEq
such that  
\bEn
\item the tuple $\big( \phi,(u_v,C_v )_{v\in \V}\big)$ is a $k$-marked genus $g$ degree $A$ $(J,\nu)$-holomorphic map as in (\ref{JnuMap_e});
\item for each $v\!\in\!\V$, $(u_v,\ze_v\equiv (\ze_{v,i})_{i\in I_v})$ satisfies (\ref{Proju_e2})-(\ref{zeI_e2}) (with $I\!=\!I_v$, $u_I\!=\!u_v$, $\ze_I\!=\!\ze_v$);
\item for each $v\!\in\!\V$, $\ord_{u_v,\ze_v}$ is \textbf{supported} at the special points $z_v \cup q_v$ in the sense that 
$$
\ord_x(u_v,\ze_v)\neq0 \quad \Rightarrow \quad x\in z_v \cup q_v;
$$
\item\label{MatchingOrders_l} $s_{\uvec{e}}\!\equiv\! \ord_{q_{\uvec{e}}}(u_v,\ze_v)\!=\! -\ord_{q_{\scz\ucev{e}}}(u_{v'},\ze_{v'})\!\equiv\!-s_{\scz\ucev{e}}$ for all $v,v'\!\in\!\V$ and $\uvec{e}\!\in\! \uvec{\E}_{v,v'}$;
\item and $\ord_{z_a}(u_v,\ze_v)=s_a\in \Z^N$ for all $v\!\in\!\V$ and $z_a\!\in\!\vec{z}_v$.
\eEn
\eDf
\noindent
In other words, a pre-log map is simply a nodal $(J,\nu)$-map with a bunch of meromorphic sections on each smooth component (with zeros and poles only at the special points), dual contact orders at the nodes, and prescribed contact orders at the marked points.
Two $k$-marked pre-log $(J,\nu)$-holomorphic maps 
$$
f_{\log}=\big(\phi, (u_v,\Si_v,\mfj_v,\ze_v=(\ze_{v,i})_{i\in I_v}\big)_{v\in \V},\vec{z}\big)~~\tn{and}~~ f'_{\log}=\big(\phi',(u'_v,\Si'_v,\mfj'_v,,\ze'_v=(\ze'_{v,i})_{i\in I_v}\big)_{v\in \V},\vec{z}\,'\big)
$$ 
are equivalent if there exists a holomorphic reparametrization $h \colon\! \Si\!\lra\!\Si'$ such that $$
u'\circ h\!=\!u,\quad \phi\!=\!\phi'\circ h,\quad h(\vec{z})=\vec{z}\,',~~\tn{and}~~h_{v}^*\ze'_{h(v),i}\!=\!c_{v,i}\ze_{v,i}
$$ 
for all $v\!\in\!\V$ and some $c_{v,i}\!\in\!\C^*$. Here $v\!\lra\!h(v)$ is the induced map on the vertices of the dual graph. 
In particular, rescaling any of the meromorphic sections by a non-zero complex number does not change the equivalence class.
The space 
$$
\cM^{\tn{plog}}_{g,\mfs}(X,D,A,\nu)_{\Gamma}
$$
of the equivalence classes of stable pre-log maps of the fixed combinatorial type $(g,A,\mfs,\Gamma)$ is too big; see~\cite[Ex~3.6]{FT1}. Similarly to \cite{FT1}, in the following, we will put some restrictions on $\Gamma$ and take out a subspace that would give us a nice compactification with the correct expected dimension. \\

\noindent
In  \cite[Sec~5]{FT1}, corresponding to a decorated dual graph $\Gamma\!=\!\Gamma(\V,\E,\L)$ as in Definition~\ref{PreLogMap_dfn}  and an arbitrary orientation $O\!\equiv\!\{\uvec{e}\}_{e\in \E} \subset \uvec{\E}$ on the edges, we constructed a homomorphism of $\Z$-modules
\bEq{DtoT_e}
\D=\D(\Gamma)\equiv \Z^\E\oplus \bigoplus_{v\in \V} \Z^{I_v}  \stackrel{\vr=\vr_{O}}{\xrightarrow{\hspace*{1.5cm}}} 
\T=\T(\Gamma)\equiv\bigoplus_{e\in \E} \Z^{I_e}
\eEq
whose kernel $\K$ and cokernel $\CK$ are independent of the choice of the orientation $O$ on $\E$ and are invariants of the decorated graph $\Gamma$. For each $v\!\in\!\V$, $e\!\in\!\E$, and $s_v\!\in\! \Z^{I_v}\! \in\! \D$, the $e$-th component of $\vr(s_v)$ is equal to $s_v\!\in\! \Z^{I_e}$, if $v$ is the starting point of $\vec{e}$ and $e$ is not a loop; $\vr(s_v)$ is equal to $-s_v\!\in\!\Z^{I_e}$,  if $v$ is the ending point of $\vec{e}$ and $e$ is not a loop; and, is zero otherwise. 
In this definition and (\ref{Kspace_e}), via the identity $I_e\!=\!I_{v}\cup I_{v'}$ (see \cite[(3.11)]{FT1}) for all $e\in \E_{v,v'}$, and the inclusion 
$$
\Z^{I_v}\!\cong\!\Z^{I_v}\!\times\!\{0\}^{I_e-I_v}\!\subset\!\ \Z^{I_e},
$$ 
we can think of $s_{v}$ as a vector  also in $ \Z^{I_e}$. For each $e\!\in\!\E$, and $\la_e\!\in\! \Z_{e}\! \in\! \Z^{\E}\!\in\! \D$, the $e$-th component of $\vr(\la_e)$ is equal to $\la_es_{\uvec{e}}$, and the rest are zero.
In particular,
\bEq{Kspace_e}
 \K=\big\{\big((\la_e)_{e\in \E},(s_v)_{v\in \V}\big)\!\in\!\ \Z^\E\oplus \bigoplus_{v\in \V} \Z^{I_v} \colon s_{v'}\!-\!s_{v}\!=\!\la_es_{\uvec{e}}\quad\forall~v,v'\!\in\!\V,~\uvec{e}\!\in\!\uvec{\E}_{v,v'}\big\}.
\eEq

\noindent
In \cite[Lmm~3.7]{FT1}, to every (equivalence class of) pre-log $J$-holomorphic map $f$ we associate a group element
\bEq{group-element_e}
\tn{ob}_\Gamma(f) \in \mc{G}(\Gamma),
\eEq
where 
\bEq{Ggroup}
\mc{G}(\Gamma)=\prod_{e\in \E}(\C^*)^{I_e}/\tn{Im}(\tn{exp}(\vr_\C))
\eEq
is the complex torus with Lie algebra $\CK_\C=\coker(\vr_\C)$, and $\vr_\C$ is the natural extension of $\vr$ over $\C$.
The reasoning of the proof of \cite[Lmm~3.7]{FT1} also applies to the $(J,\nu)$-holomorphic maps and yields a group element (\ref{group-element_e}) for every pre-log $(J,\nu)$-holomorphic map. More explicitly,
given
\bEn
\item\label{coor_l} local holomorphic coordinates $w_{\uvec{e}}$ around\footnote{i.e., $w_{\uvec{e}}(q_{\uvec{e}})=0$.} each nodal point $q_{\uvec{e}}\!\in\!\Si_v$, for all $v\!\in\!\V$ and $\uvec{e}\!\in\!\uvec{\E}_v$,
\item\label{section_l} and representatives $\ze_v$ in the $\C^*$-equivalence class $[\ze_v]$, for all $v\!\in\!\V$, 
\eEn
for each $\uvec{e}\!\in\!\uvec{\E}_v$ and $i\!\in\!I_e\!-\!I_v$, locally around $q_{\uvec{e}}\!\in\!\Si_v$, by (\ref{Proju_e}) and (\ref{zeI_e}) we have 
\bEq{Expansion_e}
u_v(w_{\uvec{e}})\cong (u_{v,i}(w_{\uvec{e}}),\wt\eta_{\uvec{e},i}(w_{\uvec{e}})w_{\uvec{e}}^{s_{\uvec{e},i}})\in \cN_XD_i,
\eEq
such that 
\bEq{Eta1_e}
0\!\neq\!\eta_{\uvec{e},i}\equiv\wt\eta_{\uvec{e},i}(0)\!\in\! \cN_XD_i|_{u_v(q_{\uvec{e}})}.
\eEq
In other words, $\eta_{\uvec{e},i}$ is the $s_{\uvec{e},i}$-th order derivative of $u_v$ with respect to $w_{\uvec{e}}$ at $q_{\uvec{e}}$ in the normal direction to $D_i$.
Also, for each $\uvec{e}\!\in\!\uvec{\E}_v$ and $i\!\in\!I_v$, locally around $q_{\uvec{e}}\!\in\!\Si_v$, we have
\bEq{Expansion_e2}
\ze_{v,i}(w_{\uvec{e}})=\wt\eta_{\uvec{e},i}(w_{\uvec{e}})w_{\uvec{e}}^{s_{\uvec{e},i}}\in u_v^* \cN_XD_i,
\eEq
where 
\bEq{Eta2_e}
0\!\neq\!\eta_{\uvec{e},i}\equiv\wt\eta_{\uvec{e},i}(0)\!\in\! \cN_XD_i|_{u_v(q_{\uvec{e}})}.
\eEq
Then, by \cite[Lmm~3.7]{FT1}, the class $\tn{ob}_{\Gamma}(f)$ of 
\bEq{ObGDfn_e}
\prod_{e\in \E} \prod_{i\in I_{e}} \frac{\eta_{\uvec{e},i}}{\eta_{\scz\ucev{e},i}}\in \prod_{e\in \E}(\C^*)^{I_e}
\eEq
in $\mc{G}$ is independent of the choice such coordinates in \ref{coor_l} and representatives in \ref{section_l}.
In other words, $\tn{ob}_\Gamma(f)\!=\!1$ if and only if there exists such coordinates in \ref{coor_l} and representatives in \ref{section_l} such that
\bEq{CofG_e}
\eta_{\uvec{e},i}=\eta_{\scz\ucev{e},i} \qquad \forall~e\!\in\!\E,~i\!\in\!I_e.
\eEq
This condition will play a major role in the construction of gluing map in \cite{FT3}.

\bDf{LogMap_dfn}
Suppose $D\!= \!\bigcup_{i\in [N]} D_i\!\subset\! (X,\om)$ is an SNC symplectic divisor and $(\om,\cR,J,\nu)\!\in\!\cH_{g,k}(X,D)$. A \textbf{log $(J,\nu)$-holomorphic} map is a pre-log $(J,\nu)$-holomorphic map  $f$ with the decorated dual graph $\Gamma$ such that 
\bEn
\item\label{Tropical_l} there exist functions 
$$
s \colon\! \V\!\lra\!\R^N, \quad v\!\lra\!s_v,\qquad\tn{and}\qquad \la \colon\! \E\!\lra\!\R_+, \quad e\!\lra\!\la_e,
$$
such that 
\bEnalph
\item $s_v\!\in\!\R_{+}^{I_v}\!\times\!\{0\}^{[N]-I_v}$ for all $v\!\in\!\V$,
\item\label{Direction_l} $s_{v_{2}(\uvec{e})}\!-\!s_{v_1(\uvec{e})}\!=\!\la_{e} s_{\uvec{e}}$  for every $\uvec{e}\!\in\!\uvec{\E}$;
\eEnalph 
\item\label{GObs_e} and $\tn{ob}_\Gamma(f)\!=\!1\!\in\!\mc{G}(\Gamma)$.
\eEn
\eDf
\noindent
The first condition is a combinatorial condition on the decorated dual graph $\Gamma$ and is equal to the condition that the subspace 
$$
\K_\R\!\subset\!\R^\E\oplus \bigoplus_{v\in \V} \R^{I_v}
$$ 
defined as in (\ref{Kspace_e}) over $\R$ has a non-empty intersection with the positive quadrant. This implies that 
$$
\si\!=\!\K_\R\cap \big(\R_{\geq 0}^\E\oplus \bigoplus_{v\in \V} \R_{\geq 0}^{I_v}\big)\subset \K_\R
$$
is a maximal strictly convex rational polyhedral cone; see \cite[Lem~3.20]{FT1}. The affine toric variety associated to $\si$ is the space of gluing parameters (up to some multiplicity) that will be used in \cite{FT3}; see \cite[Sec~3.4]{FT1}. \\

\noindent
A marked log map is \textbf{stable} if it has a finite automorphism group. It is easy to see that a marked log map is stable if and only if the underlying $(J,\nu)$-map is stable. The equivalence class of a stable marked log map is called a \textbf{stable marked log curve}. We denote the space of stable $k$-marked degree $A$ genus $g$  log $(J,\nu)$-holomorphic curves of contact type $\mfs$ by 
$$
\ov\cM_{g,\mfs}^{\log}(X,D,A,\nu).
$$
Similarly to the $J$-holomorphic case, given $\mfs\!\in\!(\Z^N)^k$, for every $k$-marked stable nodal map $f$ representing an element of $\ov\cM_{g,k}(X,A,\nu)$ with dual graph $\Gamma$ and a choice of decorations $\{s_{\uvec{e}}\}_{\uvec{e} \in \uvec{\E}}$ satisfying the necessary combinatorial conditions
\bEq{CombCond_e}
s_{\uvec{e}}\!=\!-s_{\scz\ucev{e}} \quad \forall~e\!\in\!\E,\qquad  
\sum_{\uvec{e}\in \uvec{\E}_v} s_{\uvec{e}}+\sum_{l\in \L_v} s_l=(A_v\cdot D_i)_{i\in [N]}\quad \forall~v\!\in\!\V,
\eEq
and Definition~\ref{LogMap_dfn}.\ref{Tropical_l}, there exists at most one element $f_{\log}\!\in\!\ov\cM^{\log}_{g,\mfs}(X,D,A,\nu)$ lifting $f$ with this decorated dual graph. This is because every section $\ze_{v,i}$ is uniquely determined up to the action of $\C^*$ by the location and order of its zeros/poles. While $f_{\log}$ is stable if and only if $f$ is stable, the automorphism groups are sometimes different; see \cite[Ex~3.15-3.16]{FT1}. 

\begin{remark}
If $g\!=\!0$, a pre-log lift exists iff there are vectors $\{s_{\uvec{e}}\}_{\uvec{e}\in \uvec{\E}}$ satisfying the combinatorial condition (\ref{CombCond_e}) and the lift is unique. In other words, 
$$
\ov{\cM}^{\log}_{0,\mfs}(X,D,A,\nu)\subset \ov{\cM}_{0,k}(X,A,\nu)
$$
is the subset of $(J,\nu)$-curves $\big( \phi,(u_v,C_v )_{v\in \V}\big)$ for which there exist vectors $\{s_{\uvec{e}}\}_{\uvec{e}\in \uvec{\E}}$ satisfying (\ref{CombCond_e}), 
Definition~\ref{LogMap_dfn}\ref{Tropical_l}, and \ref{GObs_e}. The first two conditions are combinatorial. The latter is a condition on the derivatives of $u$ at the nodes which depends on $u$, $\nu$, and the configuration of special points on each component. In theory, it can be stated without mentioning the meromorphic sections.

\end{remark}

\bEx{g0A0_e}
Let $g\!=\!0$, $k\!\geq\!3$, $A\!=\!0\!\in\!H_2(X,\Z)$, and $\mfs\!=\!\vec{0}\!\equiv\! (0^N)^k\in (\N^N)^k$ (no perturbation). We show that 
$$
\ov\cM_{0,\vec{0}}(X,D,0)\cong \ov\cM_{0,k}\times X.
$$
Every element of $\ov\cM_{0,\vec{0}}(X,D,0)$ is (the equivalence class of) a $k$-marked nodal domain with the constant map into a point $p\!\in\!D_I$, for some maximal $I\!\subset\!X$, and $|I|$ meromorphic functions on each component. The claim is that all $k$-marked nodal domains are allowed and (up to equivalence) there is only one possibility for the meromorphic functions. The dual graph is a tree. Starting from a vertex $v\in\V$ with only one edge $\uvec{e}\!\in \uvec{\E}_v$, since all the $s_a$ are trivial, $s_{\uvec{e}}$ should be trivial as well. Removing $v$ and continuing inductively we conclude that all the vector decorations $s_{\uvec{e}}$ should be the trivial vector. Definition~\ref{LogMap_dfn}.\ref{Tropical_l} holds with $s_v\!\equiv\! s$ for all $v\!\in\!\V$ and any fixed $s\!\in\!\Z_+^I$. The meromorphic functions $\ze_{v,i}\colon \Si_v\!\lra\! \C$ can be taken to be constant $1$. Therefore, $\tn{ob}_\Gamma(f)\!=\!1$ (no matter what $\mc{G}$ is). Furthermore, $I_v=I_e=I$ for all $v\!\in\!\V$ and $e\!\in\!\E$, the map
$$
\vr\colon \Z^\E \oplus \bigoplus_{v\in \V}\Z^{I} \lra \bigoplus_{e\in \E}\Z^{I}
$$
in (\ref{DtoT_e}) is surjective, and 
$$
\tn{ker}(\vr)=\Z^\E \oplus \De,
$$
where $\De=\{(a)_{v\in \V}\in \bigoplus_{v\in \V}\Z^{I}:~a\!\in\! \Z^I\}\cong \Z^I$ is the diagonal (note that $|\V|-|\E|=1$). We conclude that  $\mc{G}$ is trivial.
In this example, transversality holds and perturbation is not needed.
\eEx

\bRm{TofCC_rmk}
In the previous example, one may replace $\mfs\!=\!\vec{0}$ with any $\mfs=(s_1,\cdots,s_k)$ satisfying $\sum_{a\in [k]} s_a=0^N$. 
Then,
 $$
\ov\cM_{0,\vec{0}}(X,D,0)\cong  \ov\cM_{0,k}\times D_{I_0},
$$
where $[N]\!-\!I_0$ is the maximal subset such that $s_{a,i}\!=\!0$ for all $a\!\in\! [k]$ and $i\!\in\![N]\!-\!I_0$.
For every $k$-marked genus $0$ nodal domain $(\Si,\vec{z})$, there is a unique set of decorations $\{s_{\uvec{e}}\}_{\uvec{e}\in \uvec{\E}}$ such that $\Si$ can be equipped with meromorphic functions that have zeros/poles of order $s_a$ at $z_a$. The group $\mc{G}$ is still trivial, but $\tn{ker}(\vr)$ will be different. 
\eRm

\bEx{gA0_e} 
Extending Example~\ref{g0A0_e} to the higher genus case, let $2g+ k\!\geq\!3$, $A\!=\!0\!\in\!H_2(X,\Z)$, and $\mfs\!=\!\vec{0}$ (and no perturbation yet). We show that still
$$
\ov\cM_{g,\vec{0}}(X,D,0)\cong  \ov\cM_{g,k}\times X,
$$
but transversality does not hold if $g\!>\!0$. The obstruction bundle is the rank $ng$ orbi-bundle
$$
\pi_1^* \mc{E}^*_g\otimes \pi_2^*TX(-\log D)\lra \ov\cM_{g,k} \times X,
$$
where $\pi_1,\pi_2$ are the projection maps onto the first and second components, respectively, and $\mc{E}_g\!\lra\! \ov\cM_{g,k}$ is the \textbf{Hodge bundle}.
Every element of $\ov\cM_{g,\vec{0}}(X,D,0)$ is (the equivalence class of) a $k$-marked nodal domain with the constant map into a point $p\!\in\!D_I$, for some maximal $I\!\subset\!X$, together with $|I|$ meromorphic functions on each component. If the dual graph is not a tree, a priori, there are infinitely many possibilities for the vector decorations $\{s_{\uvec{e}}\}_{\uvec{e}\in \uvec{\E}}$. The claim is that all $k$-marked nodal domains are allowed, but only the trivial decoration satisfies the conditions of Definition~\ref{LogMap_dfn}. By Definition~\ref{LogMap_dfn}.\ref{Tropical_l}, there should exists vectors $s_v\!\in\!\Z_+^I$ such that 
$$
s_{v_{2}(\uvec{e})}\!-\!s_{v_1(\uvec{e})}\!=\!\la_{e} s_{\uvec{e}}\quad\tn{for some}~\la_e\!>\!0,\qquad \forall~\uvec{e}\!\in\!\uvec{\E}.
$$
For each $i\!\in\!I$, choose $v\!\in\V$ such that $s_{v,i}\!\in\!\Z_+$ is maximal. By the previous identity, $s_{\uvec{e},i}\!\leq\! 0$ for all 
$\uvec{e}\!\in\!\uvec{\E}_v$. Since
$$
\sum_{\uvec{e}\in \uvec{\E}_v} s_{\uvec{e},i}=A_v\cdot D_i=0,
$$
we conclude that $s_{\uvec{e},i}\!=\! 0$ for all $\uvec{e}\!\in\!\uvec{\E}_v$. From this we conclude that all $s_v$ should be the same; therefore, $s_{\uvec{e}}\!=\! 0$ for all $\uvec{e}\!\in\!\uvec{\E}$. The meromorphic functions $\ze_{v,i}\colon \Si_v\!\lra\! \C$ can be taken to be constant~$1$. Therefore, $\tn{ob}_\Gamma(f)\!=\!1$ (no matter what $\mc{G}$ is). In this case, the map
$$
\vr\colon \Z^\E \oplus \bigoplus_{v\in \V}\Z^{I} \lra \bigoplus_{e\in \E}\Z^{I}
$$
in (\ref{DtoT_e}) is not necessarily surjective. In this example, transversality does not hold for two reasons. The logarithmic linearization of Cauchy Riemann operator in (\ref{Ddbarnu_e}) is not surjective and $1$ is not a regular value of $\tn{ob}_\Gamma$. Passing to a regular cover 
$$
p\colon \ov{\mf{M}}_{g,k}\lra\ov\cM_{g,k}
$$
as in (\ref{Cover_e}) and taking a generic perturbation term\footnote{This argument needs some justification, as we should explain the relation between such a perturbation with perturbations in (\ref{cHgk_e2}). } $\nu_{\log}$ in 
$$
 \Gamma\big(\ov{\mf{U}}_{g,k}\times X, \pi_1^* \mc{E}^*_g\otimes \pi_2^*TX(-\log D) \big),$$
we get 
$$
\ov\cM_{g,\vec{0}}(X,D,0,\nu)\cong \nu_{\log}^{-1}(0).
$$
Therefore, the Virtual Fundamental Class of $\ov\cM_{g,\vec{0}}(X,D,0)$ is the Euler class of the orbibundle $
\pi_1^* \mc{E}^*_g\otimes \pi_2^*TX(-\log D)$ in  $\ov\cM_{g,k} \times X$ (generalizing the classical and relative examples in \cite[Sec~4.1]{FZ2}).

\eEx

\noindent
Each moduli space $\ov\cM^{\log}_{g,\mfs}(X,D,A,\nu)$ is coarsely stratified by the subspaces
$$
\cM_{g,\mfs}(X,D,A,\nu)_{\Gamma}=\tn{ob}_{\Gamma}^{-1}(1)\subset \cM_{g,\mfs}^{\tn{plog}}(X,D,A,\nu)_{\Gamma}
$$
consisting of log $(J,\nu)$-holomorphic curves with the decorated dual graph $\Gamma$. Here a decoration consists of genus and degree decorations on vertices, ordering of the marked points and $\mfs$, labelings by subsets of $[N]$ in (\ref{LabelI_e}), and vectors $\{s_{\uvec{e}}\}_{\uvec{e}\in \uvec{\E}}$ satisfying the combinatorial condition (\ref{CombCond_e}). The vectors $\{s_{\uvec{e}}\}_{\uvec{e}\in \uvec{\E}}$ are also required to satisfy Definition~\ref{LogMap_dfn}\ref{Tropical_l} but the vectors $\{s_v\}_{v\in \V}$ are not part of the decoration.
By Theorem~\ref{Compactness_th}, for each $(g,\mfs,A)$, the set of such decorated dual graphs $\Gamma$ is finite.

\bRm{SubCurve_rmk}
Suppose $f$ as in (\ref{fplogSetUp_e}) is a $(J,\nu)$-log curve in $\cM_{g,\mfs}(X,D,A,\nu)_{\Gamma}$ and $\Gamma'=(\V',\E',\L')$ is a connected subgraph of $\Gamma$. The new set of legs $\L'$ consists of those legs $l\!\in\! \L$ such that $v(l)\!\in\!\V'$, as well as those oriented edges $\uvec{e}\!\in\!\uvec{\E}$   such that $v_1(\uvec{e})\!\in\V'$ but $v_2(\uvec{e})\!\notin\V'$. Let $f'$ be the tuple as in (\ref{fplogSetUp_e}) obtained by restricting to the connected sub-nodal curve 
$$
\Si'=\bigcup_{v\in \V'} \Si_{v}.
$$
Then $f'$ is also a log $(J,\nu|_{\Si'})$-curve with 
$$
\mfs'=\Big(\big(s_l\big)_{l\in \L\colon v(l)\in \V'}, \big(s_{\uvec{e}}\big)_{\uvec{e}\in\uvec{\E}\colon v_1(\uvec{e})\in \V', v_2(\uvec{e})\notin \V'}\Big)\qquad\tn{and} \qquad A'=\sum_{v\in \V'} A_{v}.
$$
The combinatorial conditions of Definition~\ref{PreLogMap_dfn} and Definition~\ref{LogMap_dfn}\ref{Tropical_l}  are clearly satisfied. By definition, $\tn{ob}_{\Gamma'}(f')\!=\!1\!\in\!\mc{G}(\Gamma')$ iff there are local holomorphic coordinates around each nodal point of $\Si'$ and representatives $\ze_v$ in the $\C^*$-equivalence class $[\ze_v]$, for all $v\!\in\!\V'$, such that (\ref{CofG_e}) holds. The restriction to $f'$ of  such local holomorphic coordinates and representatives for $f$ has the required property.
\eRm

%-------------------------------------------------
%---------------------------------------Transversality
%---------------------------------------------------
\section{Deformation theory and transversality}\label{Deformation_s}
%--------------------------------------------------------
\subsection{Main stratum}\label{Generic_ss} 
\noindent
In this section, similarly to the classical case, we realize $\cM_{g,\mfs}(X,D,A,\nu)$ as the zero set of a section $\dbar^{\log}-\nu_{\log}$ of some infinite dimensional bundle over an appropriately defined configuration space. The linearization of this section is a logarithmic lift $D_{u}^{\log}\{\dbar-\nu\}$ of the classical linearization map $D_{u}\{\dbar\!-\!\nu\}$. Then, it follows from Sard-Smale theorem that $\cM_{g,\mfs}(X,D,A,\nu)$ is cut transversely for generic $\nu$. If $\nu\!\equiv\!0$, the same statement holds for generic $J$ if we restrict to the subspace of simple maps $\cM^\star_{g,\mfs}(X,D,A)$.\\

\noindent
We prove the following transversality statement.
\bPr{Trans1_prp}
Suppose $(X,\om)$ is a closed symplectic manifold, $D=\bigcup_{i\in [N]} D_i$ is an SNC symplectic divisor,  $A\!\in\!H_2(X,\Z)$, $g,k\!\in\!\N$,  and $\mfs\!\in\!(\N^N)^k$. 
\bEn
\item If $2g\!+\!k\!\geq\!3$, for any given choice of universal family in (\ref{UFamily_e}), there exists a Baire set of second category $\cH^{\eset}_{g,k}(X,D)_{\cR,J}\!\subset\!\cH_{g,k}(X,D)_{\cR,J}$ such that for each $\nu\!\in\!\cH^{\eset}_{g,k}(X,D)_{\cR,J}$, $\cM_{g,\mfs}(X,D,A,\nu)$ is a naturally-oriented smooth manifold of the real dimension 
$$
2\big(c_1^{TX(-\log D)}(A)+(n-3)(1-g)+k\big).
$$
The restriction of $\tn{st}\times\tn{ev}$ in (\ref{ev-st_e}) to $\cM_{g,\mfs}(X,D,A,\nu)$ is smooth.
\item  If\,\footnote{Restriction to $g\!=\!0$ is not necessary here but this is the case that we will need later. Furthermore, the generalization of this statement in part (2) of Proposition~\ref{Trans2_prp} requires the $g\!=\!0$ assumption or allowing deformations of $\cR$.} $g\!=\!0$, $\nu\!\equiv\!0$, and $A\!\neq\!0$ or $k\!\geq\! 3$, the same statement holds for $J$ in a Baire set of second category $\tn{AK}^{\eset}(X,D)_{\cR}\!\subset\!\tn{AK}(X,D)_{\cR}$ if we restrict to the subspace of simple maps $\cM^\star_{0,\mfs}(X,D,A)$.
\eEn
\ePr
\vskip.1in
\noindent
We start by setting up a suitable analytical frame work for studying the deformation theory of log $(J,\nu)$-holomorphic maps. This set up is in some sense the main step of the proof.
\bDf{CTs_dfn}
Fix a smooth $k$-marked genus $g$ curve $(\Si,\mfj,\vec{z})$,  local holomorphic coordinates\footnote{This is not needed here, but it will be needed in constructing a Banach completion of $\tn{Map}_{A,\mfs}\big((\Si,\vec{z}),(X,D)\big)$; see the end of Remark~\ref{Trivialization}.} around the marked points, $A\!\in\!H_2(X,\Z)$, a regularization $\cR$ for $D$ in $X$, and $\mfs$ as in (\ref{stuple_e}). 
With $I_a\!\subset\![N]$ as in (\ref{Ia_e}), for each $a\!\in\![k]$, we say a smooth map $u\colon\Si\!\lra\!X$ has \textbf{contact type $\mfs$ with $D$ at $\vec{z}$} if 
\bIt
\item $u^{-1}(D)=\mf{d}\subset \{z_1,\ldots,z_k\}$,
\item $u(z_a)=p_a\in D_{I_a}- \partial D_{I_a}$, for all $a\!\in\![k]$,
\item and 
\bEq{Locu_e}
\Psi_{I_a}^{-1}\circ u(w_a)= \big(u_a(w_a),\oplus_{i\in I_a}~w_a^{s_{ai}}\eta_{a,i}(w_a)\big)\in \cN_XD_{I_a},\quad~\forall a\!\in\![k],
\eEq
where $U_a$ is a sufficiently small neighborhood of the marked point $z_a$ in $\Si$, $w_a$ is the local holomorphic coordinate on $U_a$ with $w_a(z_a)\!=\!0$, $u_a\colon\!U_a\!\lra\!D_{I_a}$ is the projection of $\Psi_{I_a}^{-1}\circ u$ to $D_{I_a}$, and $\eta_a=(\eta_{a,i})_{i\in I_a}$ is a smooth section of 
$$
u_a^*\cN_XD_{I_a}\cong \bigoplus_{i\in I_a}u_a^*\cN_XD_{i}
$$
satisfying
$$
\eta_{a,i}(0)\neq 0, \qquad \forall~i\!\in\!I_a.
$$
\eIt
\eDf
\noindent
Let $\tn{Map}_{A,\mfs}\big((\Si,\vec{z}),(X,D)\big)$ denote the set of smooth degree $A$ maps of contact type $\mfs$ with $D$ at $\vec{z}$. This space is an infinite dimensional Fr\'echet manifold whose tangent space at any $u$ is the infinite dimensional vector space
$$
\Gamma(\Si,u^*TX(-\log D)_{\cR}).
$$ 
More explicitly, if $\{u_t\}_{t\in [0,\ep)}$ is a $1$-parameter family of maps in $\tn{Map}_{A,\mfs}\big((\Si,\vec{z}),(X,D)\big)$, restricted to $\Si-\mf{d}$, by the first bullet in Definition~\ref{CTs_dfn}, we get 
$$
\xi_\eset=\frac{\nd}{\nd t}u_t|_{t=0}\in \Gamma(\Si-\mf{d},\{u_0|_{\Si-\mfd}\}^*TX(-\log D))\cong\Gamma(\Si-\mfd,\{u_0|_{\Si-\mfd}\}^*TX).
$$
On the other hand, for each $a\!\in\![k]$ with $z_a\!\in\!\mfd$, restricted to the chart $U_a$ in the third bullet above, by (\ref{Locu_e}), we have 
 $$
 \Psi_{I_a}^{-1}\circ u_t= \Big(u_{t,a},\bigoplus_{i\in I_a}\;w_a^{s_{ai}}\eta_{t,a,i}\Big).
 $$
Therefore,
 \bEq{DofLog_e}
\frac{\nd}{\nd t} \big(\Psi_{I_a}^{-1}\circ u_t\big)|_{t=0}= \pi_{I_a}^* \xi_{I_a} \oplus \bigoplus_{i\in I_a} w_a^{s_{ai}}\,\eta_{0,a,i}\,c_{a,i},
 \eEq
 where $\pi^* \xi_{I_a}$ is the horizontal lift of 
$$
\xi_{I_a}= \frac{\nd}{\nd t}u_{t,a}|_{t=0}\in \Gamma(U_a,u_{0,a}^*TD_{I_a})
$$ 
to the horizontal subspace $T^{\tn{hor}}\cN_XD_{I_a}\cong \pi_{I_a}^* TD_{I_a}$, and
 $$
 c_{a,i}=\frac{\frac{\nd}{\nd t}\eta_{t,a,i}(w_a)|_{t=0}}{\eta_{0,a,i}(w_a)}\colon U_a\lra \C, \quad \forall~i\!\in\!I_a.
 $$
 By (\ref{TXlog_e})-(\ref{LogTrans_e}) and (\ref{DofLog_e}), $\xi_\eset$ and \{$\xi_{I_a} \oplus (c_{a,i})_{i\in I_a}\}_{z_a\in \mfd}$ define a global section $\xi_{\log}$ of $u^*TX(-\log D)$ that maps to 
 $$
 \xi=\frac{\nd}{\nd t}u_t|_{t=0}\!\in\! \Gamma(\Si,u_0^*TX)
 $$ 
 under the homomorphism $\iota$ in (\ref{LogToTX_e}). Conversely, given a section $\xi_{\log}\!\in\!\Gamma(\Si,u_0^*TX(-\log D))$, logarithmic exponentiation (see (\ref{Log-exp_e})) of $\xi_{\tn{log}}$, corresponding to a Hermitian metric on $TX(-\log D)$, produces a $1$-parameter family of maps in $\tn{Map}_{A,\mfs}\big((\Si,\vec{z}),(X,D)\big)$ with tangent vector $\xi\!=\!\iota(\xi_{\log})$ at $t\!=\!0$. 
 
 \noindent
 \bLm{dulog_lmm}
 For every $u\!\in\!\tn{Map}_{A,\mfs}\big((\Si,\vec{z}),(X,D)\big)$, there exist a \textit{logarithmic Cauchy-Riemann section} 
 \bEq{dbarlog_eq}
 \dbar^{\tn{log}} u\in \Gamma(\Si,\Om^{0,1}_{\Si,\mfj}\otimes_\C u^*TX(-\log D)) 
 \eEq
 such that the following diagram commutes:
 \bEq{dulog_e22}
\xymatrix{
& && u^*TX(-\log D)\ar[d]^{\iota}\\
& T^{0,1}\Si \ar[rru]^{\dbar^{\log} u}\ar[rr]^{\dbar u} && u^*TX\;.
}
\eEq
 \eLm 
 
 \bPf
 Away from $\mfd$, by the first bullet in Definition~\ref{CTs_dfn} and the identification 
 $$
 TX(-\log D)|_{X-D}\cong TX|_{X-D},
 $$
 we define $\dbar^{\log} u=\dbar u$. For each $a\!\in\![k]$ with $z_a\!\in\!\mfd$, restricted to the chart $U_a$ in Definition~\ref{CTs_dfn}, by (\ref{Locu_e}), we have 
 $$
 \Psi_{I_a}^{-1}\circ \dbar u= \pi^* \dbar u_a\oplus \bigoplus_{i\in I_a} \dbar_{u^*\cN_XD_i} (w_a^{s_{ai}}\eta_{a,i})=\pi_{I_a}^* \dbar u_a\oplus \bigoplus_{i\in I_a}  w_a^{s_{ai}}\dbar_{u^*\cN_XD_i} (\eta_{a,i}).
 $$
Restricted to $U_a$, we define $\dbar^{\log} u$ to be
 $$
 \Psi_{I_a}^{-1}\circ \dbar^{\log} u= \pi_{I_a}^*\dbar u_a\oplus \bigoplus_{i\in I_a}   \frac{\dbar_{u^*\cN_XD_i} \eta_{a,i}}{\eta_{a,i}}\;.
 $$
By (\ref{TXlog_e})-(\ref{LogTrans_e}), these local sections define a global section (\ref{dbarlog_eq}) that maps to $\dbar u$ under the homomorphism $\iota$ in (\ref{LogToTX_e}). 
 \ePf
 
 \noindent
Let $\tn{Map}_{A}\big(\Si,X\big)$ denote the space of all smooth degree $A$ maps from $\Si$ into $X$, 
\bEq{ClassicE_e}
\mc{E}_{A}\big(\Si,X\big)\lra
\tn{Map}_{A}\big(\Si,X) 
\eEq
denote the infinite dimensional vector space whose fiber over every map $u$ is  
$$
\Gamma(\Si,\Om^{0,1}_{\Si,\mfj}\otimes_\C u^*TX),
$$
and
\bEq{LogE_e}
\mc{E}_{A,\mfs}\big((\Si,\vec{z}),(X,D)\big)\lra
\tn{Map}_{A,\mfs}\big((\Si,\vec{z}),(X,D)\big) 
\eEq
denote the infinite dimensional vector space whose fiber over every map $u$ is  
$$
\Gamma(\Si,\Om^{0,1}_{\Si,\mfj}\otimes_\C u^*TX(-\log D)).
$$
The classical CR operator (\ref{Jnu-holo_e}) can be seen as a section of (\ref{ClassicE_e}). By Lemma~\ref{dulog_lmm}, the restriction of this section to $\tn{Map}_{A,\mfs}\big((\Si,\vec{z}),(X,D)\big)$ defines a section of  (\ref{LogE_e}). 
Similarly, for every $\nu\!\in\!\cH_{g,k}(X,D)_{\cR,J}$ associated to $\nu_{\log}$ as in (\ref{AsscNu_e}) and an identification $\phi$ of $(\Si,\mfj,\vec{z})$ with a fiber of the universal family $\pi\colon\ov{\mf{U}}_{g,k}\!\lra\!\ov{\mf{M}}_{g,k}$ (used to define $\nu$), the restriction of $\dbar-\nu$ to $\tn{Map}_{A,\mfs}\big((\Si,\vec{z}),(X,D)\big)$ lifts to the section $\dbar^{\log}-\nu_{\log}$ of  (\ref{LogE_e}) so that the following diagram commutes:
$$
\xymatrix{
& \mc{E}_{A,\mfs}\big((\Si,\vec{z}),(X,D)\big)\ar[rr]&& \mc{E}_{A}\big(\Si,X\big)\\
& \tn{Map}_{A,\mfs}\big((\Si,\vec{z}),(X,D)\big)\ar[u]^{\dbar^{\log}-\nu_{\log}} \ar[rr]&& \tn{Map}_{A}\big(\Si,X\big)\ar[u]^{\dbar-\nu} .
}
$$
The linearization of $\dbar^{\log}-\nu_{\log}$ along the zero set is then the restriction/lift 
$$
\tn{D}^{\tn{log}}_u\{\dbar\!-\!\nu\}\equiv \tn{D}_u\{\dbar^{\log}\!-\!\nu_{\log}\}
$$ 
of the classical linearization map   $\tn{D}_u\{\dbar-\nu\}$ to $\Gamma(\Si,u^*TX(-\log D))$ so that the following diagram commutes:
\bEq{Ddbarnu_e}
\xymatrix{
& \Gamma(\Si,u^*TX(-\log D)) \ar[rrr]^{\tn{D}^{\tn{log}}_u\{\dbar-\nu\}} \ar[d]^{\iota_1}&&& \Gamma(\Si,\Om^{0,1}_{\Si,\mfj}\otimes_\C u^*TX(-\log D))\ar[d]^{\iota_2}\\
& \Gamma(\Si,u^*TX)  \ar[rrr]^{\tn{D}_u\{\dbar-\nu\}} &&& \Gamma(\Si,\Om^{0,1}_{\Si,\mfj}\otimes_\C u^*TX)\;.
}
\eEq

\noindent
Fix a $\mfj$-Hermitian metric on $T\Si$ and a $J$-Hermitian metric on $TX(-\log D)$, an integer $\ell\!\geq\!1$, and a real number $p\!>\!2$. Via the logarithmic exponentiation map, we can then construct a completion  $W^{\ell,p}_{A,\mfs}\big((\Si,\vec{z}),(X,D)\big)$ of $\tn{Map}_{A,\mfs}\big((\Si,\vec{z}),(X,D)\big)$  which is a smooth separable Banach manifold with tangent space 
\bEq{T_u-e}
T_uW^{\ell,p}_{A,\mfs}\big((\Si,\vec{z}),(X,D)\big)=W^{k,p}\big(\Si,u^*TX(-\log D)\big).
\eEq
The completion $\cE^{\ell-1,p}_{A,\mfs}\big((\Si,\vec{z}),(X,D)\big)$ of $\cE_{A,\mfs}\big((\Si,\vec{z}),(X,D)\big)$ is a Banach complex vector bundle over $W^{\ell,p}_{A,\mfs}\big((\Si,\vec{z}),(X,D)\big)$. For every  $\nu\!\in\!\cH_{g,k}(X,D)_{\cR,J}$ and any identification $\phi$ of $(\Si,\mfj,\vec{z})$ with a smooth fiber of the universal family $\ov{\mf{U}}_{g,k}$ used to define $\nu$, $\dbar^{\log}\!-\!\nu_{\log}$ defines a smooth section of Banach bundle
\bEq{BVB_e}
\cE^{\ell-1,p}_{A,\mfs}\big((\Si,\vec{z}),(X,D)\big)\lra W^{\ell,p}_{A,\mfs}\big((\Si,\vec{z}),(X,D)\big).
\eEq

 \noindent
As in the classical case (\ref{CRLNu_e}),  $\tn{D}^{\tn{log}}_u\{\dbar-\nu\}$ can be written as the sum of a complex linear map (a~CR operator on $u^*TX(-\log D)$) and a compact operator (It is the restriction of the corresponding operators). Thus, by Riemann-Roch, it is a Fredholm operator with index
$$
\tn{dim}_\R~\Def_{\log}(u)-\tn{dim}_\R~\Obs_{\log}(u)=2\big( \tn{deg}(u^*TX(-\log D)) \!+\! \dim_\C\! X (1\!-\!g)\big),
$$
where 
$$
\Def_{\log}(u)\equiv \tn{ker}\big(\tn{D}^{\tn{log}}_u\{\dbar-\nu\}\big)\quad \tn{and}\quad 
\Obs_{\log}(u)\equiv \tn{coker}\big(\tn{D}^{\tn{log}}_u\{\dbar-\nu\}\big).
$$
From Implicit Function Theorem \cite[Thm~A.3.3]{MS2}, we deduce the following corollary.
\bCr{TrCase_cr}
If $u\!\in\!\tn{Map}_{A,\mfs}\big((\Si,\vec{z}),(X,D)\big)$ is $(J,\nu)$-holomorphic and $\Obs_{\log}(u)\!\equiv\! 0$, in a small neighborhood $B(u)$ of $u$ in $W^{\ell,p}_{A,\mfs}\big((\Si,\vec{z}),(X,D)\big)$ the set of $(J,\nu)$-holomorphic maps 
$$
V_u \!\equiv\! \{\dbar^{\log}-\nu_{\log}\}^{-1}(0)\cap B(u)
$$ 
is a smooth manifold of real dimension (\ref{RR_e}).
\eCr
\noindent 
Furthermore, by elliptic regularity, all the elements of $V_u$ and $\Def_{\log}(u)\cong T_u V_u$ are smooth (see \cite[Thm~3.1.5]{MS2}). 
The manifold $V_u$ carries a \textit{natural orientation}. Starting with the complex linear part of  $\tn{D}^{\tn{log}}_u\{\dbar-\nu\}$, both the kernel and cokernel of that are complex linear and thus naturally oriented. By deforming $\tn{D}^{\tn{log}}_u\{\dbar-\nu\}$ into its complex linear part via a  $1$-parameter family of compact operators, \cite[Prp A.2.4]{MS2} gives us a natural orientation on $\Def_{\log}(u)$.\\

\noindent
Next, we consider deformations of the marked domain $C\!=\!(\Si,\mfj,\vec{z})$.
Given a regular covering $p\colon \ov{\mf{M}}_{g,k}\!\lra\!\ov\cM_{g,k}$ and a universal family $\pi\colon\ov{\mf{U}}_{g,k}\!\lra\!\ov{\mf{M}}_{g,k}$ as in (\ref{Cover_e})-(\ref{UFamily_e}), let 
$$
\mf{M}_{g,k}\!=\!p^{-1}(\cM_{g,k}), \quad \mf{U}_{g,k}\!=\!\pi^{-1}(\mf{M}_{g,k}),\quad \pi\colon\mf{U}_{g,k}\!\lra\!\mf{M}_{g,k},
$$
be the restrictions to the subspace of smooth curves.
Choose a projective embedding  $\ov{\mf{U}}_{g,k}\!\lra\!\P^M$, for some sufficiently large $M$, and define
\bEq{UBS_e}
W^{\ell,p}_{A,\mfs}\big((\mf{U}_{g,k},\vec{\mf{z}}),(X,D)\big)\equiv\Big\{(c,u)\colon c\!\in\!\mf{M}_{g,k},~u\!\in\!W^{\ell,p}_{A,\mfs}\big((\pi^{-1}(c),\vec{\mf{z}}(c)),(X,D)\big)\Big\}
\eEq
where the metric considered on 
$$
\big((\Si,\mfj_c),\vec{z}\big)\!\equiv\!\big(\pi^{-1}(c),\vec{\mf{z}}(c)\big),
$$ 
for each $c\!\in\!\mf{M}_{g,k}$, is the restriction of Fubini-Study metric on $\P^M$ to the image of $\Si$.

\bRm{Trivialization}
By definition, each fiber of the projection map 
\bEq{SB_e}
W^{\ell,p}_{A,\mfs}\big((\mf{U}_{g,k},\vec{\mf{z}}),(X,D)\big)\lra \mf{M}_{g,k}
\eEq
has a Banach manifold structure but the total space does not a priori come with a natural smooth structure. In order to define a Banach manifold structure on (\ref{UBS_e}), we need to fix a smooth trivialization of $\mf{U}_{g,k}\!\lra\! \mf{M}_{g,k}$. Such a trivialization gives us a trivialization of $(\ref{UBS_e})$ and thus a product Banach manifold structure on that; see \cite[Sec~6.1]{FF}. The genus-$g$ surface-bundle $\mf{U}_{g,k}\!\lra\! \mf{M}_{g,k}$ does not necessarily admit a global smooth trivialization. However, for $\cB\!\subset\! \mf{M}_{g,k}$ sufficiently small around any $b\!\in\! \mf{M}_{g,k}$, 
\bEq{LocalUnivFam_e}
\cC\!=\!\pi^{-1}(\cB)\lra\cB
\eEq
is smoothly trivial. In other words, locally around every $\big((\Si,\mfj_b),\vec{z}\big)\!\equiv\!\big(\pi^{-1}(b),\vec{\mf{z}}(b)\big)$ there exists an $\tn{Aut}\big((\Si,\mfj_b),\vec{z}\big)$-equivariant diffeomorphism
\bEq{ST_e}
\varphi\colon \cC \lra \Si\times \cB 
\eEq
such that $\pi\circ \varphi^{-1}$ is the projection onto the second factor, each section $\varphi \circ \mf{z}_a$ is constant, and $\varphi|_{\pi^{-1}(b)}\!=\! \tn{id}_\Si$. The smooth trivialization $\varphi$ gives rise to a Banach manifold structure on the restriction
\bEq{UBSC_e}
W^{\ell,p}_{A,\mfs}\big((\mc{C},\vec{\mf{z}}|_{\cB}),(X,D)\big)\equiv\Big\{(c,u)\colon c\!\in\!\cB,~u\!\in\!W^{\ell,p}_{A,\mfs}\big((\pi^{-1}(c),\vec{\mf{z}}(c)),(X,D)\big)\Big\}
\eEq
which we denote by $W^{\ell,p}_{A,\mfs}\big((\mc{C},\vec{\mf{z}}|_{\cB}),(X,D)\big)_{\varphi}$.
If $\varphi_1$ and $\varphi_2$ are two such smooth trivialization maps, the map 
$$
W^{\ell,p}_{A,\mfs}\big((\mc{C},\vec{\mf{z}}|_{\cB}),(X,D)\big)_{\varphi_1}\lra W^{\ell,p}_{A,\mfs}\big((\mc{C},\vec{\mf{z}}|_{\cB}),(X,D)\big)_{\varphi_2}
$$
induced by the change of trivialization map $\varphi_2\!\circ\!\varphi_1^{-1}$ is not smooth (unless $\varphi_2\!\circ\!\varphi_1^{-1}$ is constant in $c$); see  \cite[Sec~3.1]{MW2}. Therefore, there is no natural way of putting a Banach manifold structure on (\ref{SB_e}). On the other hand, restricted to the moduli space which is the zero set of $\dbar^{\log}\!-\!\nu_{\log}$, (by elliptic regularity) the transition maps are smooth. Therefore, in the proof of Proposition~\ref{Trans1_prp} below, we cover $\mf{M}_{g,k}$ with countably many such charts, find a Baire set of regular perturbation terms for each one, and then take intersection which yields a Baire set again. It is for the same reason that we fix local coordinates around the marked points in Definition~\ref{CTs_dfn}. If $w_a$ and $w_a'$ are two local coordinates around the marked point $z_a$, they are related by a $\C^*$-valued reparametrization map $\varphi$; i.e. $w_a'=\varphi(w_a)w_a$. The Banach smooth structure $W^{\ell,p}_{A,\mfs}\big((\Si,\vec{z}),(X,D)\big)$ with respect to $w_a$ and $w_a'$ will not be the same unless $\varphi$ is constant.  In the family version (\ref{LocalUnivFam_e}), we will fix local defining equations $w_a\colon \cC\lra\C$ for Cartier divisors $\mf{z}_a(\cB)\!\subset\!\cC$. The restriction of $w_a$ to each fiber is the local coordinate needed in Definition~\ref{CTs_dfn}; see \cite[(5.20)-(5.21)]{FF} for details. 
\eRm

\noindent
Let 
\bEq{BVB_e2}
\cE^{\ell-1,p}_{A,\mfs}\big((\mf{U}_{g,k},\vec{\mf{z}}),(X,D)\big)\lra W^{\ell,p}_{A,\mfs}\big((\mf{U}_{g,k},\vec{\mf{z}}),(X,D)\big)
\eEq
denote the natural extension of (\ref{BVB_e}) over $\mf{M}_{g,k}$.
Similarly to the previous paragraph, locally over any sufficienlty small neighborhood $\mc{B}$ of $b\!\in\!\mf{M}_{g,k}$, a smooth trivialization $\varphi$ as in (\ref{ST_e}) gives rise to a Banach vector bundle structure on the restriction 
\bEq{BVB_e3}
\cE^{\ell-1,p}_{A,\mfs}\big((\cC,\vec{\mf{z}}|_{\cB}),(X,D)\big)\lra W^{\ell,p}_{A,\mfs}\big((\cC,\vec{\mf{z}}|_{\cB}),(X,D)\big)
\eEq
such that $\dbar^\tn{log}\!-\!\nu_{\log}$ is a smooth section of that. 

\bRm{rem:actual-Obs}
Let $f\!=\!(u,C)\in \cM_{g,\mfs}(X,D,A)$ (no perturbation here). Similarly to the classical case (see \cite[Sec~24.1]{Mirror} and \cite[Rmk~6.2.1]{FF}), deformation theory of $f$, i.e. if we allow deformations of both $u$ and $C$, is  described by the long exact sequence 
\bEq{equ:long-def}
\aligned
 & &0\lra &~\tn{aut}(C) \stackrel{\de}{\lra} \\
 &\Def_{\log}(u)\lra &\Def_{\log}(f) \lra &~\Def(C) \stackrel{\de}{\lra} \\
	&\Obs_{\log}(u)  \lra  &~\Obs_{\log}(f) \lra &~0,
\endaligned
\eEq
where 
$$
\tn{aut}(C)=H^0_{\dbar}(T\Si(-\log\vec{z}))\qquad\tn{and}\qquad \Def(C)\cong T_b\cB\cong H^1_{\dbar}(T\Si(-\log \vec{z})).
$$  
If $\tn{Obs}_{\log}(f)=0$, then a small neighborhood $B(f)$ of $f$ in $\cM_{g,\mfs}(X,D,A)$ is a smooth orbifold of the expected real dimension (\ref{dlog_e}). 
The long exact sequence (\ref{equ:long-def}) is the hypercohomology of a short exact sequence of complexes of fine sheaves constructed in the following way. In order to simplify the notation, for a complex vector bundle $E\!\lra\! (\Si,\mfj)$ let $\Om^0(E)$ and  $\Om^{0,1}(E)$ denote the associated fine sheaves of smooth sections of $E$ and of smooth $E$-valued $(0,1)$-forms, respectively.  The map $\nd u\colon T\Si\!\lra\!TX$ gives rise to a \textit{logarithmic derivative} map 
$$
\nd^{\log} u\colon T\Si(-\log \vec{z})\lra u^*TX(-\log D)
$$ 
such that the following diagram commutes:
 \bEq{dulog_e1}
\xymatrix{
& T\Si(-\log \vec{z}) \ar[rr]^{\nd^{\log} u}\ar[d]^{\iota_{\Si,\vec{z}}} && u^*TX(-\log D)\ar[d]^{\iota_{X,D}}\\
& T\Si \ar[rr]^{\nd u} && u^*TX\;.
}
\eEq
Away from the contact points $\mf{d}$, by the first bullet in Definition~\ref{CTs_dfn} and the identification 
 $$
 TX(-\log D)|_{X-D}\stackrel{\iota_{X,D}}{\cong} TX|_{X-D},
 $$
 we have $\nd^{\log} u=\iota_{X,D}^{-1}\circ \nd u \circ \iota_{\Si,\vec{z}}$.  For each contact point $z_a\!\in\!\vec{z}$, restricted to the chart $U_a$ in Definition~\ref{CTs_dfn}, by (\ref{Locu_e}), $\nd^{\log} u$ is given by
 $$
 \Psi_{I_a}^{-1}\circ \nd^{\log} u= \pi_{I_a}^*\nd u_a\oplus \bigoplus_{i\in I_a}  \big(s_{ai} \frac{\nd w_a}{w_a}+ \frac{\nabla^{(I_a;i)} \eta_{a,i}}{\eta_{a,i}}\big),
 $$
 which maps the local generating section $w_a\partial w_a$ to 
 \bEq{Logdu_e}
 \pi_{I_a}^*\partial u_a(w_a\partial w_a)\oplus \bigoplus_{i\in I_a}  \big(s_{ai} + w_a\frac{\nabla^{(I_a;i)}_{\partial w_a} \eta_{a,i}}{\eta_{a,i}}\big).
\eEq
The following commutative diagram has exact rows:
$$
\xymatrix{
 0\ar[r]\ar[d]& \Om^0(T\Si(-\vec{z}))\ar[r]  \ar[d]^{\nd^{\log} u\oplus \dbar} & \Om^0(T\Si(-\vec{z})) \ar[d]^{\dbar}\\
 \Om^0(u^*TX(-\log D)) \ar[r]\ar[d]^{\tn{D}^{\log}_u\dbar}  &\Om^0(u^*TX(-\log D))\oplus \Omega^{0,1}(T\Si(-\vec{z}))\ar[r] \ar[d]^{\tn{D}^{\log}_u\dbar - \nd^{\log} u}& \Omega^{0,1}(T\Si(-\vec{z}))\ar[d]\\
 \Omega^{0,1}(u^*TX(-\log D))\ar[r] & \Omega^{0,1}(u^*TX(-\log D))\ar[r] & 0 \; ;
}
$$
i.e. it is an exact sequence of chain complexes given by the columns. Then, the deformation/obstruction long exact sequence (\ref{equ:long-def}) is the  hyper-cohomology of this diagram. By (\ref{Logdu_e}) and similarly to the classical case \cite[p. 284-285]{ST}, if $u$ is an immersion away from $\vec{z}$ and 
$$
u^{-1}(D)=\mf{d}=\{z_1,\ldots,z_k\}
$$ 
(i.e. $s_{a}\neq 0$ for all $a\!\in\![k]$), then $\nd^{\log} u$ is an embedding, the quotient 
$$
\cN_X\Si(-\log D)\equiv u^*TX(-\log D)/(\nd^{\log}_u~T\Si(-\vec{z}))
$$
is a complex vector bundle, and $\tn{D}^{\log}_u\dbar$ descends to a Fredholm operator $\tn{D}^{\log}_{\cN}\dbar$ on smooth sections of $\cN_X\Si(-\log D)$ such that 
\bEq{Obs_f}
\Def_{\log}(f)=\tn{ker}(\tn{D}^{\log}_\cN\dbar) \quad \tn{and}\quad \Obs_{\log}(f)=\tn{coker}(\tn{D}^{\log}_{\cN}\dbar).
\eEq
If $\nd^{\log} u$ is not an embedding, we still obtain a short exact sequence of sheaves of $\cO_\Si$-modules
$$
0\lra \cO(T\Si(-\log \vec{z}))\stackrel{\nd^{\log}u}{\lra} \cO(u^*TX(-\log D))\lra \cN \lra 0
$$
such that 
$$
\cN= \cO(\cN_X\Si(-\log D))\oplus \cN^{\tn{tor}}
$$ 
is the direct sum of sheaf of holomorphic sections of an $(n-1)$-dimensional holomorphic vector bundle $\cN_X\Si(-\log D)$ and a skyscraper sheaf $\cN^{\tn{tor}}$. Furthermore, $\tn{D}^{\log}_u\dbar$ descends to a Fredholm operator $\tn{D}^{\log}_{\cN}\dbar$ on smooth sections of $\cN_X\Si(-\log D)$ such that 
\bEq{Obs_f2}
\Def_{\log}(f)=\tn{ker}(\tn{D}^{\log}_\cN\dbar)\oplus H^0(\cN^{\tn{tor}}) \quad \tn{and}\quad \Obs_{\log}(f)=\tn{coker}(\tn{D}^{\log}_{\cN}\dbar).
\eEq
In particular, $ \Obs_{\log}(f)\!=\!0$ whenever $\tn{dim}_\C X\!=\!1$.
\eRm
\vskip.1in

\newtheorem*{PF-Trans1_prp}{Proof of Proposition~\ref{Trans1_prp}}
\begin{PF-Trans1_prp}
Since every map $u$ in $W^{\ell,p}_{A,\mfs}\big((\mf{U}_{g,k},\vec{\mf{z}}),(X,D)\big)$ meets $D$ only at finitely many points, by substituting $D_u\{\dbar -\nu\}$ with $D^{\log}_u\{\dbar -\nu\}$, Proposition~\ref{Trans1_prp} essentially follows from restricting the arguments of the proof of \cite[Thm~3.1]{RT} and \cite[Thm~3.1.5]{MS2} to maps in $W^{\ell,p}_{A,\mfs}\big((\mf{U}_{g,k},\vec{\mf{z}}),(X,D)\big)$.\\

\noindent 
More precisely, for $m\!>\!\ell$, let $\cH^m_{g,k}(X,D)_{\cR,J}$ denote the completion of the vector space $\cH_{g,k}(X,D)$ in $C^m$-topology. The universal moduli space 
\bEq{UM_e}
\aligned
&\mf{M}_{g,\mfs}(X,D,A)=\\
&\bigg\{ \big((c,u),\nu\big)\!\in\!W^{\ell,p}_{A,\mfs}\big((\mf{U}_{g,k},\vec{\mf{z}}),(X,D)\big)\!\times\!\cH^m_{g,k}(X,D)_{\cR,J}\colon \dbar^{\log} u (x)= \nu_{\log}(x,u(x))~~~\forall~x\!\in\!\pi^{-1}(c)\bigg\}
\endaligned
\eEq
is the zero set of the section\footnote{To be precise, the right-hand side should be $\pi_1^* E^{\ell,p}_{A,\mfs}\big((\mf{U}_{g,k},\vec{\mf{z}}),(X,D)\big)$, where $\pi_1$ is projection map to the first component
$$
W^{\ell,p}_{A,\mfs}\big((\mf{U}_{g,k},\vec{\mf{z}}),(X,D)\big)\!\times\!\cH^m_{g,k}(X,D)_{\cR,J}\lra W^{\ell,p}_{A,\mfs}\big((\mf{U}_{g,k},\vec{\mf{z}}),(X,D)\big).
$$ We avoid these details to keep the notation short.}
\bEq{UCRnusec_e2}
\dbar^{\log}-\nu_{\log}\colon W^{\ell,p}_{A,\mfs}\big((\mf{U}_{g,k},\vec{\mf{z}}),(X,D)\big)\!\times\!\cH^m_{g,k}(X,D)_{\cR,J}\lra E^{\ell,p}_{A,\mfs}\big((\mf{U}_{g,k},\vec{\mf{z}}),(X,D)\big)
\eEq
and is independent of $(\ell,p)$ by the elliptic regularity.  Restricting to each sufficiently small sub-universal family $\cC\!\lra\!\cB$ and fixing a smooth trivialization $\varphi$ as above, the restricted section 
\bEq{UCRnusec_e22}
\dbar^{\log}-\nu_{\log}\colon W^{\ell,p}_{A,\mfs}\big((\cC,\vec{\mf{z}}|_\cB),(X,D)\big)\!\times\!\cH^m_{g,k}(X,D)_{\cR,J}\lra  E^{\ell,p}_{A,\mfs}\big((\cC,\vec{\mf{z}}|_{\cB}),(X,D)\big)
\eEq
is $C^{m-\ell}$-smooth. With the same reasoning as in the argument leading to the surjectivity of \cite[(3.12)]{RT}, for $c\!\in\! \mc{B}$ and $((c,u),\nu\big)$ in the universal moduli space (\ref{UM_e}), the linearization map
\bEq{LCRNu_e} 
\aligned
D^{\log}_{((c,u),\nu)}\{\dbar-\nu\}\colon T_{(c,u)}W^{\ell,p}_{A,\mfs}\big((\cC,\vec{\mf{z}}|_\cB),(X,D)\big)&\oplus T_{\nu}\cH^m_{g,k}(X,D)_{\cR,J}\lra \\
&\Gamma(\Si,\Om^{0,1}_{\pi^{-1}(c)}\otimes_\C u^*TX(-\log D))
\endaligned
\eEq
of the section (\ref{UCRnusec_e2}) is surjective. This is due to the fact that $\tn{coker}(D_u^{\log}\{\dbar-\nu\})$ can be represented by sections supported away from the contact points where everything has a classical form. Therefore, the universal moduli space
$$
\mf{M}_{g,\mfs}(X,D,A)|_{\mc{B}}
$$ 
is a separable $C^{m-\ell}$-smooth Banach manifold. Here the restriction to $\cB$ means we are only considering $(J,\nu)$-maps with domain in $\cB$ corresponding to (\ref{UCRnusec_e2}). Then by the Sard-Smale Theorem, the set of regular values  $\cH^{\tn{reg}(\cB)}_{g,k}(X,D)_{\cR,J}$ of the projection map
$$
\pi_2 \colon \mf{M}_{g,\mfs}(X,D,A)|_{\mc{B}}\lra \cH^m_{g,k}(X,D)_{\cR,J}
$$
is Baire set of second category. For every $\nu\!\in\!\cH^{\tn{reg}(\cB)}_{g,k}(X,D)_{\cR,J}$, 
$$
\mc{M}_{g,\mfs}(X,D,A,\nu)|_{\mc{B}}= \pi_2^{-1}(\nu)
$$
is a smooth manifold of the expected dimension. Cover $ \mf{M}_{g,k}$ with countably many charts $\{\cB_i\}_{i=1}^\infty$ and let 
$$
\cH^{\tn{reg}}_{g,k}(X,D)_{\cR,J}= \bigcap_{i=1}^\infty \cH^{\tn{reg}(\cB_i)}_{g,k}(X,D)_{\cR,J}.
$$ 
This is still a Baire set of second category so that for each $\nu\!\in\!\cH^{\tn{reg}}_{g,k}(X,D)_{\cR,J}$, 
$$
\mc{M}_{g,\mfs}(X,D,A,\nu)
$$
is a smooth manifold of the expected dimension. With an argument similar to the Taubes' trick in the proof of \cite[Thm~3.1.6(ii)]{MS2}, we conclude that the subset of smooth perturbations 
$$
\cH^{\eset}_{g,k}(X,D)_{\cR,J}=\cH^{\tn{reg}}_{g,k}(X,D)_{\cR,J}\cap \cH_{g,k}(X,D)_{\cR,J}
$$
satisfies the first statement in Proposition~\ref{Trans1_prp}. \\

\noindent
With the modifications above, proof of the second statement is similar to the proof of \cite[Thm~3.1.5]{MS2}. More precisely, for $g\!=\!0$ and no perturbation, if $k\!\geq\!3$, let 
$$
\mf{U}_{0,k}=\cM_{0,k+1}\lra \mf{M}_{0,k}=\cM_{0,k}
$$ 
denote  the universal curve, and if $k\!<\!3$, let $\mf{U}_{0,k}=\P^1$, $\mf{M}_{0,k}$ be a point, and $\tn{aut}_k=\tn{Aut}(\P^1,\vec{z})$. 
Let
\bEq{UM_e0}
\aligned
\wt{\mf{M}}_{0,\mfs}(X,D,A)=\bigg\{ \big((c,u),J\big)\!\in\!W^{\ell,p}_{A,\mfs}\big((\mf{U}_{0,k},\vec{\mf{z}}),(X,D)\big)\!\times\!\tn{AK}^m(X,D)_{\cR}\colon \dbar^{\log} u (x)=0\bigg\}
\endaligned
\eEq
be the zero set of the section
\bEq{UCRnusec_e20}
\dbar^{\log}\colon W^{\ell,p}_{A,\mfs}\big((\mf{U}_{0,k},\vec{\mf{z}}),(X,D)\big)\!\times\!\tn{AK}^m(X,D)_{\cR}\lra  E^{\ell,p}_{A,\mfs}\big((\mf{U}_{0,k},\vec{\mf{z}}),(X,D)\big).
\eEq
The universal moduli space $\mf{M}_{0,\mfs}(X,D,A)$ is the quotient of $\wt{\mf{M}}_{0,\mfs}(X,D,A)$ with respect to $\tn{aut}_k$.
Since $\tn{coker}\big(D_u\dbar^{\log}\big)$ can be represented by sections supported away from the contact points, the same reasoning as in the proof of \cite[Thm~3.15]{MS2} shows that the linearization map
\bEq{LCRNu_e0} 
\aligned
D_{((c,u),J)}\dbar^{\log}\colon T_{(c,u)}W^{\ell,p}_{A,\mfs}\big((\cC,\vec{\mf{z}}|_\cB),(X,D)\big)&\oplus T_{J}\tn{AK}^m(X,D)_{\cR}\lra \\
&\Gamma(\Si,\Om^{0,1}_{\pi^{-1}(c)}\otimes_\C u^*TX(-\log D))
\endaligned
\eEq
is surjective, whenever $u$ is simple. Therefore, the subset of simple maps $\mf{M}^{\star}_{0,\mfs}(X,D,A)$  in the universal moduli space is a separable $C^{m-\ell}$-smooth Banach manifold. Then by the Sard-Smale Theorem, the set of regular values  $\tn{AK}^{\tn{reg}}(X,D)_{\cR}$ of the projection map
$$
\pi_2 \colon \mf{M}^{\star}_{0,\mfs}(X,D,A)\lra \tn{AK}^{\tn{reg}}(X,D)_{\cR}
$$
is Baire set of second category. With an argument similar to the Taubes' trick in the proof of \cite[Thm~3.1.6(ii)]{MS2}, we conclude that the subset of smooth perturbations 
$$
\tn{AK}^{\eset}(X,D)_{\cR}=\tn{AK}^{\tn{reg}}(X,D)_{\cR}\cap \tn{AK}(X,D)_{\cR}
$$
satisfies the second statement in Proposition~\ref{Trans1_prp}. 

\qed
\end{PF-Trans1_prp}

\noindent
Moving to the simple nodal case in Section~\ref{NodalMap_ss}, we will need to show that certain evaluation maps on the universal moduli spaces are transverse.
In the analytical set up of \cite{MS2}, the proof of transversality of evaluation maps in \cite[Prp~6.2.8]{MS2} uses \cite[Prp~3.4.2]{MS2} and induction on the number of edges. Proposition 3.4.2 in \cite{MS2}, itself, is a consequence of \cite[Lmm~3.4.3]{MS2}. We will need the following natural generalization of \cite[Lmm~3.4.3]{MS2} to show that the evaluation maps at the nodes and $\tn{ob}_{\Gamma}$ are transverse. With $\dbar^{\log}$ and $D_u^{\log}\dbar$ in place of $\dbar$ and $D_u\dbar$, respectively, its proof is similar to the (long and explicit) proof of  \cite[Lmm~3.4.3]{MS2}.

\bLm{Ydef_lmm}
Let $A\!\neq\! 0$ and $(\Si,\vec{z})$ be a smooth $k$-marked curve\footnote{We just need the sphere case.}. With notation as above, let 
$$
(u,J) \in \wt{\mf{M}}^{\star}_{g,\mfs}(X,D,A) \subset W^{\ell,p}_{A,\mfs}\big((\Si,\vec{z}),(X, D)\big)\!\times\!\tn{AK}^m(X,D)_{\cR}.
$$
For each $a\in [k]$, let $\xi_a$ be a log tangent vector in $T_{u(z_a)}X(-\log D)_\cR$. For each open set $U\!\subset\! \Si-\vec{z}$, there exists  
$$
\xi\!\in\!W^{\ell,p}\big(\Si, u^*TX(- \log D)_{\cR}\big)\quad \tn{and}\quad Y\!\in\! T_J\tn{AK}^m(X,D)_{\cR}
$$ 
such that 
$$
\xi(z_a)=\xi_a\quad \forall~a\!\in\! [k],\qquad \tn{Supp}(Y|_{u(\Si)})\!\subset\! U, \quad\tn{and}\quad D_u\dbar^{\log} \xi+ \frac{1}{2}\,Y \circ \nd u \circ \mfj =0.
$$
\eLm

%--------------------------------------------------------
\subsection{Depth-$I$ maps}\label{NGSmooth_ss}

\noindent 
For each $\nu\!\in\!\cH_{g,k}(X,D)_{\cR,J}$, consider a log map 
$$
f\!=\!\big[u, (\ze_i)_{i\in I}, (\Si,\mfj, z_1,\ldots,z_k)\big]\!\in\!\ov\cM^{\log}_{g,\mfs}(X,D,A,\nu)
$$ 
where $\Si$ is smooth, i.e. $u(\Si)\!\subset\! D_I$ for a non-trivial maximal subset $I\!\subset\!S$, $\ord_{z_a}(u,D_i)\!=\!s_{ai}\!\geq \!0$ for all $i\!\in\![N]-I$, and $\ord_{z_a}(\ze_i)\!=\!s_{ai}$ for all $i\!\in\! I$. We allow $s_{ai}$ to be negative for $i\!\in\! I$. Let 
$$
\cM_{g,\mfs}(X,D,A,\nu)_I\subset \ov{\cM}_{g,\mfs}(X,D,A,\nu)
$$
be the stratum of such maps. The stratum $\cM_{g,\mfs}(X,D,A,\nu)_I$ is a generalization of the main stratum 
$$
\cM_{g,\mfs}(X,D,A,\nu)_\eset=\cM_{g,\mfs}(X,D,A,\nu)
$$
where the domain is still smooth, but the image could lie in a non-trivial stratum of the divisor.\\

\noindent
Forgetting the meromorphic sections $\ze_i$, for the same reason\footnote{Fixing a set of marked points, up to multiplication by a constant, there is at most one meromorphic section of any holomorphic line bundle with prescribed zeros/poles at the marked points and nowhere else.}   as in \cite[Rmk~3.1]{FT1}, we get a topological embedding
\bEq{I-embedding}
\cM_{g,\mfs}(X,D,A,\nu)_I\hookrightarrow \cM_{g,\ov\mfs}(D_I,\partial D_I,A,\nu_I), \qquad [\phi,u, (\ze_i)_{i\in I}, \Si,\vec{z}\,]\lra [\phi,u, \Si,\vec{z}\,],
\eEq
where $\partial D_I\!\subset\! D_I$ is the boundary divisor as in (\ref{parDI_e}) and
$$
\ov\mfs\!=\!\big(s_a\!=\!(s_{ai})_{i\in [N]-I}\big)_{a\in [k]}\!\in\!(\N^{[N]-I})^k.
$$
In this section we prove the following transversality argument.
\bPr{Trans2_prp}
Suppose $(X,\om)$ is a closed symplectic manifold, $D=\bigcup_{i\in [N]} D_i$ is an SNC symplectic divisor,  $A\!\in\!H_2(X,\Z)$, $g,k\!\in\!\N$,  and $\mfs\!\in\!(\Z^N)^k$, with $\ov\mfs\in(\N^{[N]-I})^k$. 
\bEn
\item If $2g\!+\!k\!\geq\!3$, for any given choice of universal family in (\ref{UFamily_e}), there exists a Baire set of second category $\cH^{I}_{g,k}(X,D)_{\cR,J}\!\subset\!\cH_{g,k}(X,D)_{\cR,J}$ such that for each $\nu\!\in\!\cH^{I}_{g,k}(X,D)_{\cR,J}$, $\cM_{g,\mfs}(X,D,A,\nu)_I$ is a naturally oriented smooth manifold of the real dimension 
$$
2\big(c_1^{TX(-\log D)}(A)+(n-3)(1-g)+k-|I|\big).
$$
\item If $g\!=\!0$, $\nu\!\equiv\!0$, and $A\!\neq\!0$ or $k\!\geq\!3$, the same statement holds for $J$ in a Baire set of second category $\tn{AK}^{I}(X,D)_{\cR}\!\subset\!\tn{AK}(X,D)_{\cR}$, if we restrict to the subspace of simple maps $\cM^\star_{0,\mfs}(X,D,A)_I$.
\eEn
\ePr

\vskip.1in
\noindent
Proposition~\ref{Trans1_prp} is a special case of Proposition~\ref{Trans2_prp}, where $I\!=\!\eset$.
By (\ref{NufromNuLog_e}) and (\ref{stdnu_e}), in a neighborhood of $D_I$, every  $\nu\!\in\!\cH_{g,k}(X,D)_{\cR,J}$ can be decomposed as $\pi_I^*\nu_I\oplus n_I$, where 
$$
\nu_I\!\in\!\cH_{g,k}(D_I,\partial D_I)_{\cR_I,J_I},
$$
and $n_I$ is determined by a family of $\C^I$-valued $(0,1)$-forms
$$
\theta_I\in \Theta_{g,k}(D_I)\equiv \bigg\{(\theta_{I,i})_{i\in I}\! \in\! \Gamma\big(\ov{\mf{U}}^\star_{g,k}\times D_I, \pi_1^*\Om^{0,1}_{g,k}\otimes_\C \C^I\big)\colon \tn{supp}(\theta_{I,i})\!\subset\big(\ov{\mf{U}}_{g,k}^\star-\bigcup_{a\in [k]}\tn{Im}(\mf{z}_a)\big)\!\times D_I \bigg\}.
$$
The map 
\bEq{NTpair_e}
\tn{Res}\colon \cH_{g,k}(X,D)_{\cR,J}\lra \cH_{g,k}(D_I,\partial D_I)_{\cR_I,J_I}\times \Theta_{g,k}(D_I), \qquad \nu \lra (\nu_I,\theta_I),
\eEq
is surjective and continuous. Therefore, in order to prove the first statement of Proposition~\ref{Trans2_prp}, it is enough to find  a Baire set of second category 
\bEq{separation_e}
\cH\Theta^{\tn{reg}}_{g,k}(D_I,\partial D_I)_{\cR_I,J_I}\!\subset\!\cH_{g,k}(D_I,\partial D_I)_{\cR_I,J_I}\times \Theta_{g,k}(D_I)
\eEq
such that for each $\nu$ with $\tn{Res}(\nu)\!:=\!(\nu_I,\theta_I)$ in this set, the statement of Proposition~\ref{Trans2_prp} holds. 
Similarly, in case (2), the map 
$$
\tn{AK}(X,D)_{\cR}\lra \tn{AK}(D_I,\partial D_I)_{\cR_I}
$$
is surjective and continuous. Below, we show that we can take $\tn{AK}^{I}(X,D)_{\cR}$ to be the preimage of $\tn{AK}^{\eset}(D_I,\partial D_I)_{\cR_I}$ given by the second part of Proposition~\ref{Trans1_prp}. Therefore, the main goal of this section is to describe the normal bundle of the embedding (\ref{I-embedding}) and $\cH\Theta^{\tn{reg}}_{g,k}(D_I,\partial D_I)_{\cR_I,J_I}$. Proposition~\ref{Trans2_prp} can also be obtained from (proof of) Proposition~\ref{Trans1_prp} by a looking $(u, (\ze_i)_{i\in I})$ as a log map into the fiber product (\ref{PP_e}) of the projectivizations of  $\cN_XD_i$. We will explain this argument in Remark~\ref{AltPf_rmk}. \\

\noindent
For each complex curve $(\Si,\mfj)$, let $\tn{Pic}^0(\Si,\mfj)$ be the group of degree $0$ holomorphic line bundles on $(\Si,\mfj)$, $\mc{O}=\mc{O}_{\Si,\mfj}\!\in\!\tn{Pic}^0(\Si,\mfj)$ be the trivial line bundle, and $\cO^I\equiv \bigoplus_{i\in I}\cO$.
Let 
$$
\tn{Pic}^0(\mf{U}_{g,k})\lra \mf{M}_{g,k}
$$ 
be the fiber bundle whose fiber over every $c\in \mf{M}_{g,k}$ is $\tn{Pic}^0(\pi^{-1}(c))$. In the following, by $\cO$ we mean the section 
$$
\cO\colon \mf{M}_{g,k}\lra \tn{Pic}^0(\mf{U}_{g,k})
$$ 
that takes $c$ to the trivial line bundle $\cO_{\pi^{-1}(c)}$. Image of $\cO$ has complex codimension $g$. By abuse of notation, we also let $\tn{Pic}^0(\mf{U}_{g,k})$ to denote the pull back of $\tn{Pic}^0(\mf{U}_{g,k})$ to $\mc{M}_{g,\mfs}(D_I,\partial D_I,A,\nu_I)$ (or any other configuration space).\\

\noindent
The next Lemma describes the (virtual) normal bundle of the embedding (\ref{I-embedding}).
\bLm{HodgeObs_lmm}
For each $(\nu_I,\theta_I)=\tn{Res}(\nu)$ as in (\ref{NTpair_e}), there exists a natural map 
\bEq{ptheta_e}
P_{\theta_I}=(P_{\theta_{I,i}})_{i\in I}\colon\cM_{g,\mfs_I}(D_I,\partial D_I,A,\nu_I)\!\lra\! \tn{Pic}^0(\mf{U}_{g,k})^I
\eEq
such that 
$$
\cM_{g,\mfs}(X,D,A,\nu)_I=P_{\theta_I}^{-1}\big(\mc{O}^I\big).
$$ 
In particular, 
\bEq{g0MI_e}
\cM_{0,\mfs}(X,D,A,\nu)_I=\cM_{0,\mfs_I}(D_I,\partial D_I,A,\nu_I).
\eEq
\eLm 

\bPf
For each $i\!\in\!I$ and $\big(\phi,u, (\Si,\mfj,\vec{z})\big)\!\in\!\cM_{g,\mfs_I}(D_I,\partial D_I,A,\nu_I)$, define 
$$
P_{\theta_{I,i}}\big(\phi,u, (\Si,\mfj,\vec{z})\big)=u^*\cN_XD_i \otimes \cO_\Si\big(-\sum_{a\in [k]} s_{ai} z_a\big) \in \tn{Pic}^0(\Si,\mfj),
$$ 
where the holomorphic structure on the line bundle $u^*\cN_XD_i$ is given by the $\dbar$-operator
$$
\tn{D}_{u}^{\cN_i}\{\dbar\!-\!\nu\}= \dbar_{u^*\cN_XD_i}-(\phi,u)^*\theta_{I,i}
$$ 
and  the second term is the line bundle corresponding to the divisor $\sum_{a=1}^k s_{ai} z_a$. By definition, $P_{\theta_{I,i}}(\phi,u, \Si,\vec{z})\!=\!\cO$ if and only if there exists a non-trivial meromorphic section $\ze_i$ of $u^*\cN_XD_i$ with zeros/poles of order $s_{ai}$ and $z_a$, for all $a\!=\!1,\ldots,k$ (and nowhere else).
\ePf

\bRm{Noze_e}
In light of Lemma~\ref{HodgeObs_lmm}, the moduli space $\ov\cM^{\log}_{g,\mfs}(X,D,A,\nu)$ can be described without mentioning the meromorphic sections $\ze_{v,i}$ in the following way. This explains the absence of these sections in the proof of Proposition~\ref{Trans2_prp} part (1) and other proofs.
An element of $\ov\cM^{\log}_{g,\mfs}(X,D,A,\nu)$ is the equivalence class of a stable $(J,\nu)$-holomorphic map 
$$
\big(u_v,\Si_v,\mfj_v,\vec{z}_v\cup q_{v}\big)_{v\in \V},
$$
together with a choice of decorations $\{s_{\uvec{e}}\}_{\uvec{e}\in \uvec{\E}}$ on the nodal points such that 
such that 
\bEn
\item 
$$
s_{\uvec{e}}\!=\!-s_{\scz\ucev{e}} \quad \forall~e\!\in\!\E,\qquad  
\sum_{\uvec{e}\in \uvec{\E}_v} s_{\uvec{e}}+\sum_{l\in \L_v} s_l=(A_v\cdot D_i)_{i\in [N]}\quad \forall~v\!\in\!\V;
$$
\item for each $v\!\in\! \V$, $\uvec{e}\!\in\!\uvec{\E}_v$, and $i\!\not\in \!I_v$, $u_v$ has a tangency of order $s_{\uvec{e},i}$ with $D_i$ at $q_{\uvec{e}}$;
\item for each $v\!\in\! \V$, $l\!\in\!\L_v$, and $i\!\not\in \!I_v$, $u_v$ has a tangency of order $s_{a_l,i}$ with $D_i$ at $z_{a_l}$; 
\item there exists a vector-valued function $s\colon\!\V\!\lra\!\R^N$ such that $s_v\!=\!s(v)\!\in\! \R_{+}^{I_v}\!\times \{0\}^{[N]-I_v}$ for all $v\!\in\!\V$, and 
$$
s_{v_{2}(\uvec{e})}\!-\!s_{v_1(\uvec{e})}\!=\!\la_{e} s_{\uvec{e}}\quad\tn{for some}~\la_e\!>\!0,\qquad \forall~\uvec{e}\!\in\!\uvec{\E};
$$
\item\label{OC_l}$u_v^*\cN_XD_i \cong \cL_{v,i}\equiv \cO_{\Si_v}\big(\sum_{l\in \L_v} s_{a_l,i} \,z_{a_l}+\sum_{\uvec{e}\in \uvec{\E}_v} s_{\uvec{e},i}\, q_{\uvec{e}}\big)$, for all $v\!\in\!\V$ and $i\!\in\!I_v$;
\item and, $\tn{ob}_\Gamma(f)\!=\!1$.
\eEn
The last condition can (in theory) be expressed in terms of $u_v$, the canonical sections of $\cL_{v,i}$, and the isomorphisms of the holomorphic line bundles in \ref{OC_l}.
\eRm

\noindent
The following statements are immediate corollaries of the first and second statements of Lemma~\ref{HodgeObs_lmm}, respectively.

\noindent
\bCr{Ist-dim_cr}
Replacing $(X,D)$ with $(D_I,\partial D_I)$ in Proposition~\ref{Trans1_prp}, let $\cH^{\eset}_{g,k}(D_I,\partial D_I)_{\cR_I,J_I}$ and $\tn{AK}^{\eset}(D_I,\partial D_I)_{\cR_I}$ be the resulting sets of regular perturbations and almost complex structures on $D_I$, respectively.\vskip.1in

\noindent
$\tn(1)$ If $\nu_I\!\in\!\cH^{\eset}_{g,k}(D_I,\partial D_I)_{\cR_I,J_I}$ and the image of $P_{\theta_I}$ is transverse to $\cO^I$  then $\cM_{g,\mfs}(X,D,A,\nu)_I$ is a naturally oriented manifold of the expected real dimension 
\bEq{IStratum_e}
\aligned
&2\big(c_1^{TD_I(-\log \ov{D})}(A)+(n-|I|-3)(1-g)+k-|I| \dim_\C \tn{Pic}^0(\Si)\big)=\\
&2\big(c_1^{TX(-\log D)}(A)+(n-3)(1-g)+k-|I|\big).
\endaligned
\eEq
$\tn{(2)}$ Similarly, if $J_I\!=\!J|_{TD_I}\!\in\!\tn{AK}^{\eset}(D_I,\partial D_I)_{\cR_I}$ then $\cM^\star_{0,\mfs}(X,D,A)_I$ is a naturally oriented manifold of the expected real dimension 
$$
2\big(c_1^{TD_I(-\log \ov{D})}(A)+(n-|I|-3)+k\big)=2\big(c_1^{TX(-\log D)}(A)+(n-3)+k-|I|\big).
$$
\eCr
\noindent
\textit{The second statement establishes part (2) of Proposition~\ref{Trans2_prp}}. It also shows that in the genus $0$ case, Condition~\ref{OC_l} in Remark~\ref{Noze_e} is automatically satisfied. By the first statement, in order to prove Proposition~\ref{Trans2_prp}.(1), we need to show that for generic $\theta_I$ and $\nu_I$,  the image of $P_{\theta_I}$ is transverse to $\cO^I$.
\vskip.1in

\newtheorem*{PF-Trans2_prp}{Proof of Proposition~\ref{Trans2_prp} part (1)}
\begin{PF-Trans2_prp}
Let $\Theta^m_{g,k}(D_I)$ be the completion of $\Theta^m_{g,k}(D_I)$ in the $C^m$-norm.
Fix a sub-universal family $\cC\!\lra\!\cB$ of $\mf{U}_{g,k}$ around $C\!=\!(\Si,\mfj, z_1,\ldots,z_k)$  and a smooth trivialization $\varphi$ of $\cC$ as in (\ref{LocalUnivFam_e}) and (\ref{ST_e}), respectively.
Consider the configuration space
\bEq{WIConf_e}
W_I=W^{\ell,p}_{\ov\mfs,A}\big((\cC,\vec{\mf{z}}|_{\cB}),(D_I,\partial D_I)\big)_{\varphi}\!\times\!\cH^m_{g,k}(D_I,\partial D_I)_{\cR_I,J_I}\times \Theta^m_{g,k}(D_I).
\eEq
The map $P$ in  Lemma~\ref{HodgeObs_lmm}, extends to a map 
\bEq{GenPMap_e}
P\colon W_I\lra \tn{Pic}^0(\mf{U}_{g,k})^I, \quad ((c,u),\nu_I,\theta_I)\quad \tn{to} \quad P_{\theta_I}\big(\tn{id}_{\pi^{-1}(c)},u,(\pi^{-1}(c),\vec\mfz(c))\big).
\eEq
The universal moduli space
\bEq{UM_e2}
\aligned
\mf{M}_{g,\mfs}(X, D,A)_I=\bigg\{ &\big((c,u),\nu_I,\theta _I\big)\!\in\!W_I \colon \\
&\dbar^{\log} u (x)= \nu_{\log,I}(x,u(x))~~~\forall~x\!\in\!\pi^{-1}(c),~P\big((c,u),\nu_I,\theta _I\big)=\cO^I\bigg\}
\endaligned
\eEq
is the $(0\!\oplus\! \cO^I)$-level set of 
\bEq{CRP_e}
(\dbar^{\log}-\nu_{\log,I})\!\times\! P \colon W_I  \lra  E^{\ell,p}_{\ov\mfs,A}\big((\cC,\vec{\mf{z}}|_{\cB}),(D_I,\partial D_I)\big) \!\times\! (\tn{Pic}^0(\mf{U}_{g,k}))^I.
\eEq
We show that $0\times \cO^I$ is a regular value.\\

\noindent
Since $\Theta^m_{g,k}(D_I)$ is a linear space we have 
$$
T_{\theta_I}\Theta^m_{g,k}(D_I)\cong \Theta^m_{g,k}(D_I),\qquad \forall~\theta_I\!\in\!\Theta^m_{g,k}(D_I).
$$
For every $((c,u),\nu_I,\theta_I)\in \mf{M}_{g,\mfs}(X, D,A)_I$,
the normal component of the linearization of (\ref{CRP_e}) at $((c,u),\nu_I,\theta_I)$ has the form
\bEq{DCRNuP_e}
\aligned
&D_{((c,u),\nu_I,\theta_I)}(\dbar^{\log}\!-\!\nu_{\log,I})\!\times\!P\colon  T_{(c,u)}W^{\ell,p}_{\ov\mfs,A}\big((\cC,\vec{\mf{z}}|_\cB),(D_I,\partial D_I)\big)\oplus T_{\nu_I}\cH^m_{g,k}(D_I,\partial D_I)_{\cR,J}\oplus T_{\theta_I}\Theta^m_{g,k}(D_I)\\
 & \lra  W^{\ell-1,p}\big(\Si,\Om^{0,1}_{\pi^{-1}(c)}\otimes_\C u^*TD_I(-\log \partial D_I)\big) \oplus  H^{0,1}(\Si,\mfj)^I\;.
 \endaligned
\eEq
It is given by (\ref{LCRNu_e}) on the first two components and sends $0\oplus 0\oplus \wt\theta_I$ to the cohomology class 
$$
0\oplus [(\tn{id}_{\pi^{-1}(c)},u)^* \wt{\theta}_I]\!\in\!H^{0,1}(\Si,\mfj)^I.
$$ 
For the same reason as in the proof of Proposition~\ref{Trans1_prp}, every element of 
$$
W^{\ell-1,p}(\Si,\Om^{0,1}_{\pi^{-1}(c)}\otimes_\C u^*TD_I(-\log \partial D_I))
$$ 
is in the image of the restriction of (\ref{DCRNuP_e}) to the first two summands of the domain, followed by the projection to the first component of the target. Since a representative of every cohomology class in $H^{0,1}(\Si,\mfj)^I$ can be extended to a global $(0,1)$-form on $\cU_{g,k}$, the restriction of (\ref{DCRNuP_e}) to the third summand of the domain is a map onto the second summand of the target. Therefore, (\ref{DCRNuP_e}) is surjective. Consequently, by Implicit-Function Theorem, the universal moduli space $\mf{M}_{g,\mfs}(X,D,A)_I$ is a separable $C^{m-\ell}$-smooth Banach manifold. Then by the Sard-Smale Theorem, the set of regular values  $\cH\Theta^{\tn{reg,m}}_{g,k}(D_I,\partial D_I)_{\cR_I,J_I}$ of the projection map
$$
\tn{proj}\colon \mf{M}_{g,\mfs}(X,D,A)_I\lra \cH^{m}_{g,k}(D_I,\partial D_I)_{\cR_I,J_I}\times\Theta^{m}_{g,k}(D_I)
$$
is Baire set of second category. By construction, for every $\nu$ with
$$
\tn{Res}(\nu)=(\nu_I,\theta_I) \in \cH\Theta^{\tn{reg,m}}_{g,k}(D_I,\partial D_I)_{\cR_I,J_I}
$$ 
the stratum
$$
\mf{M}_{g,\mfs}(X,D,A,\nu)_I=\tn{proj}^{-1}(\nu_I,\theta_I)
$$
is a naturally oriented smooth manifold of the expected dimension (\ref{IStratum_e}). With an argument similar to the Taubes' trick in the proof of \cite[Thm~3.1.6(ii)]{MS2}, we conclude that the subset of smooth perturbations 
$$
\cH\Theta^{\tn{reg}}_{g,k}(D_I,\partial D_I)_{\cR_I,J_I}=\cH\Theta^{\tn{reg,m}}_{g,k}(D_I,\partial D_I)_{\cR_I,J_I}\cap \Big(\cH_{g,k}(D_I,\partial D_I)_{\cR_I,J_I}\times \Theta_{g,k}(D_I)\Big)
$$
satisfies the requirement of (\ref{separation_e}). The set $\cH^{I}_{g,k}(X,D)_{\cR,J}=\tn{Res}^{-1}(\cH\Theta^{\tn{reg}}_{g,k}(D_I,\partial D_I)_{\cR_I,J_I})$ satisfies Proposition~\ref{Trans2_prp}.(1).
\qed

\end{PF-Trans2_prp}

\begin{remark}
Suppose $f=[\phi,u, \ze=(\ze_i)_{i\in I}, \Si,\vec{z}]\!\in\!\cM_{g,\mfs}(X,D,A,\nu)_I$. Restricted to $D_I$, by (\ref{TXlog_e}), we have 
\bEq{2parts_e}
TX(-\log D)|_{D_I}= TD_I(-\log \partial D_I)\oplus  D_I\times\C^I.
\eEq
Replacing $(X,D,\nu)$ in (\ref{Ddbarnu_e}) with $(D_I,\partial D_I,\nu_I)$, we get the linearized CR operator
\bEq{Ddbarnu_e2}
\tn{D}^{\log}_u\{\dbar-\nu_I\}\colon W^{\ell,p}(\Si,u^*TD_I(-\log \partial D_I))\lra W^{\ell-1,p}(\Si,\Om^{0,1}_{\Si}\otimes_\C u^*TD_I(-\log \partial D_I)).
\eEq
Let
\bEq{Ddbarnu_e3}
\dbar_{\tn{std}}\colon W^{\ell,p}(\Si,\C^I)\lra W^{\ell-1,p}(\Si,\Om^{0,1}_{\Si,\mfj}\otimes_\C \C^I)
\eEq
denote the standard $\dbar$-operator on the trivial bundle $\Si\times\C^I$.
In order to extend (\ref{Ddbarnu_e}) to the case of maps with smooth domain but image in a stratum $D_I$, using the decomposition (\ref{2parts_e}), define 
\bEq{Ddbarnu_e4}
\aligned
&\tn{D}^{\log}_{u}\{\dbar-\nu\}\colon W^{\ell,p}(\Si,u^*TX(-\log D))\lra W^{\ell-1,p}(\Si,\Om^{0,1}_{\Si}\otimes_\C u^*TX(-\log D)),\\
&\tn{D}^{\log}_{u}\{\dbar-\nu\}(\xi\oplus \eta\!=\!(\eta_i)_{i\in I})=\tn{D}^{\log}_u\{\dbar-\nu_I\}(\xi)\oplus \bigg(\dbar_{\tn{std}}(\eta_i)- \tn{D}_\xi \tn{D}_u^{\cN_i}\{\dbar-\nu\} \bigg)_{i\in I},
\endaligned
\eEq
where
$$
\tn{D}_\xi  \tn{D}_u^{\cN_i}\{\dbar-\nu\} \in \Gamma\big(\Si, \Om^{0,1}_{\Si,\mfj}\big)
$$
is the derivative of the $\dbar$-operator $\tn{D}_u^{\cN_i}\{\dbar-\nu\} $ in (\ref{CRLNuNormal_ei}) in the direction of $\xi$. Note that the derivative of a $1$-parameter family of $\dbar$-operators on a complex line bundle is a $(0,1)$-form. In other words, with respect to the decomposition (\ref{2parts_e}), $\tn{D}^{\log}_{(u,\ze)}\{\dbar-\nu\}$ is the lower diagonal operator 
$$
\begin{bmatrix} 
\tn{D}^{\log}_u\{\dbar-\nu_I\} & 0 \\
* & \dbar_{\tn{std}}
\end{bmatrix}
$$
with diagonal entires (\ref{Ddbarnu_e2}) and (\ref{Ddbarnu_e3}).   Cohomology class of 
$$
\big(\tn{D}_\xi  D_u^{\cN_i}\{\dbar-\nu\} \big)_{i\in I}\in \Gamma\big(\Si, \Om^{0,1}_{\Si,\mfj}\big)^I
$$ 
is 
$$
\tn{D}_uP(\xi)\in H^{0,1}(\Si,\mfj)^I.
$$
In light of the proof of Proposition~\ref{Trans2_prp}, Corollary~\ref{Ist-dim_cr} can be rephrased in the following way.

\bCr{Ist-dim_cr2} 
If $\tn{D}^{\log}_{u}\{\dbar-\nu\}$ is surjective then $\cM_{g,\mfs}(X,D,A,\nu)_I$ is cut transversely in a neighborhood of $f$. Furthermore, if $\mfs\!\in\!(\N^{[N]})^k$, then the whole moduli space $\ov\cM_{g,\mfs}(X,D,A,\nu)$ is cut transversely in a neighborhood of $f$.
\eCr

\noindent
The elements of the form $0\oplus (c_i)_{i\in I}$ in the kernel of $\tn{D}^{\log}_{u}\{\dbar-\nu\}$, where $c_i$ is a constant section of the trivial line bundle $\Si\!\times\!\C$, correspond to those deformations of $u$ that push the image of $u$ out of $D_i$ in the direction of $\ze_i$ by $c_i\ze_{i}$.
\end{remark}

\noindent
Proposition~\ref{Trans2_prp} can essentially be obtained from Proposition~\ref{Trans1_prp} by looking at $(u, (\ze_i)_{i\in I})$ as a log map into the fiber product (\ref{PP_e}) of the projectivizations of  $\cN_XD_i$, in the following way.

\bRm{AltPf_rmk}
Given $(X,D,\om,\cR)$, let 
\bEq{PP_e2}
\wt{X}\equiv \prod_{i\in I} \P(\cN_XD_i|_{D_I}\oplus \C) \lra D_I
\eEq
be the associated $(\P^1)^I$-fiber bundle over $D_I$.
For $i\!\in\! [N]\!-\!I$, let $\wt{D}_i=\pi^{-1} D_{i\cup I}$. For $i\!\in\! I$, let
$$
\wt{D}_i= (D_{i,0}\cup D_{i,\infty}) \times \prod_{j\in I-i} \P(\cN_XD_j|_{D_I}\oplus \C),
$$
where $D_{i,0}$ and $D_{i,\infty}$ are the zero and infinity divisors of $\P(\cN_XD_i|_{D_I}\oplus \C)$, respectively. The inclusion 
$$
\wt{D}=\bigcup_{i\in [N]} \wt{D}_i\subset \wt{X}
$$
is an SNC symplectic divisor with respect to the symplectic structure obtained from the standard symplectic structure on $\cN_XD_I$. As explained in \cite[p.~11]{FZ2}, from any $\cR$-compatible almost complex structure $J$ on $(X,D)$ and any $\nu\!\in\!\cH_{g,k}(X,D)$ we obtain a compatible $(\C^*)^I$-equivariant almost complex structure $\wt{J}$ on $(\wt{X},\wt{D})$ and a  $(\C^*)^I$-equivariant perturbation term $\wt\nu\!\in\!\cH_{g,k}(\wt{X},\wt{D})$. The latter only depends on $\nu_I$ and $\theta_I$. Every tuple $\big(\phi,u, \ze=(\ze_i)_{i\in I}, \Si,\vec{z}\big)$ representing an element of $\cM_{g,\mfs}(X,D,A,\nu)_I$ can be seen as a log $(\wt{J},\wt{\nu})$-map $\wt{u}\colon \Si\lra \wt{X}$ representing an element of 
$$
\cM_{g,\mfs}(\wt{X},\wt{D},\wt{A},\wt{\nu})
$$
where, for each $i\!\in\! I$, $s_{ai}>0$ denotes a tangency of order $s_{ai}$ with $D_{i,0}$ and $s_{ai}<0$ denotes a tangency of order $|s_{ai}|$ with $D_{i,\infty}$. We have 
\bEq{QSpace_e}
\cM_{g,\mfs}(X,D,A,\nu)_I\cong \cM_{g,\mfs}(\wt{X},\wt{D},\wt{A},\wt{\nu})/(\C^*)^I.
\eEq
Following the proof of Proposition~\ref{Trans1_prp}, we can show that for generic $(\wt{J},\wt{\nu})$, $\cM_{g,\mfs}(\wt{X},\wt{D},\wt{A},\wt{\nu})$ is a smooth oriented $(\C^*)^I$-equivariant manifold of real dimension 
$$
2\big( c_1^{T\wt{X}(-\log \wt{D})}(\wt{A}) + (n-3)(1-g)+k\big)= 2\big( c_1^{TX(-\log D)}(A) + (n-3)(1-g)+k\big).
$$
Then Proposition~\ref{Trans1_prp} follows from (\ref{QSpace_e}), and $\tn{D}^{\log}_{u}\{\dbar-\nu\}$ in (\ref{Ddbarnu_e4}) is equivalent to $\tn{D}^{\log}_{\wt{u}}\{\dbar-\wt{\nu}\}$ defined in  (\ref{Ddbarnu_e}).

\eRm

%--------------------------------------------------------
\subsection{Simple nodal maps}\label{NodalMap_ss} 
\noindent
Moving to the nodal case, let 
$$
\cM_{g,\mfs}(X,D,A,\nu)_\Gamma\subset \ov\cM_{g,\mfs}(X,D,A,\nu)
$$
be the stratum of stable nodal log $(J,\nu)$-holomorphic curves with the decorated dual graph $\Gamma$ (and $|\V|\!\geq\! 2$). With a set up similar to \cite[Sec~6.3]{FF}, the deformation/obstruction theory of $\cM_{g,\mfs}(X,D,A)_\Gamma$ around any $f\!=\!(u, [\ze], \Si,z_1,\ldots,z_k)$ is given by (1) the sum of $\tn{D}^{\log}_{u_v}\{\dbar-\nu\}$ over the irreducible components $\Si_v$, and (2) the obstruction map in (\ref{group-element_e}), i.e.,
\bEq{FP_e}
\tn{ob}_\Gamma\colon \times_{v\in \V} \cM_{g_v,\mfs_v}(X,D,A_v,\nu)_{I_v}\lra \mc{G}(\Gamma),
\eEq
where $\times_{v\in \V}$ denotes the fiber-product over the evaluation maps at the nodes.\\

\noindent
We  write $\V=\V_{\tn{p}}\cup \V_{\tn{b}}$, where $\V_{\tn{p}}$ corresponds to set of non-contracted or \textbf{principal} components and $ \V_{\tn{b}}$ corresponds to set of contracted or \textbf{bubble} components.
For each $v\!\in\!\V$, let $\cM^\star_{g_v,\mfs_v}(X,D,A_v,\nu)_{I_v}$ be the space of simple maps corresponding to the $v$-th component $\Si_v$ of $\Si$ in (\ref{fplogSetUp_e}). An element of $\cM^\star_{g_v,\mfs_v}(X,D,A_v,\nu)_{I_v}$ is the equivalence class of a tuple 
$$
f_v=\big(\phi_v=\phi|_{\Si_v},u_v,\ze_v=(\ze_{v,i})_{i\in I_v},C_v=(\Si_v,\mfj_v,\vec{z}_v\cup {q}_v)\big)
$$
where $\vec{z}_v$ are the marked points on $\Si_v$, $q_v\!=\!\{q_{\uvec{e}}\}_{\uvec{e}\in \uvec{\E}_v}$ is the set of nodal points on $\Si_v$, $\mfs_v$ is the set of contact orders at $\vec{z}_v\cup {q}_{v}$,  $u_v$ is a map into $D_{I_v}$ satisfying 
$$
\dbar u_v (x)= \nu(\phi_v(x),u_v(x)),\quad \forall~x\!\in\!\Si_v,
$$
and $\ze_v\!=\!(\ze_{v,i})_{i\in I_v}$ is a meromorphic section of $\cN_XD_{I_v}$ with zeros/poles of orders determined by $\mfs_v$ at $\vec{z}_v\cup {q}_{v}$.
Locally around any $f$, we will fix a random ordering $\vec{q}_v$ of $q_v$ that we will forget at the end. If $\Si_v$ is not a bubble component, then 
$$
\cM^\star_{g_v,\mfs_v}(X,D,A_v,\nu)_{I_v}=\cM_{g_v,\mfs_v}(X,D,A_v,\nu)_{I_v};
$$ 
otherwise, $\Si_v$ is a sphere, $\phi_v$ is a constant map (that we will drop from the notation), $\dbar u_v=0$, and by (\ref{g0MI_e})
$$
\cM^\star_{g_v,\mfs_v}(X,D,A_v,\nu)_{I_v}=\cM^{\star}_{0,\ov\mfs_v}(D_{I_v},\partial D_{I_v},A_v).
$$
For each $\uvec{e}\in \uvec{\E}_v$, let
$$
\tn{ev}_{\uvec{e}}\colon \cM_{g_v,\mfs_v}(X,D,A_v,\nu)_{I_v}\lra D_{I_{e}}, \quad [u_v,\ze_v,C_v]\lra u_v(q_{\uvec{e}}),
$$
denote the evaluation map at the nodal point $q_{\uvec{e}}$. Recall that for $I\!=\!\eset$ the convention is $D_I=X$. Let 
$$
\tn{ev}_\E\equiv \prod_{e\in \E}(\tn{ev}_{\uvec{e}}\times \tn{ev}_{\scz\ucev{e}})\colon \prod_{v\in \V} \cM_{g_v,\mfs_v}(X,D,A_v,\nu)_{I_v}\lra \prod_{e\in \E} (D_{I_e})^2
$$
denote the overall evaluation maps at the nodal points.
Then the the fiber product space in (\ref{FP_e}) is
\bEq{PlogAsFiber_e}
\cM^{\tn{plog}}_{g,\mfs}(X,D,A,\nu)_\Gamma=\times_{v\in \V} \cM_{g_v,\mfs_v}(X,D,A_v,\nu)_{I_v}= \tn{ev}_\E^{-1}\big(\prod_{e\in \E} \De_e\big) 
\eEq
where 
$$
\De_e\subset D_{I_e}\times D_{I_e},\qquad \forall~e\!\in\!\E,
$$
is the diagonal subspace. The obstruction map in (\ref{group-element_e}) is the map (\ref{FP_e}) from this fiber product into the obstruction group $\mc{G}$ and
$$
\cM_{g,\mfs}(X,D,A,\nu)_\Gamma=\tn{ob}_{\Gamma}^{-1}(1).
$$
Let 
\bEq{StarProduct_e}
\prod^\star_{v\in \V} \cM^\star_{g_v,\mfs_v}(X,D,A_v,\nu)_{I_v}\subset \prod_{v\in \V} \cM^\star_{g_v,\mfs_v}(X,D,A_v,\nu)_{I_v}
\eEq
be the subset of tuples where the images of every two non-constant bubble components in $X$ are distinct. Also, let
\bEq{PlogAsFiber_e-star}
\cM^{\tn{plog},\star}_{g,\mfs}(X,D,A,\nu)_\Gamma=\cM^{\tn{plog}}_{g,\mfs}(X,D,A,\nu)_\Gamma\cap \prod^\star_{v\in \V} \cM^\star_{g_v,\mfs_v}(X,D,A_v,\nu)_{I_v}.
\eEq
By (\ref{cHgk_e}) and Definition~\ref{UnivFamily_dfn}, the map $\phi$ has image in a product of universal families
\bEq{UFamily_e2}
\bigg(\pi_v\colon\ov{\mf{U}}_{g_v,k_v+\ell_v}\lra \ov{\mf{M}}_{g_v,k_v+\ell_v}, \mf{z}_v\cup \mf{q}_v\bigg)_{v\in\V_{\tn{p}}}.
\eEq
The restriction of a perturbation term $\nu$ (or $\nu_{\log}$) in $\cH_{g,k}(X,D)$ to $\ov{\mf{U}}_{g_v,|z_v|+|q_v|}$ defines an element of 
$\cH_{g_v,|z_v|+|q_v|}(X,D)$. Furthermore, recall from (\ref{NTpair_e}) that if $\nu\!\in\!\cH_{g,k}(X,D)_{\cR,J}$, then the restriction of that to $\mf{U}_{g_v,|z_v|+|q_v|}$ and $D_{I_v}$ is made of components 
\bEq{Nuv_e}
\nu_v\in \cH_{g_v,|z_v|+|q_v|}(D_{I_v},\partial D_{I_v})_{\cR_{I_v},J_{I_v}}
\eEq
and 
\bEq{Thv_e}
\theta_v=(\theta_{v,i})_{i\in I_v}\in\Theta_{g_v,|z_v|+|q_v|}(D_{I_v})
\eEq
such that 
\bEq{RecallDef_e}
\dbar u_v =(\phi_v,u_v)^*\nu_v\quad \tn{and}\quad \dbar_{u_v^*\cN_XD_i} \ze_{I_v,i} =\theta_{v,i}\ze_{I_v,i}\qquad \forall~i\!\in\!I_v.
\eEq

\bLm{Dim_lmm} With notation as above, if 
\bEn
\item $J\!\in\!\tn{AK}^{I_v}(X,D)_\cR$ in the sense of Proposition~\ref{Trans2_prp}.(2) for all $v\in \V_{\tn{b}}$,
\item $\nu\!\in\!\cH^{I_v}_{g_v,|z_v|+|q_v|}(X,D)_{\cR,J}$ in the sense of Proposition~\ref{Trans2_prp}.(1) for all $v\in \V_{\tn{p}}$,
\item the map $\tn{ev}_\E$ restricted to $\prod^\star_{v\in \V} \cM^\star_{g_v,\mfs_v}(X,D,A_v,\nu)_{I_v}$ is transverse, and
\item\label{ObGamma_l} the map $\tn{ob}_{\Gamma}$ restricted to $\cM^{\tn{plog},\star}_{g,\mfs}(X,D,A,\nu)_\Gamma$ is transverse,
\eEn
then $\cM^\star_{g,\mfs}(X,D,A,\nu)_\Gamma$ is a naturally oriented smooth manifold of the real dimension (\ref{dGamma_e}).
\eLm
\vskip.1in
\bPf
By Proposition~\ref{Trans2_prp}, under the first two conditions, each $\cM^\star_{g_v,\mfs_v}(X,D,A_v,\nu)_{I_v}$ is a naturally oriented smooth manifold of the real dimension 
$$
2\big(c_1^{TX(-\log D)}(A_v)+(n-3)(1-g_v)+k_v+\ell_v-|I_v|\big), \quad\tn{where}~~k_v=|\vec{z}_v|,~\ell_v=|q_v|.
$$
By the third bullet, since each $D_I$ is oriented, the fiber product space $\cM^{\tn{plog},\star}_{g,\mfs}(X,D,A,\nu)_\Gamma$ is a naturally oriented smooth manifold of the real dimension 
$$
\aligned
&2\bigg(\sum_{v\in \V}\big(c_1^{TX(-\log D)}(A_v)+(n-3)(1-g_v)+k_v+\ell_v-|I_v|\big)-\sum_{e\in \E} (n-|I_e|)\bigg)=\\
&2\bigg(c_1^{TX(-\log D)}(A)+(n-3)(1-g)+k-|\E|-\sum_{v\in \V}|I_v|+\sum_{e\in \E} |I_e|\bigg).
\endaligned
$$
By (\ref{DtoT_e}), 
$$
\dim_\R \K_\R(\Gamma)-\tn{dim}_\C(\mc{G})=|\E|+\sum_{v\in \V}|I_v|-\sum_{e\in \E}|I_e|.
$$
Therefore, by the fourth bullet, and since $\mc{G}$ is a complex manifold, $\cM^\star_{g,\mfs}(X,D,A,\nu)_\Gamma$ is a naturally oriented smooth orbifold of the real dimension
$$
\aligned
&2\bigg(c_1^{TX(-\log D)}(A)+(n-3)(1-g)+k-|\E|-\sum_{v\in \V}|I_v|+\sum_{e\in \E} |I_e|- \tn{dim}_\C(\mc{G})\bigg)=\\
&2\bigg(c_1^{TX(-\log D)}(A)+(n-3)(1-g)+k- \dim_\R \K_\R(\Gamma)\bigg).
\endaligned
$$
\ePf
\noindent
Recall that, by the first condition in Definition~\ref{LogMap_dfn},  
$$
 \dim_\R \K_\R(\Gamma)>0
 $$ 
 unless $\Gamma$ is the trivial one vertex graph ($\V\!=\!\{v\}, \E\!=\!\eset$ with $I_v\!=\!\eset$) which corresponds to the main stratum. In the classical case (no $D$), the map $\vr$ in (\ref{DtoT_e}) is the trivial map $\Z^\E\lra 0$. Therefore, $\dim_\R \K_\R(\Gamma)=|\E|$ is the number of the nodes. In the logarithmic case, there are configurations with arbitrary large number of nodes and $\dim_\R \K_\R(\Gamma)=1$; see \cite[Ex~3.12]{FT1}.
 
\newtheorem*{PFTransGamma_prp}{Proof of Theorem~\ref{TransGamma_thm}}
\begin{PFTransGamma_prp}
With modifications as in Sections~\ref{Generic_ss}~and~\ref{NGSmooth_ss}, proof of this proposition is similar to the proof of \cite[Prp~3.16]{RT}, \cite[Thm~6.2.6]{MS2}, and \cite[Prp~4.3]{Z2}. The main difference is the extra evaluation-type map $\tn{ob}_\Gamma$ that needs to be transversed as in Lemma~\ref{Dim_lmm}.\ref{ObGamma_l}.\\

\noindent
For each $v\!\in\!\V_{\tn{p}}$, let $\ov{\mf{U}}_{g_v,k_v+\ell_v}$, $\nu_v$, and $\theta_v$ be as in (\ref{UFamily_e2}), (\ref{Nuv_e}), and (\ref{Thv_e}), respectively. For each $v\!\in\!\V_{\tn{b}}$, since $\Si_v\!=\!\P^1$, if $C_v$ is stable, let  $\ov{\mf{U}}_{0,k_v+\ell_v}$ denote the universal curve, and if $C_v$ is not stable, i.e., if it is a $\P^1$ with less than $3$ points, then let $\ov{\mf{U}}_{0,k_v+\ell_v}\!=\!C_v$ and $\ov{\mf{M}}_{0,k_v+\ell_v}$ be just a point. For $v\!\in\!\V_{\tn{b}}$, $\nu_v\!=\!\phi^*\nu|_{\Si_v}$ is zero and we are dealing with log $J$-holomorphic curves.\\

\noindent 
For each $v\!\in\!\V$, fix a local family $\pi_v|_{\cC_v}\colon \cC_v\lra \cB_v$ around $C_v$ as in (\ref{LocalUnivFam_e}) and a smooth trivialization $\varphi_v$ of that as in (\ref{ST_e}).\\

\noindent
For $v\!\in\!\V_{\tn{p}}$, let
$$
W_v\equiv W^{\ell,p}_{\ov\mfs_v,A_v}\big((\cC_v,\{\mf{z}_v\cup \mf{q}_v\}|_{\cB_v}),(D_{I_v},\partial D_{I_v})\big)_{\varphi_v}\!\times\!\cH^m_{g_v,k_v+\ell_v}(D_{I_v},\partial D_{I_v})_{[\om]}\times_J \Theta^m_{g_v,k_v+\ell_v}(D_{I_v}),
$$
where 
$$
\ov\mfs_{v}\!=\!\big(s_{v;a}\!=\!(s_{v;ai})_{i\in [N]-I_v}\big)_{a\in [k_v+\ell_v]}\!\in\!(\N^{[N]-I_v})^{k_v+\ell_v}
$$ 
is as in (\ref{I-embedding}),
be the corresponding configuration space as in (\ref{WIConf_e}). This time, $(\cR,J)$ is not fixed and can change as well. Thus, the notation $\times_J$ means fiber product over $\tn{AK}(X,D)_{[\om]}$.  For each $v\!\in\!\V_{\tn{b}}$, the configuration space is 
$$
W_v\equiv W^{\ell,p}_{\ov\mfs_{v},A_v}\big((\cC_v,\{\mf{z}_v\cup \mf{q}_v\}|_{\cB_v}),(D_{I_v},\partial D_{I_v})\big)_{\varphi_v}
\!\times\!\tn{AK}^m(X,D)_{[\om]}
$$
as in (\ref{UM_e0}) with $(D_{I_v},\partial D_{I_v})$ in place of $(X,D)$.
The evaluation map $\tn{ev}_\E$ at the nodal points extends to the product $\prod_{v\in \V} W_v$. For $v\!\in\!\V$, let 
$$
\cE_v\equiv \cE^{\ell,p}_{\ov\mfs_v,A_v}\big((\cC_v,\{\mf{z}_v\cup \mf{q}_v\}|_{\cB_v}),(D_{I_v},\partial D_{I_v})\big)\lra W_v
$$
be the Banach bundles in (\ref{BVB_e3}).
Let 
$$
W_\Gamma\!\equiv \!\times_{v\in \V} W_v= \tn{ev}_\E^{-1}\big(\prod_{e\in \E} \De_e\big),
$$
$\pi_v\colon W_\Gamma\lra W_v$ denote the projection map into the $v$-th component, and
$$
E_\Gamma\equiv \bigoplus_{v\in \V} \pi_v^* \cE_v  \lra W_\Gamma
$$
denote the obstruction Banach bundle.
The tangent space of  $W_\Gamma$ has the form 
$$
\aligned
T_fW_\Gamma& \cong T^{\tn{ver}}_fW_\Gamma \oplus \bigoplus_{v\in \V} T_{C_v} \cB_v\oplus T_{(\cR,J)}\tn{AK}^m(X,D)_{[\om]} \oplus\\ 
&\bigoplus_{v\in \V_{\tn{p}}}\big(T_{\nu_v} \cH^m_{g_v,k_v+\ell_v}(D_{I_v},\partial D_{I_v})_{\cR_{I_v},J_{I_v}}\oplus T_{\theta_v} \Theta^m_{g_v,k_v+\ell_v}\big),
\endaligned
$$
where $T^{\tn{ver}}_fW_\Gamma$ has the following description. By (\ref{T_u-e}),
$$
T_{u_v}W^{\ell,p}_{\ov\mfs_v,A_v}\big((\Si_v,\vec{z}_v\cup q_v),(D_{I_v},\partial D_{I_v})\big)=W^{\ell,p}\big(\Si_v,u_v^*TD_{I_v}(-\log D_{I_v})\big).
$$
For every $\xi\!\in\!W^{k,p}\big(\Si_v,u_v^*TD_{I_v}(-\log D_{I_v})\big)$, the section $\iota(\xi)\!\in\!C^0\big(\Si_v,u_v^*TD_{I_v}\big)$ defined via (\ref{LogToTX_e}) satisfies 
$$
\iota(\xi)(q_{\uvec{e}})\in TD_{I_e}\qquad \forall~\uvec{e}\!\in\!\uvec{\E}_v.
$$
Then 
\bEq{NMxi_e}
\aligned
T^{\tn{ver}}_fW_\Gamma= \bigg\{(\xi_v)_{v\in \V} \in \bigoplus &T_{u_v}W^{\ell,p}_{\ov\mfs_v,A_v}\big((\Si_v,\vec{z}_v\cup q_v),(D_{I_v},\partial D_{I_v})\big)\colon \\
& \iota(\xi_v)(q_{\uvec{e}})=\iota(\xi_{v'})(q_{\scz\ucev{e}}) \quad \forall~v,v'\!\in\!\V,~\uvec{e}\!\in\!\E_{v,v'}\bigg\}.
\endaligned
\eEq
Summing the maps (\ref{GenPMap_e}) over all principal components we get
\bEq{GaPMap_e}
P=\prod_{v\in \V_{\tn{p}}} P_v \colon W_{\Gamma}\lra  \tn{Pic}^0_{\Gamma}\equiv \prod_{v\in \V_{\tn{p}}} \tn{Pic}^0(\mf{U}_{g_v,k_v+\ell_v})^{I_v}.
\eEq
Let $\wt{\mf{M}}_{g,\mfs}^{\tn{plog}}(X, D,A)_\Gamma $ denote the $(0\times \cO_\Gamma)$-level set  of 
\bEq{CRPGamma_e}
(\dbar^{\log}-\nu_{\log})\times P \colon W_\Gamma \lra  E_\Gamma \times   \tn{Pic}^0_{\Gamma},
\eEq
where $\cO_\Gamma\!\equiv\! \prod_{v\in \V_{\tn{p}}} \cO^{I_v}$. The universal pre-log moduli space $\mf{M}_{g,\mfs}^{\tn{plog}}(X, D,A)_\Gamma $ is the quotient of $\wt{\mf{M}}_{g,\mfs}^{\tn{plog}}(X, D,A)_\Gamma $ by the automorphism group. The latter is the product of automorphism groups of the bubble components $\Si_v\!\cong\! \P^1$ with $k_v+\ell_v <3$.  \\

\noindent
\textbf{Claim 1}.  Restricted to the subset of simple maps $\wt{\mf{M}}_{g,\mfs}^{\tn{plog},\star}(X, D,A)_\Gamma $, $0\!\times\! \cO_\Gamma$ is a regular value of (\ref{CRPGamma_e}).\\

\noindent
\textbf{Proof}. For every $(f,J,\nu)\!\in\! \wt{\mf{M}}_{g,\mfs}^{\tn{plog}}(X, D,A)_\Gamma $, the linearization $\tn{D}^{\log}_{f,J,\nu}(\{\dbar-\nu\}\!\times\! P)$ of (\ref{CRPGamma_e}) is the direct sum linearization map
\bEq{BigDf_e}
\tn{D}_{f,J,\nu}^{\log}\big(\{\dbar-\nu\}\!\times\! P\big)\!\equiv\! 
\bigoplus_{v\in \V_{\tn{b}}} \tn{D}^{\log}_{u_v,J} \dbar \oplus \bigoplus_{v\in \V_{\tn{p}}} \tn{D}_{u_v,J,\nu_v}^{\log}\big(\{\dbar-\nu_v\}\!\times\! P_v\big) 
\eEq
with summands as in (\ref{LCRNu_e0}) and (\ref{DCRNuP_e}). By the proof of Proposition~\ref{Trans2_prp} part (1), fixing $J$, for each $v\!\in\!\V_{\tn{p}}$, $\tn{D}_{u_v,\nu_v}^{\log}\big(\{\dbar-\nu_v\}\times P_v\big)$
is surjective onto 
$$
\cE_v\oplus T_{\cO^{I_v}}\tn{Pic}^0(\mf{U}_{g_v,k_v+\ell_v})^{I_v}.
$$
Furthermore, transversality can be achieved by sections supported away from the nodes. Therefore, for any fixed $J$, the direct sum 
$$
\bigoplus_{v\in \V_{\tn{p}}} \tn{D}_{u_v,\nu_v}^{\log}\big(\{\dbar-\nu_v\}\times P_v\big) 
$$
is surjective onto 
$$
\bigoplus_{v\in \V_{\tn{p}}} \cE_v\oplus T_{\cO_\Gamma} \tn{Pic}^0_{\Gamma}.
$$
By the proof of Proposition~\ref{Trans2_prp} part (2), for each $v\!\in\!\V_{\tn{b}}$, $\tn{D}_{u_v,J}^{\log} \dbar$
is surjective onto $\cE_v$. 
Furthermore, if the bubble components are all simple with mutually different images, as in the proof of \cite[Prp~6.2.7]{MS2}, by Proposition~\ref{Trans2_prp} part (2) and Lemma~\ref{Ydef_lmm}, transversality can be achieved simultaneously by deformations of $J$ that are supported along different images of these components. We conclude that (\ref{BigDf_e}) is surjective along the subset of simple maps.\qed
\\

\noindent
\textbf{Claim 2.} \textit{The map $\tn{ob}_{\Gamma}\colon \wt{\mf{M}}^{\tn{plog},\star}_{g,\mfs}(X,D,A)_\Gamma\lra \mc{G}(\Gamma)$ is transverse.}\\

\noindent
\textbf{Proof}. Fix a tuple
$$
f=\bigg(\big(c,u_v,\ze_v=(\ze_{v,i})_{i\in I_v},C_v=(\Si_v,\mfj_v,\vec{z}_v\cup {q}_v)\big), J, \nu\bigg)\in \tn{ob}_{\Gamma}^{-1}(1).
$$
By definition, we can choose local holomorphic coordinates $w_{\uvec{e}}$ around each nodal point $q_{\uvec{e}}\!\in\!\Si_v$ and representatives $\ze_v$ such that 
\bEq{CofG_e2}
\eta_{\uvec{e},i}/\eta_{\scz\ucev{e},i}=1, \qquad \forall~e\!\in\!\E,~i\!\in\!I_e;
\eEq
see (\ref{CofG_e}).
Let $q_{\uvec{e}}$ be a nodal point on $\Si_v$ connecting that to $\Si_{v'}$. 
By Lemma~\ref{TofCC_lmm}, we may assume that either $v\!\in\!\V_{\tn{p}}$ or $v\!\in\!\V_{\tn{b}}$ and $A_v\!\neq\!0$. Equally, one may use the method of proof of \cite[Thm~6.3.1]{MS2} in \cite[p.155]{MS2} to address the nodes connecting two ghost bubbles. Choose any $i\!\in\!I_e$. \\

\noindent
(i) If $v\!\in\!\V_{\tn{p}}$ and $i\!\in\!I_v$, we have a meromorphic section $\ze_{v,i}$ included in $f$ satisfying 
$$
\dbar_{u_v^*\cN_XD_i} \ze_{v,i}=\theta_{v,i} \ze_{v,i};
$$
see (\ref{RecallDef_e}) above. Recall from (\ref{Expansion_e2}) that 
$$
\ze_{v,i}(w_{\uvec{e}})=\wt\eta_{\uvec{e},i}(w_{\uvec{e}})w_{\uvec{e}}^{s_{\uvec{e},i}}\in u_v^* \cN_XD_i,
$$
such that 
$$
0\!\neq\!\eta_{\uvec{e},i}\equiv\wt\eta_{\uvec{e},i}(0)\!\in\! \cN_XD_i|_{u_v(q_{\uvec{e}})}.
$$
Let
$$
\beta\colon \Si_v \lra \C
$$
be a smooth function that is supported in a neighborhood of $q_{\uvec{e}}$, is constant in a smaller neighborhood of $q_{\uvec{e}}$, and satisfies 
\bEq{viDirection}
\frac{\nd}{\nd t} (e^{t\beta (q_{\uvec{e}})}\eta_{\uvec{e},i}/\eta_{\scz\ucev{e},i})|_{t=0}=1.
\eEq
Let 
$$
(f_t)_{t\in [0,1]}\in  \wt{\mf{M}}^{\tn{plog},\star}_{g,\mfs}(X,D,A,\nu)_\Gamma  
$$
be the path obtained by deforming $\ze_{v,i}$ to $e^{t\beta}\ze_{v,i}$ and $\theta_{v,i}$ to $\theta_{v,i}+t \dbar \beta$. Then, by (\ref{viDirection}), 
\bEq{1ei_e}
\frac{\nd}{\nd t} \tn{ob}_{\Gamma}(f_t)= [1_{e,i}] \in T_{1}\mc{G}_\Gamma= \tn{coker} (\vr_{\C}),
\eEq
where $[1_{e,i}]$ denotes the image of $1_{e,i}\!\in\!\bigoplus_{e\in \E} \C^{I_e}$ in $\tn{coker} (\vr_{\C})$.\\

\noindent
(ii) If $v\!\in\!\V_{\tn{p}}$ and $i\!\in\!I_e-I_v$, near $q_{\uvec{e}}$, $u_v$  has the local form 
$$
u_v(w_{\uvec{e}})\cong (u_{v,i}(w_{\uvec{e}}),\wt\eta_{\uvec{e},i}(w_{\uvec{e}})w_{\uvec{e}}^{s_{\uvec{e},i}})\in \cN_XD_i,
$$
such that 
$$
0\!\neq\!\eta_{\uvec{e},i}\equiv\wt\eta_{\uvec{e},i}(0)\!\in\! \cN_XD_i|_{u_v(q_{\uvec{e}})};
$$
see (\ref{Expansion_e}). Therefore, instead, we can deform $u_v$ by deforming $\wt\eta_{\uvec{e},i}(w_{\uvec{e}})$ to  $e^{t\beta}\wt\eta_{\uvec{e},i}(w_{\uvec{e}})$ as above and get the same conclusion as in (\ref{1ei_e}). \\

\noindent
(iii) If $v\!\in\!\V_{\tn{b}}$ and $i\!\in\!I_e-I_v$, we need to consider a deformation of $J$. By Lemma~\ref{Ydef_lmm}, for any open set $U\!\subset\!\Si-\{\vec{z}_v\cup q_{v}\}$, there exists  
$$
\xi\!\in\!W^{\ell,p}\big(\Si_v, u_v^*TD_{I_v}(- \log \partial D_{I_v})\big)\quad \tn{and}\quad Y\!\in\! T_{J_{I_v}}\tn{AK}^m(D_{I_v},\partial D_{I_v})_{\cR_{I_v}}
$$ 
such that 
$$
\aligned
&\xi(q_{\uvec{e}})=0 \oplus 1_{e,i} \in T_{u_v(q_{\uvec{e}})}D_{I_v}(-\log \partial D_{I_v})\cong T_{u_v(q_{\uvec{e}})}D_{I_e}(-\log \partial D_{I_e})\oplus \C^{I_e-I_v},\\
&\xi(z_{a})\!=\!0\quad \forall~z_a\!\in\! \vec{z}_{v},
\quad \xi(q_{\uvec{e}'})\!=\!0\quad \forall~\uvec{e}'\!\in\! \uvec{\E}_{v}-\uvec{e},\quad \tn{Supp}(Y|_{\Si_v})\!\subset\! U
\quad D_{u_v}\dbar^{\log} \xi\!+\! \frac{1}{2}\,Y \circ \nd u_v \circ \mfj =0 .
\endaligned
$$
By trivial extension of $\xi$ to other components and horizontal extension of $Y$ to a deformation of $J$ on the entire $X$, such a pair $(\xi,Y)$ defines a tangent vector in 
$$
T_f\wt{\mf{M}}^{\tn{plog},\star}_{g,\mfs}(X,D,A,\nu)_\Gamma
$$
such that 
$$
D_f\tn{ob}_\Gamma (\xi,Y)= [1_{e,i}] \in T_{1}\mc{G}_\Gamma.
$$ 
\vskip.1in

\noindent
(iv) If $v\!\in\!\V_{\tn{b}}$ and $i\!\in\!I_v$, we need to consider a deformation of $J$ in the normal direction to $D_i$. That is, we need a deformation of $\cR$ by deforming the connection $\na^{(i)}$ in Definition~\ref{omRegV_dfn}. This is the only step where we need to deform $\cR$. Deformation of $\na^{(i)}$ by a $1$-form results in a deformation of the corresponding $\dbar$-operator $\dbar_{\cN_XD_i}$ on $\cN_XD_i$ in \cite[Lmm~2.1]{FT1}. Deformation of $\dbar_{\cN_XD_i}$, then, yields a deformation of the associated almost complex structure on the total space of $\cN_XD_i$; see \cite[Lmm~2.2]{Z1} or \cite[Sec~2.1]{FT1}. In other words, the isomorphism
$$
T\cN_XD_i \cong \pi_i^*TD_i \oplus \pi^*\cN_XD_i,
$$
and thus the construction of $J$ on $T\cN_XD_i$ via $J|_{TD_i}$ and $\mfi_i$ on $\cN_XD_i$, depends on $\na^{(i)}$. Deforming the latter results in a deformation of the former. A deformation 
$$
u^*\dbar_{\cN_X D_i} \lra u^*\dbar_{\cN_X D_i} + t\beta 
$$ 
(supported on an open set of $\Si_v$ whose image in $X$ is disjoint from the image of the rest of non-trivial bubble components) as in (i) such that (\ref{viDirection}) holds yields a path $(f_t)_{t\in [0,1]}$ as in (i)  such that (\ref{1ei_e}) holds.\qed 
\vskip.1in

\noindent
Finally, by Claims~1 and 2, and Implicit-Function Theorem, the universal moduli space 
$$
\mf{M}^{\star}_{g,\mfs}(X,D,A)_\Gamma
$$
is a separable $C^{m-\ell}$-smooth Banach manifold. Then by the Sard-Smale Theorem, the set of regular values  $\cH\Theta^{\tn{reg,m}}_{g,k}(X,D)_{[\om]}$ of the projection map
$$
\tn{proj}\colon \mf{M}_{g,\mfs}^\star(X,D,A)_I\lra \prod_{v\in \V_{\tn{p}}}\cH^m_{g_v,k_v+\ell_v}(D_{I_v},\partial D_{I_v})_{[\om]}\times_J \Theta^m_{g_v,k_v+\ell_v}(D_{I_v})
$$
is Baire set of second category. By definition and Lemma~\ref{Dim_lmm}, for every $(\om',\cR,J,\nu)\!\in\!\cH_{g,k}(X,D)_{[\om]}$ such that 
$$
\tn{Res}(\om',\cR,J,\nu)=\bigg(\om',\cR,J,\big((\nu_v,\theta_v)=\nu|_{\Si_v}\big)_{v\in \V_{\tn{p}}}\bigg) \in  \cH\Theta^{\tn{reg,m}}_{g,k}(X,D)_{[\om]},
$$ 
the stratum
$$
\mf{M}_{g,\mfs}(X,D,A,\nu)_I=\tn{proj}^{-1}\bigg(\om',\cR,J,\big((\nu_v,\theta_v)=\nu|_{\Si_v}\big)_{v\in \V_{\tn{p}}}\bigg) 
$$
is a naturally oriented smooth manifold of the expected dimension (\ref{dGamma_e}). With an argument similar to the Taubes' trick in the proof of \cite[Thm~3.1.6(ii)]{MS2}, we conclude that the subset of smooth perturbations 
$$
\cH^{\Gamma}_{g,k}(X,D)_{[\om]}=\tn{Res}^{-1}\big(\cH\Theta^{\tn{reg,m}}_{g,k}(X,D)_{[\om]}\big) \cap \cH^{\Gamma}_{g,k}(X,D)_{[\om]}
$$
satisfies Theomre~\ref{TransGamma_thm}.(1). If we restrict this proof to genus $0$ $J$-holomorphic log curves and the resulting set of regular values in $\tn{AK}(X,D)_{[\om]}$, we get Theorem~\ref{TransGamma_thm}.(2).
\qed

\end{PFTransGamma_prp}

%--------------------------------------------------------
\subsection{Genus zero multiple-cover maps}\label{RMCMap_ss} 
\noindent
A main step in proving Proposition~\ref{MC_pr} (for arbitrary $D$) is to address the transversality issue at multiple-cover log $J$-holomorphic spheres. In this section we show that under the positivity/semi-positivity conditions in Definition~\ref{NCsemipos_dfn}, multiple-cover log spheres do not happen in families of larger than the expected dimension.\\

\noindent
Suppose 
$$
[u,\P^1,z_1,\ldots,z_\ell]\!\in\!\cM^\star_{0,\mft}(X,D,B)
$$
with $\mft\!=\!(t_1,\ldots,t_\ell)$ and $t_a\!\neq\! 0$ for all $a\!\in\![\ell]$. All the marked points are contact points; therefore, none of them can be ignored (as in the classical case) to decrease the expected dimension. The other cases can be reduced to this case by ignoring those marked points that have trivial contact with $D$. Every such point has a finite number of pre-images in any multiple-cover of $u$. Let  $d$ and $k_1,\ldots,k_\ell\!\leq\! d$ be   positive integers and set $k\!=\!k_1+\ldots+k_\ell$. For each $a\!\in\![\ell]$ let  
$$
\al_{1,a}+\cdots+\al_{a,k_a}=d
$$
be an ordered partition of $d$ into a sum of $k_a$ positive numbers. We are interested in those tuples 
\bEq{Partition_e}
\al\equiv \big(\al_{a,b}\big)_{a\in [\ell], b\in [k_a]}
\eEq
such that there exists a degree $d$  covering map 
\bEq{h_e}
h\colon \P^1\lra \P^1
\eEq
satisfying
$$
h^{-1}(z_a)=\{z_{ab}\}_{b\in [k_a]}, \qquad \ord_{z_{ab}}h=\al_{a,b}.
$$ 
If $h$ is such a covering map, the  $k$-marked degree $d$ map
$$
f_h=\big(h,\P^1, (z_{ab})_{a\in [\ell], b\in [k_a]}\big)
$$
(with the lexicographic order on the marked points $z_{ab}$) defines a point of the log moduli space 
$$
\cM_{0,\al}(\P^1,D_z,[d])
$$
where $D_z\!=\!\{z_1,\ldots,z_\ell\}$ and we are treating $\al\!\in\! \N^k$ as the tangency order data with $D_z$ at those $k$ points. By (\ref{Obs_f2}), if $\cM_{0,\al}(\P^1,D_z,[d])$ is non-empty, it is a smooth manifold of the expected dimension
\bEq{MC-dim_e}
d_{\tn{fiber}}\equiv\tn{dim}_\C~\cM_{0,\al}(\P^1,D_z,[d])= (d-1)(2-\ell)+k-\ell.
\eEq
Also, the  $k$-marked degree $A=dB$ composition map
\bEq{MCu_e}
\Big(u\circ h,\P^1, (z_{ab})_{a\in [\ell], b\in [k_a]}\Big)
\eEq
defines a point of the log moduli space $\cM_{0,\mfs}(X,D,A)$, with 
$$
\mfs=(s_{ab})_{a\in [\ell], b\in [k_a]}, \quad  \N^N\ni s_{ab}= \al_{a,b}s_a. 
$$
Let 
$$\cM^\al_{0,\mfs}(X,D,A)\!\subset\! \cM_{0,\mfs}(X,D,A)$$ 
denote the subspace of multiple-cover maps of type $\al$. There is a projection map 
\bEq{alTostar_e}
\cM^\al_{0,\mfs}(X,D,A)\lra \cM^\star_{0,\mft}(X,D, B)
\eEq
whose fiber over $[u,\P^1,z_1,\ldots,z_\ell]$ is $\cM_{0,\al}(\P^1,D_z,[d])$.  We have 
$$
\aligned
d_{\tn{down}}\equiv \tn{exp-dim}_\C~\cM^\star_{0,\mft}(X,D,B)&=c_1^{TX(-\log D)}(B)+n-3+\ell, \\
d_{\tn{up}}\equiv\tn{exp-dim}_\C~\cM^\star_{0,\mfs}(X,D,A)&=c_1^{TX(-\log D)}(A)+n-3+k,
\endaligned
$$
In order for the image of $\cM^\al_{0,\mfs}(X,D,A)$ under $\tn{ev}$ to have a smaller (resp. smaller or equal) dimension than the image of the main stratum $\cM^\star_{0,\mfs}(X,D,A)$, whenever $d_{\tn{fiber}},d_{\tn{down}}\geq 0$, we need $d_{\tn{down}}\leq d_{\tn{up}}$, i.e.
\bEq{SP-cond_e}
c_1^{TX(-\log D)}(B)+n-3+\ell \geq 0  ~~\Rightarrow ~~(d-1)\, c_1^{TX(-\log D)}(B)+k-\ell\!>\! 0~~(\tn{resp.}\,\geq\!0).
\eEq
In other words, we want to \textbf{avoid} a situation where $d_{\tn{fiber}}\geq 0$, $B\!\in\!\pi_2(X)$ with $\om(B)>0$, and 
$$
3-n-\ell \leq  c_1^{TX(-\log D)}(B) \leq  \frac{\ell-k}{d-1}~~\bigg(\tn{resp. } < \frac{\ell-k}{d-1}\bigg).
$$
For $\ell=0,1,2$, the condition $d_{\tn{fiber}}\geq 0$ automatically holds.
For $\ell>2$, the condition $d_{\tn{fiber}}\geq 0$ implies
$$
\displaystyle{\frac{\ell-k}{d-1}}\leq 2-\ell.
$$
Therefore, (\ref{SP-cond_e}) can be replaced with the stronger requirement
\bEq{GF_e}
\begin{cases}
c_1^{TX(-\log D)}(B) \notin [3-n-\ell,0]~~\big(\tn{resp. } [3-n-\ell,0)\big) & \tn{if}\quad \ell=0,1,2\\
c_1^{TX(-\log D)}(B) \notin [3-n-\ell,2-\ell]~~\big(\tn{resp. } [3-n-\ell,2-\ell)\big) & \tn{if}\quad \ell>2
\end{cases}
\eEq
More generally, we will consider maps with image in $D_I$. Then $n$ should be replaced with $n-|I|$.
For $N\!>\!1$ and $(I,\ell)\!\neq\! (\emptyset,0)$, we will further need  $d_{\tn{fiber}}+d_{\tn{down}}\leq d_{\tn{up}}$. The latter is equivalent to $c_1^{TX(-\log D)}(B)\geq 2-\ell$. This explains Definition~\ref{NCsemipos_dfnN}. The following lemma summarizes the outcome of these calculations.

\bLm{MCLemma}
With $D_I$ and $\partial D_I$ in place of $(X,D)$ in (\ref{alTostar_e}), if $[X,D,\om]$ is semi-positive (resp. positive) in the sense of Definition~\ref{NCsemipos_dfn} and $\ell\leq 2$, then 
\bEq{IfSimple_e}
\tn{exp-dim}\; \cM^\star_{0,\mft}(D_I,\partial D_I, B)\geq 0
\eEq
implies 
$$
\tn{exp-dim}\, \cM^\al_{0,\mfs}(D_I,\partial D_I,dB) \geq \tn{exp-dim}\, \cM^\star_{0,\mft}(D_I,\partial D_I, B)\qquad (\tn{resp}. >),
$$
for all $d\!\geq\!1$ and $(\al,\mfs,\mft)$ as above.
Furthermore, if $[X,D,\om]$ is strongly-semi-positive in the sense of Definition~\ref{NCsemipos_dfnN}, then 
(\ref{IfSimple_e}) implies 
$$
\tn{exp-dim}\, \cM^\al_{0,\mfs}(D_I,\partial D_I,A) \geq \tn{exp-dim}\, \cM^\star_{0,\mft}(D_I,\partial D_I, B)+ \tn{exp-dim}\, \cM_{0,\al}(\P^1,D_z,[d]).
$$

\eLm

\noindent
Regarding~(\ref{SP-cond_e}), the example below illustrates a non-positive situation where (\ref{SP-cond_e}) does not hold, $\cM^\star_{0,\mfs}(X,D,A)$ is empty, and $\cM^\al_{0,\mfs}(X,D,A)$ is always positive dimensional.
\bEx{P2Q4_ex}
Let $X\!= \!\P^2$, $D$ be a smooth quartic hypersurface, $B\!=\![1]\!\in\! H_2(\P^2,\Z)\!\cong\! \Z$, $\ell\!=\!2$, $k\!=\!3$, $\mft\!=\!(2,2)$, and $\mfs=(2a,2b,2d)$ with $a,b\!>\!0$ and $a\!+\!b\!=\!d$. For generic $J$, 
$$
\cM^\star_{0,\mft}(X,D,B)=\ov\cM^{\log}_{0,\mft}(X,D,B)
$$
is the (zero-dimensional) moduli space of lines with $2$ intersections of order $2$ with $D$ and has $160$ points. Also, $d_{\tn{fiber}}=1$; if $(z_1,z_2)=(0,\infty)$, the map $h$ in (\ref{h_e}) is of the form 
$$
h\big(z\big)=\la\frac{ (z-z_{11})^a(z-z_{12})^b}{(z-z_{21})^d}\quad\tn{for some } \la\in \C^*.
$$
On the other hand, $d_{\tn{up}}= 2-d$; therefore, if $d\!>\!2$, $\cM^\star_{0,\mfs}(X,D,A)$ is empty for generic $J$ while $\cM^\al_{0,\mfs}(X,D,A)$ will always be positive dimensional.\qed
\eEx

\noindent
For $(J,\nu)$-holomorphic curves, every connected cluster $\Si'$ of contracted components is a tree of spheres, with a total of at most $2$ special points at least one of which is a nodal point. Here, by a special point we mean either a marked point or a nodal point connecting the cluster to an irreducible component of $\Si$ outside the cluster \footnote{By Remark~\ref{SubCurve_rmk}, such a cluster defines a log $J$-holomorphic curve $f'\!\in\!\cM_{0,\mfs'}(X,D,A')_{\Gamma'}$ where $\mfs'$ records the contact data at those (one or two) special points (which now act as marked points for $f'$), $\Gamma'$ is the decorated sub-graph of the cluster, and $A'$ is the total homology class of the cluster}. Because of this restriction on the number of special points, it is natural to expect that the condition $\ell\!\leq\!2$ in Lemma~\ref{MCLemma} to be always satisfied. The following lemma shows that this is indeed the case under the Nef condition of Section~\ref{SPLP_ss}.\\

\noindent
For each bubble component $\P^1\!\cong\!\Si_v\!\subset\!\Si'$, the restriction of $f'$ to $\Si_v$ is a tuple 
$$
f_v=\big(u_v,\ze_v,C_v=(\Si_v,\mfj_v,z_v\cup q_v)\big),
$$
where $z_v$ is empty or is the only $1$ marked point allowed on $\Si'$, $q_v\!=\!\{q_{\uvec{e}}\}_{\uvec{e}\in \uvec{\E}_v}$ is the set of nodal points on $\Si_v$, $\mfs_v$ is the set of contact orders at $z_v\cup q_{v}$,  $u_v$ is a $J$-holomorphic map into $D_{I_v}$, and $\ze_v\!=\!(\ze_{v,i})_{i\in I_v}$ is a meromorphic section of $u_v^*\cN_XD_{I_v}$ with zeros/poles of orders determined by $\mfs_v$ at $z_v\cup q_{v}$. Note that while $f'$ has at most  $2$ special points, $q_{v}$ can be arbitrary large. The case we are interested in is when $u_v$ is a multiple-cover map of the form $u_v=\ov{u}_v\circ h$ as in (\ref{MCu_e}), where $\ov{u}_v$ represents the homology class $B_v$ and has contact type $\mft_v$. We say a point $p\in \Si_v$ has \textbf{non-trivial} contact with $D$, or $p$ is a \textbf{contact point} with $D$, if $I_v\neq[N]$ and 
$$
p\in u_v^{-1}\big(\bigcup_{i\in [N]-I_v} D_i \big).
$$
Let $\de_v$ (resp. $\ov\de_{v}$) denote the number of contact points of $u_v$ (resp. $\ov{u}_v$) with $D$.
We say $p\in \Si_v$ is a \textbf{positive} point if there exists $i\!\in\![N]$ such that
$$
\tn{ord}^{i}_{u_v,\ze_v}(p)>0.
$$
Let $\de^{+}_v$ denote the number of positive points on $\Si_v$. Each positive point is either a nodal point or a marked point, and $\de_v^{+}\geq \de_{v}$, as every contact point is a positive point.

\bLm{Simple_lmm}
If $D$ is Nef, then $\de_v^{+}\!\leq\! 2$ for all $v\!\in\!\V_{\tn{b}}$.
\eLm

\bPf
Assume  more than two points in $z_v\cup q_{v}$ are positive. Since $\Si'$ has at most  $2$ special points, removing $\Si_v$ from $\Si'$ we get some sub-clusters, at least one of which, say $\Si''$  has the following properties:
\bIt
\item it is not connected to the non-contracted (principal) part of $\Si$,
\item it does not carry any of the marked points,
\item it is connected to $\Si_v$ at a node $q_{\uvec{e}}\!\in\!\Si_{v}$ which is a positive point. 
\eIt
Let $\Si_{v_1}$ be the component of $\Si''$ connected to $\Si_v$ at 
$\Si_{v_1}\!\ni\!q_{\scz{\ucev{e}}}\sim q_{\uvec{e}}\!\in\!\Si_{v}$. Since, by assumption, $\tn{ord}^{i}_{u_v,\ze_v}(q_{\uvec{e}})>0$, for some $i\!\in\![N]$, the image of $u_{v_1}$ most lie in $D_{i}$ and the meromorphic section $\ze_{v_1,i}$ on $\Si_{v_1}$ corresponding to that should have a non-trivial pole (of the same order as $s_{\uvec{e},i}$)  at $q_{\scz{\ucev{e}}}$. By the Nef assumption, the line bundle 
$$
u_{v_1}^*\cN_XD_i \lra \Si_{v_1} \cong \P^1
$$
has a non-negative degree. Therefore, $\ze_{v_1,i}$ should have a non-trivial zero at another nodal point $q_{\uvec{e}_1}\!\in\! \Si_{v_1}$ and there is another component $\Si_{v_2}$ of $\Si''$ connected to $\Si_{v_1}$ at $\Si_{v_2}\!\ni\!q_{\scz{\ucev{e}}_1}\sim q_{\uvec{e}_1}\!\in\!\Si_{v_1}$. Continuing inductively we will see an infinite chain of irreducible components in $\Si''$, which is a contradiction.\ePf

\noindent
As a conclusion of this lemma, for every bubble component $f_v$ in $f$, since non-contact points (with $\partial D_{I_v}$) in $z_v\cup q_{v}$ are not relevant to the argument leading to (\ref{SP-cond_e}), by ignoring these points, replacing $n$ with $n-|I_v|$, and assuming $\ell=\ov\de_v\leq k=\de_v \leq 2$ in (\ref{GF_e}),  the semi-positivity condition (\ref{SP1_e}) in Definition~\ref{NCsemipos_dfn} guaranties that
\bEq{IfThen_e}
c_1^{TX(-\log D)}(B_v)+(n-I_v)-3+\ell\geq 0\; \Rightarrow\; c_1^{TX(-\log D)}(B_v)\geq 0.
\eEq

\noindent
The following statement summarizes the main result of this section.
\noindent
\bCr{non-simple_cr}
Suppose $(X,D,\om)$ is semi-positive in the sense of Definition~\ref{NCsemipos_dfn}, and $f$ is a log $(J,\nu)$-curve in $\ov\cM_{g,\mfs}(X,D,A,\nu)$. If $J$ belongs to the Baire sets $\tn{AK}^{I_v}(X,D)_{\cR}$ associated with  $\cM^\star_{0,\mft_v}(X,D,B_v)_{I_v}$ in Propositions~\ref{Trans2_prp}, for all $v\!\in\!\V_{\tn{b}}$, then $c_1^{TX(-\log D)}(B_v)\geq 0$ for all $v\!\in\!\V_{\tn{b}}$.
\eCr

\bPf
By the assumption on $J$, each moduli space $\cM^\star_{0,\mft_v}(X,D,B_v)_{I_v}$ is a non-empty smooth manifold of the expected complex dimension 
$$
c_1^{TX(-\log D)}(B_v)+ \!(n-I_v)-3\!+\!\ov\de_v\geq 0.
$$ 
The conclusion follows from (\ref{IfThen_e}).
\ePf

\noindent
Lemma~\ref{Simple_lmm} puts a major restriction on a contracted cluster $\Si'$. In general, each cluster will be of one of the following types. \\

\noindent
(i) A cluster with one node $q_{\uvec{e}_0}\!\in\!\Si_{v_0}$ that connects $\Si'$ to the principal part and no marked points. In this case $\Si'$ is a rooted tree with the root $v_0$. Distance from $v_0$ defines a partial order  $\prec$ on the vertices of $\Si'$ with $v_0$ being the minimal vertex. By the same inductive reasoning as in the proof of Lemma~\ref{Simple_lmm}, if $v\prec v'$ and $q_{\uvec{e}}\!\in\!\uvec{\E}_{v,v'}$ then $q_{\uvec{e}}$ can not be a positive point. Therefore each $\Si_{v'}$ in the cluster has at most one positive point, that will be the unique nodal point $q_{\scz\ucev{e}}$ connecting $\Si_{v'}$ to the unique component $\Si_{v}$ with  $v\prec v'$. Furthermore, $I_v\supseteq I_{v'}$ whenever $v\prec v'$, and $I_v\nsupseteq I_{v'}$ only if $q_{\scz\ucev{e}}$ is a contact point. Figure~\ref{TypeiCluster_fig}.(Left) illustrates the situation with the $``+"$ sign indicating a positive point. At a positive nodal point $q_{\uvec{e}}$, we have $0\!\neq\!s_{\uvec{e}}\in \N^N$. At nodal points without a $``+"$ sign on any side, $s_{\uvec{e}}$ must be zero. We put a $0$ to indicate those points. 
\begin{figure}
\begin{pspicture}(2,-1)(11,2)
\psset{unit=.3cm}

\rput(23,9){\small{$\Si'$}}
\pscircle(20,2){3}\rput(20,2){\small{$\Si_{v_0}$}}
\pscircle(19.13,6.92){2} 
\pscircle(23.53,5.53){2}
\pscircle(24.92,1.13){2}
\pscircle(27.70,3.9){2}
\pscircle(28.4,-.8){2}

\pscircle*(18.3,0.3){.25}\rput(19.2,0.8){\small{$q_{\uvec{e}_0}$}}

\pscircle*(19.48,4.954){.25}\rput(19.5,5.5){\tiny{+}}
\pscircle*(22.12,4.12){.25}\rput(22.55,4.4){\tiny{0}}
\pscircle*(22.954,1.48){.25}\rput(23.5,1.5){\tiny{0}}
\pscircle*(26.3,2.5){.25}\rput(26.7,2.85){\tiny{+}}
\pscircle*(26.7,0.3){.25}\rput(27.1,0){\tiny{+}}

\pscircle(35,4){2}\pscircle*(33.7,4){.25}\rput(34.5,5){\small{$q_{\uvec{e}_0}$}}
\rput(35,3){\small{$\Si_{v_{01}}$}}
\pscircle(39,4){2}
\rput(39,4){\small{$\Si_{v_{02}}$}}
\pscircle(43,4){2}
\rput(43,4){\small{$\Si_{v_{03}}$}}
\pscircle*(43,2){.25}\pscircle[linestyle=dashed](43,-1){3}
\rput(43,-1){\tiny{cluster $\Si_3'$}}
\rput(43,-2){\tiny{of type (i)}}
\pscircle(47,4){2}
\rput(47,3){\small{$\Si_{v_{04}}$}}
\pscircle*(48.3,4){.25}\rput(47.9,4.8){\small{$z$}}
\rput(42,7){\small{$\Si'$}}
\pscircle*(37,4){.25}
\pscircle*(41,4){.25}
\pscircle*(45,4){.25}

\end{pspicture}
\caption{(Left): A cluster of type (i). (Right): A cluster of type (ii).}
\label{TypeiCluster_fig}
\end{figure}
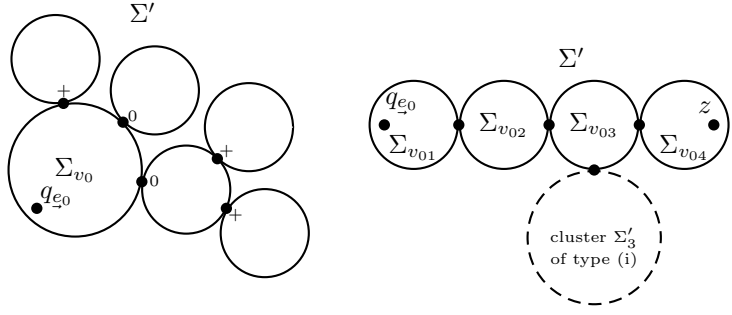
If each $D_i$ is positive, in the sense that $A\cdot D_i\!>\! 0$ for all $A\!\in\! \pi_2(X)$ such that $\om(A)\!>\!0$, then each bubble component   will have exactly one positive point.\\

\noindent
(ii) A cluster with one node $q_{\uvec{e}_0}\!\in\!\Si_{v_{01}}$ that connects $\Si'$ to the principal part and one marked point $z$ on $\Si_{v_{0k}}$ (possibly $v_{01}\!=\!v_{0k}$). In this case, $\Si'$ includes a chain of components $\Si_{v_{01}},\ldots,\Si_{v_{0k}}$.  For each $\Si_{v_{0i}}$, $\Si'$ might include a cluster $\Si_i$ (or more)  of type (i) that is attached to $\Si_{v_{0i}}$ at a unique nodal point $q_{\uvec{e}_{i0}}\in \Si_i$; see Figure~\ref{TypeiCluster_fig}.(Middle). \\

\noindent
(iii) A cluster with two nodes $q_{\uvec{e}_0}\!\in\!\Si_{v_{01}}$ and $q_{\uvec{e}_k}\!\in\!\Si_{v_{0k}}$ that connect $\Si'$ to the principal part. With $q_{\uvec{e}_k}$ in place of $z$, this case is like case (ii) above. 

\bRm{de2_rmk}
If $\de_v\!\geq\! 2$ for some geometric reason, then the only possibility for $\Si'$ is a chain of bubble components as in (ii) or (iii) between the two special points.
This is for example the case if $X$ is a toric variety and $D$ is the boundary divisor.
\eRm

%--------------------------------------------------------
\subsection{Non-simple maps}\label{MCMap_ss}

\noindent
In \cite{RT}, in order to prove the classical analogue of Proposition~\ref{MC_pr} (i.e., \cite[Thm~3.11]{RT}), in a process which we will call it \textbf{RT-process} here, they replace a non-simple map $f$ with an underlying simple map $f'$ with \textbf{multi-nodes}. A multi-node $m$ is a point at which more than two components of the domain are connected to each other. In a nodal domain, a node $q_{e}$ is obtained by attaching two irreducible components $\Si_{v_1}$ and $\Si_{v_2}$ at nodal points $q_{\uvec{e}_1}\!\in\!\Si_{v_1}$ and $q_{\uvec{e}_2}\!\in\!\Si_{v_1}$. A multi-node $q_m$ is obtained by attaching more than two components $\Si_{v_1},\ldots,\Si_{v_{\ell}}$ at nodal points $q_{\uvec{e}_i}\!\in\!\Si_{v_i}$, with $i\!=\!1,\ldots, \ell$. To keep the notation inline with the rest of the paper, we let $\E$ to denote the set of nodes and multi-nodes, and $\uvec{\E}$ to denote the set of nodal points on different components. When there is no multi-node we have $|\uvec{\E}|\!=\!2|\E|$; otherwise, $|\uvec{\E}|\!>\!2|\E|$. In the presence of an SNC divisor $D=\bigcup_{i\in [N]}D_i$, the sets $\E$ and $\uvec{\E}$ admit decompositions 
$$
\E=\bigcup_{I\subset [N]} \E_I\qquad \tn{and}\qquad \uvec{\E}=\bigcup_{I\subset [N]} \uvec{\E}_I
$$ 
by the type of these points. Similarly to Section~\ref{NodalMap_ss}, we  write $\V=\V_{\tn{p}}\cup \V_{\tn{b}}$, where $\V_{\tn{p}}$ corresponds to set of non-contracted or \textbf{principal} components and $ \V_{\tn{b}}$ corresponds to set of contracted or \textbf{bubble} components. Furthermore, we  write $\V_{\tn{b}}=\V_{\tn{b},\bullet}\cup \V_{\tn{b},\circ}$, where $A_v\!\neq\!0$ for  $v\!\in\!\V_{\tn{b},\bullet}$ and $ \V_{\tn{b},\circ}$ is the set of ghost bubbles.\\

\noindent
Ruan-Tian proved that the dimension of the space of such $f'$ is at least $2$ real dimension less than the dimension of the main stratum. The same argument does not directly work for log maps. If we apply the same process to the underlying $(J,\nu)$-map $\ov{f}$ of a log map $f$, the resulting simple map $\ov{f}'$ may not lift to a log map $f'$ for at least two reasons: (1) since we replace a multiple-cover map with its underlying simple map, the matching condition in Definition~\ref{PreLogMap_dfn}.\ref{MatchingOrders_l} may no longer hold, 
(2) since we identify different components with the same image, the vectors $s_v$ satisfying Definition~\ref{LogMap_dfn}.\ref{Tropical_l} may no longer exist. \\

\noindent
Suppose 
$$
f=\bigg(\phi,\big(u_v,\ze_v=(\ze_{v,i})_{i\in I_v},C_v=(\Si_v,\mfj_v,\vec{z}_v\cup \vec{q}_v)\big)_{v\in \V}, J, \nu\bigg)\in \cM^{\tn{ns}}_{g,\mfs}(X,D, A,\nu)_\Gamma.
$$
The RT-process, described after \cite[Dfn 3.10]{RT}, changes the underlying $(J,\nu)$-map
$$
\ov{f}=\bigg(\phi,\big(u_v,C_v=(\Si_v,\mfj_v,\vec{z}_v\cup \vec{q}_v)\big)_{v\in \V}, J, \nu\bigg)\in \cM^{\tn{ns}}_{g,k}(X, A,\nu)_\Gamma.
$$
to another $(J,\nu)$-map  
$$
\ov{f}''=\bigg(\phi,\big(u_v,C_v=(\Si_v,\mfj_v,\vec{z}_v\cup \vec{q}_v)\big)_{v\in \V_{\tn{p}}}, \big(u_{v''},C_{v''}=(\Si_{v''},\mfj_{v''},\vec{z}_{v''}\cup \vec{q}_{v''})\big)_{v''\in \V''_{\tn{b}}}, J, \nu\bigg).
$$ 
with multi-nodes in the following way.\\

\noindent
 (i) It collapses the ghost bubbles (and any marked point on it will be thrown away). As a result we get some multi-nodes. \\
(ii) It replaces each multiple-cover bubble component by its image. Since some of the special (marked or nodal) points may have the same image, this step may produce further multi-nodes. \\
(iii) It collapses each sub-tree of the bubbles whose components have the same image. \\

\noindent
None of these three steps changes the genus though. Let $\Gamma'$ be the resulting combinatorial type of the domain with components indexed by $\V'$, nodes and multi-nodes indexed by $\E'$, and nodal points indexed by $\uvec{\E}'$. With notation as in Section~\ref{NodalMap_ss}, $\V'$ decomposes as $\V'_{\tn{p}}\!\cup\!\V'_{\tn{b}}$. The first component is identical to $\V_{\tn{p}}$. There is a collapsing map 
\bEq{red_e}
\mf{red}\colon \V_{\tn{b},\bullet}\lra \V'_{\tn{b}}
\eEq 
and a multiplicity map 
\bEq{degree_e}
d\colon \V'_{\tn{b}}\lra \Z_+
\eEq
such that 
$$
\sum_{v\in \mf{red}^{-1}(v')} A_v =  d_{v'} A_{v'} \qquad \forall~v'\!\in\!\V'_{\tn{b}}.
$$
Here $A_{v'}$ is the homology class of the resulting simple curve. 
After these three steps, there might still be bubble components (not adjacent to each other) which have the same image or their nodal points have the same image.\\

\noindent
(iv) We identify components with the same image and nodal points with the same image.\\

\noindent
After step (iv), we a get a domain with possibly further multi-nodes and higher genus. Let $\Gamma''$ be the resulting combinatorial type of this domain with components indexed by $\V''$, nodes and multi-nodes indexed by $\E''$, and nodal points indexed by $\uvec{\E}''$. The maps (\ref{red_e}) and (\ref{degree_e}) descend to
$$
\V_{\tn{b},\bullet}\stackrel{\mf{red}}{\lra} \V_{\tn{b}}' \stackrel{\mf{red}}{\lra} \V_{\tn{b}}'', \quad d\colon \V''_{\tn{b}}\lra \Z_+,
$$
with the same properties as above. For the same reason as in \cite[Cr.~3.17]{RT}, we have
\bEq{Evec-E_e}
|\uvec{\E}''_I|-|\E''_I|=|\uvec{\E}'_I|-|\E'_I|.
\eEq
Each nodal point of $f''$ is still decorated by a well-defined subset $I\!\subset\![N]$ such that all of its pre-image nodal points in $f$ have the same decoration $I$. Also, if a bubble component $\Si_{v''}$ of  $f''$ has image in $D_{I_{v''}}$, then all of its pre-image bubble components have image in the same stratum, i.e.
$$
I_{v}=I_{v''}\qquad \forall~v\!\in\!\mf{red}^{-1}(v'').
$$
The combinatorial type $\gamma$ of this process is encoded in the triple $(\Gamma,\Gamma',\Gamma'')$, and the associated maps $\mf{red}$ and $d$. By Theorem~\ref{Compactness_th}, the set of such $\gamma$ is finite. Let $\cM_{g,k}(X,A,\nu)_{\Gamma''}$ denote the classical moduli space of such $(J,\nu)$ curves $\ov{f}''$. 
By \cite[Prp~3.21]{RT}, if $(X,\om)$ is semi-positive, for generic $(J,\nu)$, $\cM_{g,k}(X,A,\nu)_{\Gamma''}$ is a smooth moduli space of the $\C$-dimension at most 
\bEq{UL_e}
c_1^{TX}(A)+(n-3)(1-g)+k-(|\uvec{\E}''|-|\E''|).
\eEq
Comparing (\ref{UL_e}) with \cite[Prp~3.21]{RT}, note that 
$$
|\uvec{\E}''|-|\E''|=|\uvec{\E}'|-|\E'|=n_{\phi}+|\V'_{\tn{b}}|,
$$
where $n_{\phi}$ is the number of nodes of the stable domain $\phi(\Si)\!\in\!\ov{\mf{U}}_{g,k}$.\\

\noindent
Away from the principal components, the map $\ov{f}''$ does not lift to a log curve. However, some of the information still passes to $\ov{f}''$. First, let us consider the pre-log space 
$$
\cM^{\tn{plog}}_{g,\mfs}(X,D,A,\nu)_\Gamma,
$$
that is we forget about the Condition~\ref{GObs_e} in Definition~\ref{LogMap_dfn}. Instead of $\cM_{g,k}(X,A,\nu)_{\Gamma''}$, we consider the set $\cM^{\tn{plog}}_{g,\mfs}(X,D,A,\nu)_{\Gamma''}$ of tuples 
$$
f''=\bigg(\phi,\big(u_v,\ze_v,C_v=(\Si_v,\mfj_v,\vec{z}_v\cup \vec{q}_v)\big)_{v\in \V_{\tn{p}}}, \big(u_{v''},C_{v''}=(\Si_{v''},\mfj_{v''},\vec{z}_{v''}\cup \vec{q}_{v''})\big)_{v''\in \V''_{\tn{b}}}, J, \nu\bigg)
$$ 
obtained from the elements of $\cM^{\tn{plog}}_{g,\mfs}(X,D,A,\nu)_\Gamma$, where the principal components still carry the information of the meromorphic sections $\ze_v$. 
Let $\cM_{g,\mfs}^{\tn{plog},\gamma}(X,D,A,\nu)_{\Gamma}$ denote its pre-image in $\cM^{\tn{plog}}_{g,\mfs}(X,D,A,\nu)_{\Gamma}$, i.e. those pre-log maps for which the RT-process is of type $\gamma$. The projection
\bEq{pigamma_e}
\pi_\gamma\colon \cM_{g,\mfs}^{\tn{plog},\gamma}(X,D,A,\nu)_{\Gamma}\lra \cM^{\tn{plog}}_{g,\mfs}(X,D,A,\nu)_{\Gamma''}
\eEq
is a surjective fiber bundle. The key point is that by (\ref{g0MI_e}), if $\Si_v$ is genus $0$, as long as the  second combinatorial condition in (\ref{CombCond_e}) is satisfied, for each $i\!\in\! I_v$, there are meromorphic sections $\ze_{v,i}$ of $u_v^*\cN_XD_i$ with zeros/poles of orders $s_{\uvec{e},i}$ and $s_{a}$ at $q_{\uvec{e}}$ and $z_a$, respectively, for all $\uvec{e}\!\in\!\uvec{\E}_v$ and $z_a\!\in\!z_v$. The fiber $\cM_{f''}$ of (\ref{pigamma_e}) over any $f''$ is a product of the manifolds of the form described below and 
$$
\tn{st}\times\tn{ev}\colon   \cM^{\tn{plog},\gamma}_{g,\mfs}(X,D,A,\nu)_\Gamma \lra \ov\cM_{g,k}\!\times\!X^{\mfs}
$$
factors through $\pi_\gamma$. \\

\noindent 
(1) For each $v\!\in\!\V_{\tn{b},\circ}$, we have the configuration space $\cM_{0,k_v+\ell_v}$ of the special points on the ghost bubble $\Si_v$ in $\pi^{-1}(f'')$. \\

\noindent
(2) For all $v''\!\in\!\V''$ and $v\!\in\!\mf{red}^{-1}(v'')$, we have 
$u_v=u_{v''}\circ h_{v}$ for some degree $d_v$ covering map $h_{v}\colon \Si_v\lra \Si_{v''}$ as in (\ref{h_e}). Note that 
$$
d_{v''}=\sum_{v\in \mf{red}^{-1}(v'')} d_{v}.
$$
The combinatorial type $\al_v$ of $h_v$ is determined by the image $z_{v''}\cup q_{v''}$ of $z_{v}\cup q_{v}$ in $\Si_{v''}$ and the branching order of $h_v$ at $z_{v}\cup q_{v}$ (i.e. the partition of $d_v$ as in (\ref{Partition_e}) at the contact points 
$$
u_{v''}^{-1}(\partial D_{I_v''})\subset z_{v''}\cup q_{v''}
$$
and the branching orders at the rest of the points). Therefore, as in (\ref{alTostar_e}), in the fiber over $f''$ we get the relative moduli space 
\bEq{Fiberv_e}
\cM_{0,\al_v}(\Si_{v''}=\P^1,z_{v''}\cup q_{v''},[d_v])
\eEq
of tuples 
$$
\big(h_v\colon \Si_{v}=\P^1\lra \Si_{v''}=\P^1, \vec{z}_{v}\cup \vec{q}_{v}\big)
$$ 
relative to the divisor $z_{v''}\cup q_{v''}\subset \Si_{v''}$, with the ramification/tangency order data $\al_{v}$. The diffeomorphism type of (\ref{Fiberv_e}) is independent of the location of $z_{v''}\cup q_{v''}$ (and thus $f''$).    \\

\noindent
We conclude that 
\bEq{DimMfdp_e}
\cM_{f''}=\prod_{v\in \V_{\tn{b},\circ}} \cM_{0,k_v+\ell_v} \times \prod_{v''\in \V''} \prod_{v\in \tn{red}^{-1}(v'')} \cM_{0,\al_v}(\P^1,z_{v''}\cup q_{v''},[d_v]).
\eEq

\noindent 
Similarly to \cite[Thm~3.16]{RT}, and with a similar proof as in Section~\ref{NodalMap_ss}, for $(\om',\cR,J,\nu)$ in a subset of second category  $\cH^\gamma_{g,k}(X,D)_{[\om]}\subset \cH_{g,k}(X,D)_{[\om]}$, $\cM^{\tn{plog}}_{g,\mfs}(X,D,A,\nu)_{\Gamma''}$ is a smooth manifold of $\C$-dimension 
$$
\aligned
k''+|\uvec{\E}''|+&\sum_{v\in \V_{\tn{p}}} \big(c_1^{TX(-\log D)}(A_v)+ (n-3)(1-g_v)-|I_v|\big)\\
+&\sum_{v\in \V''_{\tn{b}}} \big(c_1^{TX(-\log D)}(A_{v''})+ (n-3)-|I_{v''}|\big)\\
-&\sum_{I\subset [N]} (n-|I|)(|\uvec{\E}''_I|-|\E''_I|).
\endaligned
$$
Note the number $h_{\ov{D}}$ in \cite[p.~485-486]{RT} is $|\uvec{\E}''|$ in our notation and  the number $t_{\ov{D}}$ there is $|\E''|$. Also $k''$ denotes the number of surviving\footnote{which will be $k$ or $k-1$ since the contracted part carries at most one of the marked points.} marked points.
For generic $(\om',\cR,J)$, 
\bEq{PosD_e}
c_1^{TX(-\log D)}(A_{v''})+ (n-3)-|I_{v''}|+\de_{v''}\geq 0,
\eEq
where $\de_{v''}$ is the number of contact points of $u_{v''}$ with $\partial D_{I_{v''}}$ as in Corollary~\ref{non-simple_cr}.
If some bubble in $\Gamma''$ happens to be the image of two or more bubbles in $\Gamma'$, by adding (\ref{PosD_e}) to the dimension, and by (\ref{Evec-E_e}), we get 
$$
\aligned
\tn{dim}_\C\; \cM^{\tn{plog}}_{g,\mfs}(X,D,A,\nu)_{\Gamma''} \leq k'+|\uvec{\E}'|+&\sum_{v\in \V_{\tn{p}}} \big(c_1^{TX(-\log D)}(A_v)+ (n-3)(1-g_v)-|I_v|\big)\\
+&\sum_{v\in \V'_{\tn{b}}} \big(c_1^{TX(-\log D)}(A_{v'})+ (n-3)-  |I_{v'}|\big) \\
-&(n-1)(|\uvec{\E}'|-|\E'|)+\sum_{I\subset [N]} (|I|-1)(|\uvec{\E}''_I|-|\E''_I|).
\endaligned
$$
Since 
$$
1-g=\sum_{v\in \V'} (1-g_v) -( |\uvec{\E}'|-|\E'|),
$$
we get
\bEq{GenDimFormula_e}
\aligned
\tn{dim}_\C\; \cM^{\tn{plog}}_{g,\mfs}(X,D,A,\nu)_{\Gamma''} &\leq  k'+ (n-3)(1-g)+\sum_{v'\in \V'} \big(c_1^{TX(-\log D)}(A_{v'})- |I_{v'}|\big)\\\
+&(2|\E'|-|\uvec{\E}'|)+\sum_{I\subset [N]} (|I|-1) (|\uvec{\E}'_I|-|\E'_I|).
\endaligned
\eEq
By the semi-positivity condition, we have $c_1^{TX(-\log D)}(A_{v'})\geq 0$ for all $v'\in \V_{\tn{b}}'$, see Corollary~\ref{non-simple_cr}. 
Also $|\uvec{\E'}|\geq 2|\E'|$. Therefore, the last equation is less than or equal to 
\bEq{SimEst_e}
\aligned
&c_1^{TX(-\log D)}(A)+ (n-3)(1-g)+k\\
&- \sum_{v'\in \V'} |I_{v'}|+\sum_{I\subset [N]} (|I|-1)(|\uvec{\E}'_I|-|\E'_I|).
\endaligned
\eEq

\noindent
Also, note that if the equality happens, then
$$
\sum_{v''\in \V''} |I_{v''}|=\sum_{v'\in \V'} |I_{v'}|\qquad \tn{and}\qquad|\uvec{\E'}|= 2|\E'|.
$$
The first equality implies that step (iv) in RT-process is trivial. The second one implies that step (i) is trivial and no multi-node is created in steps (ii) and (iii) of the process.
Then it is easy to see that
$$
\sum_{v''\in \V''} |I_{v''}|+\sum_{I\subset [N]} (|I|-1)(|\uvec{\E}''_I|-|\E''_I|)= -|\E|-\sum_{v\in \V} |I_{v}|+\sum_{e\in \E} |I_e|= \tn{dim}_\C(\mc{G})-\dim_\R \K_\R(\Gamma).
$$

\newtheorem*{PFMC_th}{Proof of Proposition~\ref{MC_pr}}
\begin{PFMC_th}
If $D$ is smooth, i.e. $N\!=\!1$, then the second line (\ref{SimEst_e}) is negative. We conclude that
$$
\tn{dim}_\C\; \cM^{\tn{plog}}_{g,\mfs}(X,D,A,\nu)_{\Gamma''}<\tn{dim}_\C\;\cM_{g,\mfs}(X,D,A,\nu).
$$

\noindent
In this sense, for $D$ smooth, Proposition~\ref{MC_pr} 
essentially follows from the classical result of Ruan-Tian by looking at the image of non-simple maps in $\ov\cM_{g,k}(X,A,\nu)$.\qed
\end{PFMC_th}

\noindent
For an arbitrary SNC divisor $D$, in order to take care of the extra term 
\bEq{Deficit_e}
(2|\E'|-|\uvec{\E}'|)-\sum_{v'\in \V'} |I_{v'}|+\sum_{I\subset [N]} (|I|-1)(|\uvec{\E}'_I|-|\E'_I|)
\eEq
in (\ref{GenDimFormula_e}), we need to use  Condition~\ref{GObs_e} in Definition~\ref{LogMap_dfn} to reduce the dimension. We encounter the following two problems.
\bEn
\item\label{Problem1} The map $\tn{ob}_\Gamma$ into $\mc{G}(\Gamma)$ is defined on $\cM^{\tn{plog},\gamma}_{g,\mfs}(X,D,A,\nu)_\Gamma$. As the examples below show, unlike $\tn{st}\times\tn{ev}$, $\tn{ob}_\Gamma$ does not necessarily factor through $\pi_\gamma$. Therefore, we need to work with the larger space $\cM^{\tn{plog},\gamma}_{g,\mfs}(X,D,A,\nu)_\Gamma$.
\item\label{Problem2}  Since the elements of $\cM^{\tn{plog},\gamma}_{g,\mfs}(X,D,A,\nu)$ are not simple, transversality of $\tn{ob}_\Gamma$ in the sense of Claim 2 of proof of Theorem~\ref{TransGamma_thm} might not be achievable. We need to replace $\mc{G}(\Gamma)$ with a smaller group $\mc{G}(\gamma)$ that admits a surjective homomorphism $h\colon \mc{G}(\Gamma)\lra\mc{G}(\gamma)$, and such that
$$
\tn{ob}_\gamma=h\circ \tn{ob}_{\Gamma}\colon \cM^{\tn{plog},\gamma}_{g,\mfs}(X,D,A,\nu)\lra\mc{G}(\gamma)
$$ 
can be transversed (i.e. $1\!\in\!\mc{G}(\gamma)$ is a regular value of that).
\eEn

\noindent
The two problems above are not specific to the particular compactification considered in this work and should appear, either explicitly or implicitly, in any other analytical approach.\\

\noindent
The first example below illustrates a simple situation where $\tn{ob}_\Gamma$ depends on the location of special points on a ghost bubble. In Appendix~\ref{GH_ss}, we study this dependence in details. The second example below illustrates a situation where $\tn{ob}_\Gamma$ depends on the choice of the covering map $h_v$ for some multiple cover bubble $\Si_v$. The third example below illustrates a situation where different bubbles have the same image, and $\tn{ob}_\Gamma$ can not be transversed. Finally, the fourth example below illustrates the necessity of the extra condition (\ref{SP2_eN}), whenever $N>1$. 

\begin{example}\label{GhostDep_ex}
Let $X\!=\!\P^2$, $D\!=\!D_1\!\cup\! D_2$ be the transverse union of two coordinate lines, 
$$
g\!=\!0,\quad k\!=\!2,\quad A\!=\![3],\quad\tn{and}\quad \mfs\!=\!(s_1,s_2)\!=\!((2,1),(1,2)).
$$ 
Let $\Gamma$ be the configuration in Figure~\ref{3line2pt_fig}; it is the nodal configuration obtained by three lines  (not contained in $D$) passing through the point $D_{12}$, together with a ghost bubble mapped to $D_{12}$ that connects the domains of these three lines and carries the marked points $z_1,z_2$.
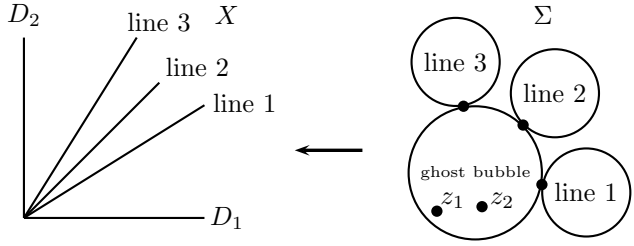
\begin{figure}
\begin{pspicture}(-3,-.5)(11,2.5)
\psset{unit=.3cm}
\psline(0,0)(8,0)
\psline(0,0)(0,8)
\psline(0,0)(8,5)\rput(9.8,5.1){\small{line $1$}}
\psline(0,0)(6,6)\rput(7.7,6.7){\small{line $2$}}
\psline(0,0)(5,8)\rput(6,8.7){\small{line $3$}}

\rput(9,9){\small{$X$}}
\rput(9,0){\small{$D_1$}}\rput(0,9){\small{$D_2$}}

\rput(23,9){\small{$\Si$}}
\pscircle(20,2){3}\rput(20,2){\tiny{ghost bubble}}
\pscircle(19.13,6.92){2} \rput(19.13,6.92){\small{line $3$}}
\pscircle(24.92,1.13){2}\rput(24.92,1.13){\small{line $1$}}
\pscircle(23.53,5.53){2}\rput(23.53,5.53){\small{line $2$}}
\pscircle*(18.3,0.3){.25}\rput(19,0.8){\small{$z_1$}}
\pscircle*(20.3,0.5){.25}\rput(21.2,0.8){\small{$z_2$}}
\pscircle*(19.48,4.954){.25}
\pscircle*(22.12,4.12){.25}
\pscircle*(22.954,1.48){.25}
\psline[linewidth=.15]{->}(15,3)(12,3)

\end{pspicture}
\caption{A nodal configuration in the boundary of $\ov\cM^{\log}_{0,((3,3))}(\P^2,D,[3])$.}
\label{3line2pt_fig}
\end{figure}
Let $\V=\{v_0,v_1,v_2,v_3\}$ where $v_1,v_2,v_3$ correspond to the lines $1$, $2$, and $3$, respectively, and $v_0$ corresponds to the ghost bubble.  Let $\E\!=\!\{e_1,e_2,e_3\}$ where $e_i$ corresponds to the node connecting the domain of the $i$-th line to the ghost bubble. We can choose the orientation $\uvec{e}_i$ to be the one ending at $v_0$. We have 
$$
I_{v_0}=\{1,2\}, \quad 
I_{v_i}=\eset, \quad 
I_{e_i}=\{1,2\}, \quad
s_{\uvec{e}_i}=(1,1)\qquad \forall~i=1,2,3.
$$
The map 
$$
\vr\colon \Z^{\E} \oplus \Z^{I_{v_0}}\lra \bigoplus_{i=1}^3 \Z^{I_{e_i}}
$$
has a $1$-dimensional kernel and a $2$-dimensional cokernel. The kernel is generated by 
$$
(\la_{e_1},\la_{e_2},\la_{e_3},s_{v_0})=(1,1,1,(1,1)).
$$
Therefore, Condition~\ref{Tropical_l} of Definition~\ref{LogMap_dfn} is satisfied.
The obstruction group $\mc{G}(\Gamma)$ is $2$-dimensional. The homomorphism 
$$
 (\C^*)^{I_{e_1}} \times (\C^*)^{I_{e_2}} \times (\C^*)^{I_{e_3}}\lra (\C^*)^2,\quad 
 ((x_1,y_1),(x_2,y_2),(x_3,y_3))\lra \Big(\frac{x_1}{y_1}\frac{y_2}{x_2}, \frac{x_2}{y_2}\frac{y_3}{x_3}\Big)
$$
descends to an isomorphism $\mc{G}(\Gamma)\lra (\C^*)^2$.
We take the marked point $z_1$ to be $0$, $z_2$ to be $\infty$, and the nodal points $q_{\scz\ucev{e}_1}$, $q_{\scz\ucev{e}_2}$, and $q_{\scz\ucev{e}_3}$ to be $1$, $\al_2$, and $\al_3$, respectively. Thus, $\al_2$ and $\al_3$ are parametrizing the two-dimesnional configuration space $\cM_{0,5}$ of the five special points on $\Si_{v_0}$. The meromorphic sections (functions) $\ze_{v_0,1}$ and  $\ze_{v_0,2}$ are given by 
$$
\ze_{v_0,i}(z)=\frac{\beta_i\,z^{s_{1,i}}}{(z-1)^{s_{\uvec{e}_1,i}}\,(z-\al_2)^{s_{\uvec{e}_2,i}}\,(z-\al_3)^{s_{\uvec{e}_3,i}}}=
\begin{cases}
\frac{\beta_1\,z^2}{(z-1)\,(z-\al_2)\,(z-\al_3)} &\qquad i=1,\\
\\
\frac{\beta_2\,z}{(z-1)\,(z-\al_2)\,(z-\al_3)} &\qquad i=2,
\end{cases}
$$
where $\beta_i\!\in\!\C^*$. We conclude that 
$$
\aligned
\tn{ob}_{\Gamma}(f)= \Big(\frac{m_1}{m_2\al_2},\frac{m_2\al_2}{m_3\al_3}\Big)\in (\C^*)^2,
\endaligned
$$
where $m_i$ is the slope of the $i$-th line. The RT process removes the ghost bubble $\Si_{v_0}$ and creates a map $f'$ with one multi-node. In general, as we will show in Section~\ref{GH_ss},  the element $\tn{ob}_{\Gamma}(f)$ can be expressed as a product of two terms 
$$
\tn{ob}_{\Gamma}(f)=\tn{ob}_{\Gamma}(f')\cdot o_{v_0}(\P^1,z_1,z_2,q_{\scz\ucev{e}_1},q_{\scz\ucev{e}_2},q_{\scz\ucev{e}_3}).
$$ 
Here, the decomposition is 
\bEq{fTofo_e}
\Big(\frac{m_1}{m_2\al_2},\frac{m_2\al_2}{m_3\al_3}\Big)=
\Big(\frac{m_1}{m_2},\frac{m_2}{m_3}\Big)\cdot 
\Big(\al_2^{-1},\frac{\al_2}{\al_3}\Big)
\eEq
However, as (\ref{fTofo_e}) shows, the second term on righthand side can be a non-trivial function from $\cM_{0,k_{v_0}+\ell_{v_0}}$ into $\mc{G}(\Gamma)$.\qed
\end{example}

\noindent
Next, we consider a similar configuration in one dimension higher with a multiple-cover map in place of the ghost bubble above. 

\begin{example}\label{MCDep_ex}
Let $X\!=\!\P^3$, $D\!=\!D_1\!\cup\! D_2$ be the transverse union of two hyperplanes, 
$$
g\!=\!0,\quad k\!=\!2,\quad A\!=\![5],\quad\tn{and}\quad \mfs\!=\!(s_1)\!=\!((5,0),(0,5)).
$$ 
Let $\Gamma$ be the configuration in Figure~\ref{3linewMC_fig}; it is the nodal configuration obtained by three lines  (not contained in $D$) passing through $D_{12}\cong \P^1$, together with a degree $2$ multiple-cover map 
$$
h_{v_0}\colon \Si_{v_0}=\P^1 \lra D_{12}
$$ 
that connects the domains of these three lines and carries the marked points $z_1,z_2$. 
\begin{figure}
\begin{pspicture}(-4,-.5)(11,3)
\psset{unit=.3cm}
\psline(0,0)(8,0)(4,3)(-4,3)(0,0)
\psline(0,0)(0,8)(-4,11)(-4,3)(0,0)
\psline(-2,1.5)(6,6.5)\rput(7.8,6.6){\small{line $1$}}
\psline(-2,1.5)(4,7.5)\rput(5.7,8.2){\small{line $2$}}
\psline(-3,2.25)(2,10.25)\rput(3,10.95){\small{line $3$}}

\rput(4,2){\small{$D_1$}}\rput(-2,8){\small{$D_2$}}
\pscircle*(-2,1.5){.25}\rput(-2.75,1.5){\small{$\al$}}
\pscircle*(-3,2.25){.25}\rput(-3.75,2.05){\small{$\al_3$}}
\pscircle*(-1,0.75){.25}\rput(-1.75,.65){\small{$\beta$}}

\rput(23,9){\small{$\Si$}}
\pscircle(20,2){3}\rput(20,2){\tiny{double cover}}
\pscircle(19.13,6.92){2} \rput(19.13,6.92){\small{line $3$}}
\pscircle(24.92,1.13){2}\rput(24.92,1.13){\small{line $1$}}
\pscircle(23.53,5.53){2}\rput(23.53,5.53){\small{line $2$}}
\pscircle*(18.3,0.3){.25}\rput(19,0.8){\small{$z_1$}}
\pscircle*(20.3,0.5){.25}\rput(21.2,0.8){\small{$z_2$}}
\pscircle*(19.48,4.954){.25}
\pscircle*(22.12,4.12){.25}
\pscircle*(22.954,1.48){.25}
\psline[linewidth=.15]{->}(15,3)(12,3)

\end{pspicture}
\caption{A nodal configuration in the boundary of $\ov\cM^{\log}_{0,((5,0),(0,5))}(\P^3,D,[5])$.}
\label{3linewMC_fig}
\end{figure}
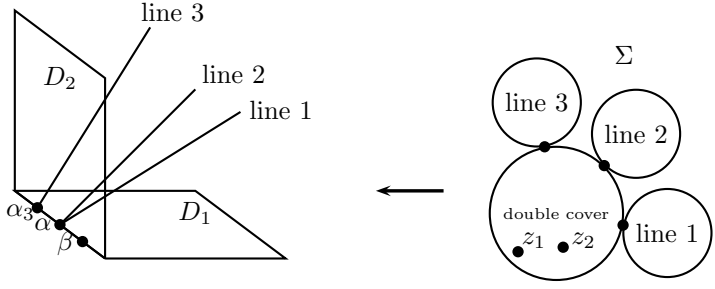
The dual graph $\Gamma$ and the isomorphism $\mc{G}\lra (\C^*)^2$ is the same as in Example~\ref{GhostDep_ex}. After a reparametrization, we may assume  that the double points are at $a,b\!\in \D_{12}=\P^1$ and $0,\infty \!\in\!\Si_{v_0}=\P^1$. Therefore, 
$$
h_{v_0}(z)=\frac{a-bz^2}{1-z^2}.
$$
The location of each special point on $\Si_{v_0}$ is determined up to a sign by its image in $D_{12}$, i.e. 
$$
h_{v_0}(q_{\scz\ucev{e}_i})=\al_i \leftrightarrow q_{\scz\ucev{e}_i}\in \pm\sqrt{\frac{a-\al_i}{b-\al_i}}, \quad h_{v_0}(z_i)=\beta_i \leftrightarrow z_i\in \pm\sqrt{\frac{a-\beta_i}{b-\beta_i}}.
$$
The RT-process replaces $u_{v_0}$ with the identification $u_{v'_0}\colon \Si_{v'_0}=\P^1\lra D_{12}$ together with the marked points $\{\al_1,\al_2,\al_3,\beta_1,\beta_2\}\!\in\!\Si_{v'_0}=\D_{12}$ (if two of these are the same, they will be identified).
In the fibers $\cM_{f'}$ of (\ref{pigamma_e}) over $f'$, the parameters $\al_1,\al_2,\al_3,\beta_1,\beta_2$ will be fixed, but $a,b$ are allowed to change and they parametrize the moduli space of multiple covers that yield the same  underlying simple marked curve $f'$. The type of $\Gamma'$ depends on the location of $\{\al_1,\al_2,\al_3,\beta_1,\beta_2\}$. For example, we may assume that
$$
\al_1=\al_2=\al \quad \tn{and} \quad \beta_1=\beta_2=\beta
$$
as in Figure~\ref{3linewMC_fig}. Then $\Gamma'$ has a multi-node at $\al$, just one marked point at $\beta$, and a regular node at $\al_3$.\\

\noindent
In any case, the meromorphic sections $\ze_{v_0,1}$ and  $\ze_{v_0,2}$ of $u_{v_0}^*\cN_XD_i=\cO(2)$ are given by 
$$
\ze_{v_0,i}(z)=\frac{c_i\,(z-z_1)^{s_{1,i}}(z-z_2)^{s_{2,i}}}{ z^{s_{\uvec{e}_1,i}}\,(z-1)^{s_{\uvec{e}_2,i}}\,(z-\al)^{s_{\uvec{e}_3,i}}}=
\begin{cases}
c_1 (z-z_1)^5\big((z-q_{\scz\ucev{e}_1})\,(z-q_{\scz\ucev{e}_2})\,(z-q_{\scz\ucev{e}_3})\big)^{-1} \quad &\tn{if}~~ i=1,\\
\\
c_2 (z-z_2)^5 \big((z-q_{\scz\ucev{e}_1})\,(z-q_{\scz\ucev{e}_2})\,(z-q_{\scz\ucev{e}_3})\big)^{-1} \quad &\tn{if}~~ i=2.
\end{cases}
$$
where $c_i\!\in\!\C^*$. We conclude that 
$$
\aligned
\tn{ob}_{\Gamma}(f)= \Bigg(\frac{m_1}{m_2} \bigg(\frac{(q_{\scz\ucev{e}_1}-z_1)(q_{\scz\ucev{e}_2}-z_2)}{(q_{\scz\ucev{e}_1}-z_2)(q_{\scz\ucev{e}_2}-z_1)}\bigg)^5,
\frac{m_2}{m_3}\bigg(\frac{(q_{\scz\ucev{e}_2}-z_1)(q_{\scz\ucev{e}_3}-z_2)}{(q_{\scz\ucev{e}_2}-z_2)(q_{\scz\ucev{e}_3}-z_1)}\bigg)^5\Bigg)\in (\C^*)^2,
\endaligned
$$
where $m_i$ is the slope of the $i$-th line with respect to $D_1$ and $D_2$.
The fractions 
$$
\frac{(q_{\scz\ucev{e}_1}-z_1)(q_{\scz\ucev{e}_2}-z_2)}{(q_{\scz\ucev{e}_1}-z_2)(q_{\scz\ucev{e}_2}-z_1)}\quad \tn{and} \quad  \frac{(q_{\scz\ucev{e}_2}-z_1)(q_{\scz\ucev{e}_3}-z_2)}{(q_{\scz\ucev{e}_2}-z_2)(q_{\scz\ucev{e}_3}-z_1)}
$$
are the cross-ratios of $(q_{\scz\ucev{e}_1},q_{\scz\ucev{e}_2},z_1,z_2)$ and $(q_{\scz\ucev{e}_2},q_{\scz\ucev{e}_3},z_1,z_2)$, respectively. Their dependence on $a$ and $b$, and thus on $h_{v_0}$,  is non-trivial.
\end{example}

\begin{example}\label{GoodEx1}
Let $X\!=\!\P^2$, $D\!=\!D_1\!\cup\! D_2$ be the transverse union of two coordinate lines, 
$$
g\!=\!0,\quad k\!=\!1,\quad A\!=\![3],\quad\tn{and}\quad \mfs\!=\!(s_1)\!=\!((3,3)).
$$ 
Let $\Gamma$ be the configuration in Figure~\ref{3line_fig}; it is the nodal configuration obtained by three lines  (not contained in $D$) passing through the point $D_{12}$, together with a ghost bubble mapped to $D_{12}$ that connects the domains of these three lines and carries the marked point $z_1$.
\begin{figure}
\begin{pspicture}(-3,-.5)(11,3)
\psset{unit=.3cm}
\psline(0,0)(8,0)
\psline(0,0)(0,8)
\psline(0,0)(8,5)\rput(9.8,5.1){\small{line $1$}}
\psline(0,0)(6,6)\rput(7.7,6.7){\small{line $2$}}
\psline(0,0)(5,8)\rput(6,8.7){\small{line $3$}}

\rput(9,9){\small{$X$}}
\rput(9,0){\small{$D_1$}}\rput(0,9){\small{$D_2$}}

\rput(23,9){\small{$\Si$}}
\pscircle(20,2){3}\rput(20,2){\tiny{ghost bubble}}
\pscircle(19.13,6.92){2} \rput(19.13,6.92){\small{line $3$}}
\pscircle(24.92,1.13){2}\rput(24.92,1.13){\small{line $1$}}
\pscircle(23.53,5.53){2}\rput(23.53,5.53){\small{line $2$}}
\pscircle*(18.3,0.3){.25}\rput(19,0.8){\small{$z_1$}}
\pscircle*(19.48,4.954){.25}
\pscircle*(22.12,4.12){.25}
\pscircle*(22.954,1.48){.25}
\psline[linewidth=.15]{->}(15,3)(12,3)

\end{pspicture}
\caption{A nodal configuration in the boundary of $\ov\cM^{\log}_{0,((3,3))}(\P^2,D,[3])$.}
\label{3line_fig}
\end{figure}
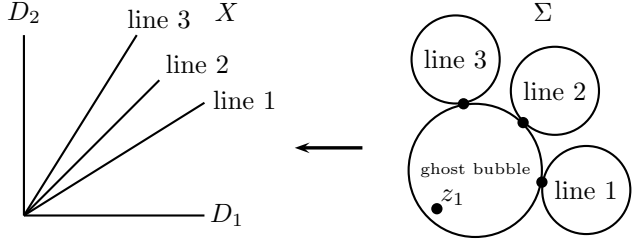
The dual graph $\Gamma$ is as in Example~\ref{GhostDep_ex} and thus the isomorphism $\mc{G}\lra (\C^*)^2$ is induced by 
$$
 (\C^*)^{I_{e_1}} \times (\C^*)^{I_{e_2}} \times (\C^*)^{I_{e_3}}\lra (\C^*)^2,\quad 
 ((x_1,y_1),(x_2,y_2),(x_3,y_3))\lra \Big(\frac{x_1}{y_1}\frac{y_2}{x_2}, \frac{x_2}{y_2}\frac{y_3}{x_3}\Big).
$$
We take the marked point $z_1$ to be $\infty$ and the nodal points $q_{\scz\ucev{e}_1}$, $q_{\scz\ucev{e}_2}$, and $q_{\scz\ucev{e}_3}$ to be $0$, $1$, and $\al$, respectively. Thus, $\al$ is parametrizing the configuration space $\cM_{0,4}$ of the four special points on $\Si_{v_0}$. The meromorphic sections (functions) $\ze_{v_0,1}$ and  $\ze_{v_0,2}$ are given by 
$$
\ze_{v_0,i}(z)=\frac{\beta_i}{z^{s_{\uvec{e}_1,i}}\,(z-1)^{s_{\uvec{e}_2,i}}\,(z-\al)^{s_{\uvec{e}_3,i}}}=\frac{\beta_i\,}{ z\,(z-1)\,(z-\al)} \qquad i=1,2,
$$
where $\beta_i\!\in\!\C^*$. We conclude that 
$$
\aligned
\tn{ob}_{\Gamma}(f)= \Big(\frac{m_1}{m_2},\frac{m_2}{m_3}\Big)\in (\C^*)^2,
\endaligned
$$
where $m_i$ is the slope of the $i$-th line. In particular, it does not depend on $\al$ (which is good). However, if we assume $m_1=m_2=m_3$, i.e. the $3$-lines have the same image, we obtain a configuration for which the step (iii) of the RT-process will be non-trivial. Step (i) will collapse the ghost bubble and $\tn{ob}_{\Gamma}$ does not depend on the location of $4$ special points on that. Step (ii) is trivial, and step (iii) yields just one line with the slop $m=m_1=m_2=m_3\!\in\!\C^*$. Then 
$$
\tn{ob}_\Gamma\colon \cM^{\tn{plog},\gamma}_{g,\mfs}(X,D,A)_\Gamma\cong \C^*\times \cM_{0,4}\lra (\C^*)^2
$$
is the constant map $1$ and $1$ is not a regular value of $\tn{ob}_\Gamma$.
\qed
\end{example}

\begin{example}\label{MCissue_ex}
Let $X\!=\!\P^4$, and $D\!=\!D_1\!\cup\! D_2$ where $D_1$ is a degree $3$ hypersurafce and $D_2$ is a hyperplane. 
The intersection is a cubic surface with $27$ lines in it. Let
$$
g\!=\!0,\quad k\!=\!1,\quad A\!=\![d+a],\quad\tn{and}\quad \mfs\!=\!(s_1)\!=\!((3(d+a),d+a)).
$$ 
Let $\Gamma$ be the configuration in Figure~\ref{alines_fig}; it is the nodal configuration obtained by $a$ lines  (not contained in $D$) each of which intersects $D_1$ and $D_2$ at a single point along $D_{12}$ with tangency order $(3,1)$, together with a degree $d$ rational curve in $D_{12}$ that connects these $a$ lines and carries the marked point $z_1$.
If we assume that the latter is a $d$-fold multiple-cover of a line in $D_{12}$, we obtain a non-simple configuration $\gamma$.
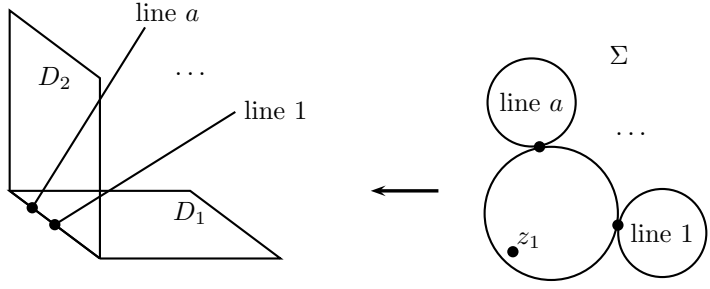
\begin{figure}
\begin{pspicture}(-3,-.5)(11,3)
\psset{unit=.3cm}
\psline(0,0)(8,0)(4,3)(-4,3)(0,0)
\psline(0,0)(0,8)(-4,11)(-4,3)(0,0)
\psline(-2,1.5)(6,6.5)\rput(7.8,6.6){\small{line $1$}}
\rput(4,8.2){\small{$\ldots$}}
\psline(-3,2.25)(2,10.25)\rput(3,10.95){\small{line $a$}}

\rput(4,2){\small{$D_1$}}\rput(-2,8){\small{$D_2$}}
\pscircle*(-2,1.5){.25}
\pscircle*(-3,2.25){.25}

\rput(23,9){\small{$\Si$}}
\pscircle(20,2){3}
\pscircle(19.13,6.92){2} \rput(19.13,6.92){\small{line $a$}}
\pscircle(24.92,1.13){2}\rput(24.92,1.13){\small{line $1$}}
\rput(23.53,5.53){\small{$\ldots$}}
\pscircle*(18.3,0.3){.25}\rput(19,0.8){\small{$z_1$}}

\pscircle*(19.48,4.954){.25}
\pscircle*(22.954,1.48){.25}
\psline[linewidth=.15]{->}(15,3)(12,3)

\end{pspicture}
\caption{A nodal configuration in the boundary of $\ov\cM^{\log}_{0,((3(d+a),d+a))}(\P^4,D,[d+a])$.}
\label{alines_fig}
\end{figure}
The associated map 
$$
\vr\colon \D=\Z^a \oplus \Z^2 \lra \Z^{2a}
$$
has a $1$-dimensional kernel generated by $\big((1,\ldots,1),(4,1)\big)$. Therefore, $\mc{G}(\Gamma)$ is $(a-1)$-dimensional.
A simple calculation shows that the expected dimension of $\cM^{\gamma}_{0,\mfs}(X,D,A)_{\Gamma}$ is $2d+a$. On the other hand, the expected dimension of the main stratum is $d+a+2$. For $d\geq 2$, the former is not smaller than the latter.
In the classical case, $D_{12}$ is a positive manifold and multiple-covers cause no issue. Here, though, the map $\tn{ob}_{\Gamma}$ depends on the covering map and can not be ignored. Note that Condition~(\ref{SP2_eN}) is not satisfied for the line class in $D_{12}$. 
\qed
\end{example}

\noindent
Regarding Problem~\ref{Problem1} in Page~\pageref{Problem1}, the stronger condition of Definition~\ref{NCsemipos_dfnN} allows us to get the upper\footnote{Better upper bounds can be achieved.}  bound
\bEq{UpperDimM_e}
\tn{dim}_\C~\cM_{g,\mfs}^{\tn{plog},\gamma}(X,D,A,\nu)_{\Gamma}\leq Q(\Gamma)= c_1^{TX(-\log D)}(A)+(n-3)(1-g)+k + \tn{dim}_\C(\mc{G})-\dim_\R \K_\R(\Gamma).
\eEq
We get (\ref{UpperDimM_e}), by following steps (i)-(iv) of the RT-process and tracking the change of the quantity
$$
Q(\Gamma)=\sum_{v\in \V} c_1^{TX(-\log D)}(A_v)+k+ (2|\E|-|\uvec{\E}|)- \sum_{v\in \V} |I_{v}|+\sum_{I\subset [N]} (|I|-1)(|\uvec{\E}_I|-|\E_I|).
$$

\noindent
\bIt
\item Collapsing a ghost bubble $C_v=(\Si_{v}=\P^1,z_{v}\cup q_{v})$, with $k_v=|z_{v}|$, $\ell_{v}=|q_v|$, creates a dual graph $\ov\Gamma$ with 
\bEq{GH-Reduction_e}
Q(\Gamma)-Q(\ov\Gamma)= k_v+\ell_v-3=\tn{dim}_\C~\cM_{0,k_v+\ell_v}.
\eEq
Let $\Gamma_{(i)}$ be the result of collapsing all the ghost bubbles in the step (i) of the RT-process. Applying (\ref{GH-Reduction_e}) inductively, we conclude that 
$$
Q(\Gamma)=Q(\Gamma_{\tn{(i)}})+ \tn{dim} \prod_{v\in \V_{\tn{b},\circ}} \cM_{0,k_v+\ell_v}.
$$
\item Replacing a multiple cover bubble $(u_v,C_v=(\Si_{v}=\P^1,z_{v}\cup q_{v}))$, with the underlying simple map $(u_{\ov{v}},C_{\ov{v}}=(\Si_{\ov{v}}=\P^1,z_{\ov{v}}\cup q_{\ov{v}}))$, with $A_v\!=\!d_vA_{\ov{v}}$, creates a dual graph $\ov\Gamma$ with the same set of vertices  and 
\bEq{GH-Reduction_e2}
Q(\Gamma)-Q(\ov\Gamma)= (d_v-1)c_1^{TX(-\log D)}(A_{\ov{v}})+ (k_v+\ell_v)-(k_{\ov{v}}-\ell_{\ov{v}}).
\eEq
By (\ref{SP2_eN}) and (\ref{MC-dim_e}), the righthand side is smaller than or equal to 
$$
\tn{dim}_{\C}\cM_{0,\al_v}(\P^1,z_{\ov{v}}\cup q_{\ov{v}},[d_v])=(d_v-1)(2-\de_{\ov{v}})+ \de_v-\de_{\ov{v}},
$$
where $\de_v$ (resp. $\de_{\ov{v}}$) is the number of contact points of $u_v$ (resp. $u_{\ov{v}}$) with $D$.
Let $\Gamma_{\tn{(ii)}}$ be the result of reducing all the multiple-cover bubbles of $\Gamma_{\tn{(i)}}$ in the step (ii) of the RT-process. Starting with $\Gamma_{\tn{(i)}}$ and applying (\ref{GH-Reduction_e2}) inductively, we conclude that 
$$
Q(\Gamma)\geq Q(\Gamma_{\tn{(ii)}})+ \tn{dim}_{\C} ~\cM_{f''},
$$
where $\cM_{f''}$ is the fiber of (\ref{pigamma_e}) in (\ref{DimMfdp_e}).
\eIt
Step (iii) will further decrease this quantity; therefore, 
$$
Q(\Gamma)\geq Q(\Gamma_{\tn{(iii)}})+ \tn{dim}_{\C} ~\cM_{f''}.
$$
The inequality (\ref{UpperDimM_e}) follows.\\

\noindent
Regarding Problem~\ref{Problem2} in Page~\pageref{Problem2}, we need to replace $\mc{G}(\Gamma)$ with a smaller group $\mc{G}(\gamma)$ for which the obstruction map $\tn{ob}_{\gamma}\colon \cM_{g,\mfs}^{\tn{plog},\gamma}(X,D,A,\nu)_{\Gamma}\lra \mc{G}(\gamma) $ can be transversed and 
$$
\tn{dim}_\C~\cM_{g,\mfs}^{\tn{plog},\gamma}(X,D,A,\nu)_{\Gamma}- \tn{dim}_{\C}~\mc{G}(\gamma)<c_1^{TX(-\log D)}(A)+(n-3)(1-g)+k.
$$
This should be done by selecting a subset $\E^{\star}$ of $\E$ where the equation (\ref{CofG_e}) can be transversed at the nodes corresponding to $\E^{\star}$.
We plan to address these in a future work.

\appendix

%--------------------------------------------------------
\section{Dependence of $\tn{ob}_{\Gamma}$ on ghost bubbles}\label{GH_ss} 
In this section, we study the dependence of the obstruction map $\tn{ob}_{\Gamma}$ on the configuration of special points on a ghost bubble. \\

\noindent
First, we discuss the effect of collapsing a tree of ghost bubbles into one ghost bubble, and further into a multi-node.
Consider a decorated dual graph $\wt\Gamma=(\wt\V,\wt\E)$ associated to a nodal pre-log map 
$$
\wt{f}=\big(\phi_v=\phi|_{\Si_v},u_v,\ze_v=(\ze_{v,i})_{i\in I_v},C_v=(\Si_v,\mfj_v,\vec{z}_v\cup \vec{q}_v)\big)_{v\in \wt\V}
$$ 
as in Definition~\ref{PreLogMap_dfn}. Let $\Gamma'\!=\!(\V',\E')$ be a sub-tree of contracted components in $\wt{\V}_{\tn{b}}$ with $A_v\!=\!0$ for all $v\!\in\!\V'$. We can replace $\Gamma'$ with a single vertex $v_0$ (keeping the decorations unchanged at the rest of the edges) to obtain a new decorated dual graph ${\Gamma}\!=\!({\V},{\E})$ satisfying the combinatorial conditions (\ref{CombCond_e}) (but not necessarily Definition~\ref{LogMap_dfn}.\ref{Tropical_l}). Conversely, starting from $\Gamma$ we obtain various $\wt\Gamma$ by replacing a ghost bubble $v_0$ with a tree of ghost bubbles $\Gamma'$.
Note that 
\bEq{VErel_e}
\wt\V= ({\V} -\{v_0\}) \cup \V'\quad\tn{and}\quad \wt\E={\E} \cup \E'\,.
\eEq
Also, 
$$
I_e=I_{v}=I_{v'}\qquad \forall~v,v'\!\in\!\V',e\!\in\!\E'.
$$
This common value is what $I_{v_0}$ is in ${\Gamma}$. For each $v\!\in\!\V'$, $u_v$ is the constant map into some $p\!\in\!D_{I_{v_0}}\!-\!\partial D_{I_{v_0}}$. The map $u_{v_0}$ is defined to be the constant map into $p$ as well.
Since $v_0\!\in\!\V_{\tn{b}}$ and the combinatorial condition (\ref{CombCond_e}) is satisfied, for any arrangement of special point $z_{v_0}\cup q_{v_0}$ on $\Si_{v_0}$, there exists a set of meromorphic functions
$$
\ze_{v_0}=(\ze_{v_0,i})_{i\in I_{v_0}}
$$
such that 
$$
{f}=\big(\phi_v,u_v,\ze_v,C_v)_{v\in \V-\{v_0\}}\cup\big(u_{v_0},\ze_{v_0},C_{v_0}=(\Si_{v_0},\mfj_{v_0},\vec{z}_{v_0}\cup \vec{q}_{v_0}))  \in \cM^{\tn{plog}}_{g,\mfs}(X,D,A,\nu)_{\Gamma}.
$$

\noindent
The reduced graph $\Gamma$ may not satisfy Condition~\ref{Tropical_l} of Definition~\ref{LogMap_dfn}, however, since $\mc{G}(\Gamma')$ is trivial (see Example~\ref{g0A0_e} and remark~\ref{TofCC_rmk}), Lemma~\ref{TofCC_lmm} below shows that the natural homomorphism from the obstruction group $\mc{G}(\Gamma)$ associated with $\Gamma$ to $\mc{G}(\wt\Gamma)$ is onto. Expanding $v_0$ will increase the kernel of $\vr$ in (\ref{DtoT_e}) and reduces its cokernel.\\

\noindent
Let
$$
\iota_{\D}\colon \D=\Z^{\E}\oplus \bigoplus_{v\in {\V}} \Z^{I_v} \lra \wt\D=\Z^{\wt\E}\oplus \bigoplus_{v\in \wt\V} \Z^{I_v}
$$
be the embedding which maps $\Z^{I_{v_0}}$ diagonally into $ \bigoplus_{v\in \V'} \Z^{I_v}$, and is the identity map on the rest of the terms with respect to the identifications (\ref{VErel_e}). Let 
$$
\iota_{\T}\colon {\T}=\bigoplus_{e\in {\E}} \Z^{I_e} \lra \wt\T=\bigoplus_{e\in \wt\E} \Z^{I_e}
$$
denote the natural inclusion map corresponding to the second identification in (\ref{VErel_e}).

\bLm{TofCC_lmm}
With notation as above, for compatible choices of orientations on $\E$ and $\ov{\E}$, the commutative diagram
$$
\xymatrix{
\D  \ar[d]^{{\vr}} \ar[rr]^{\iota_{\D}} &&  \wt\D  \ar[d]^{\wt\vr} \\
 \T   \ar[rr]^{\iota_{\T}} &&  \wt\T \;.
}
$$
induces a surjective homomorphism $\varphi_{\Gamma,\wt\Gamma}\colon\mc{G}(\Gamma)\!\lra\! \mc{G}(\wt\Gamma)$.
\eLm
\bPf
Since $I_v\!=\!I_{v_0}$ for all $v\in\V'$, the map $\vr'\colon \D'\lra \T'$ corresponding to $\Gamma'$ descends to the similarly denoted map
$$
\vr'\colon \ov\D'=\Z^{\E'}\oplus \frac{\bigoplus_{v\in \V'} \Z^{I_v}}{\Z^{I_{v_0}}}\lra \T'= \bigoplus_{e\in \E'} \Z^{I_e},
$$
where $\Z^{I_{v_0}}\hookrightarrow  \bigoplus_{v\in \V'} \Z^{I_v}$ is the diagonal embedding.
Let 
$$
\pi_{\D}\colon \wt\D\lra \ov\D'\qquad\tn{and}\qquad \pi_{\T}\colon \wt\T\lra \T'
$$
denote the natural projection maps. 
For each $e\!\in\!\wt\E$, let $\uvec{e}$ denote the chosen orientation on $e$ for defining $\vr$. By restriction, this gives us the choice of orientations on $\E$ and $\E'$ used to define $\vr$ and $\vr'$.
The commutative  diagram
$$
\xymatrix{
\D  \ar[rr]^{\iota_{\D}} \ar[d]^{{\vr}} &&  \wt\D \ar[rr]^{\pi_{\D}}  \ar[d]^{\wt\vr} && \ov\D' \ar[d]^{\vr'}\\
 \T   \ar[rr]^{\iota_{\T}} &&  \wt\T \ar[rr]^{\pi_{\T}}  && \T' 
}
$$
has exact rows. Therefore, since $\tn{coker}(\vr')\!=\!0$, we get the long exact sequence,
$$
0\lra \tn{ker}({\vr})\lra \tn{ker}(\wt\vr) \lra \tn{ker}(\vr') \lra \tn{coker}(\vr)\lra \tn{coker}(\wt\vr)\lra 0.
$$
The last map gives us the surjective homomorphism 
$\varphi_{\Gamma,\wt\Gamma}\colon\mc{G}(\Gamma)\!\lra\! \mc{G}(\wt\Gamma)$.
\ePf

\noindent
Suppose $\{f_t\}_{t=1}^{\infty}$ is a sequence of pre-log maps $f_t$ are obtained by the deformations
\bEq{DMConv_e}
(\Si_{v_0}, \vec{z}_{v_0,t}\cup \vec{q}_{v_0,t})_{t=1}^{\infty}\in \cM_{0,k_{v_0}+\ell_{v_0}}
\eEq
of $(\Si_{v_0}, \vec{z}_{v_0}\cup \vec{q}_{v_0})$ and keeping the other components fixed. Suppose that, as $t\!\lra\! \infty$, $f_t$ converges to $\wt{f}\!\in\!\cM^{\tn{plog}}_{g,\mfs}(X,D,A,\nu)_{\wt\Gamma}$ that has the nodal configuration $C'$ in place of $C_{v_0}$. For each $i\!\in\!I_{v_0}$, let $\ze_{v_0,t,i}$ be a meromorphic function on $\Si_{v_0}$ with zeros/poles of the given order at $\vec{z}_{v_0,t}\cup \vec{q}_{v_0,t}$. By \cite[Cr~4.12]{FT1}, for each  $v\in \V'$, as $t\!\lra\!\infty$, restricted to $\Si_{v}$ and up to scaling, $\ze_{v_0,t,i}$ converges to $\ze_{v,i}$. 
We conclude that 
\bEq{limOb_e}
\lim_{t\to \infty} \varphi_{\Gamma,\wt\Gamma}\circ \tn{ob}_{\Gamma}(f_t)=\tn{ob}_{\wt\Gamma}(\wt{f})\in\mc{G}(\wt\Gamma).
\eEq
\vskip.1in

\noindent
Next, we study the dependence of $\tn{ob}_{\Gamma}$ on the location of the special points  $\vec{z}_{v_0}\cup q_{v_0}$ on a single ghost bubble $\Si_{v_0}$. For any pre-log map
$$
f=\big(\phi_v,u_v,\ze_v,C_v)_{v\in \V-\{v_0\}}\cup\big(u_{v_0},\ze_{v_0},C_{v_0}=(\Si_{v_0},\mfj_{v_0},\vec{z}_{v_0}\cup \vec{q}_{v_0}))  \in \cM^{\tn{plog}}_{g,\mfs}(X,D,A,\nu)_{\Gamma}
$$
with domain $\Si=\bigcup_{v\in \V} \Si_v$, let 
\bEq{Ovf_e}
\ov{f}=\big(\phi_v,u_v,\ze_v,C_v)_{v\in \ov\V=\V-\{v_0\}}  \in \cM^{\tn{plog}}_{g,\mfs}(X,D,A,\nu)_{\ov\Gamma}
\eEq
denote the tuple obtained by forgetting the $v_0$-th component, defined on the domain $\ov\Si$ obtained by removing $\Si_{v_0}$ from $\Si$. This is no longer a nodal domain; $\ov\Si$ is a domain with a multi-node $m$ in place of $v_0$. Let $\ov\Gamma=(\ov\V,\ov\E)$ denote the combinatorial type of $\ov\Si$. We have 
$$
\ov\V=\V-\{v_0\},\qquad \ov\E=\big(\E-\{e_1,\ldots,e_{\ell_0}\}\big)\cup \{m\}, \quad \uvec{\ov\E}=\ov\E-\{\ucev{e}_1,\ldots,\ucev{e}_{\ell_0}\}.
$$ 
Note that $I_{v_0}=I_{e_1}=\ldots=I_{e_{\ell_0}}$. We let $I_m$ to also denote this common value. Associated to $\ov\Gamma$ (or in general, any domain with such multi-nodes) we consider the linear map
\bEq{MNvr_e}
\ov\vr\colon \ov\D=\Z^{\E}\oplus \bigoplus_{v\in \ov\V} \Z^{I_v} \lra \ov\T=\bigoplus_{e\in \E-\{e_1,\ldots,e_{\ell_0}\}}\Z^{I_e} \oplus \frac{\bigoplus_{e\in \{e_1,\ldots,e_{\ell_0}\} } \Z^{I_e}}{\Z^{I_m}},
\eEq
where 
$$
\Z^{I_{m}}\lra  \bigoplus_{e\in \{e_1,\ldots,e_{\ell_0}\} } \Z^{I_e}
$$
is the diagonal embedding, and $\ov\vr$ is the composition of the natural inclusion\footnote{Since $\ov\D$ is just missing a summand in $\D$, we get both natural inclusion and projection maps $\ov\D\!\hookrightarrow\! \D$ and $\pi_{\D}\colon\D\!\lra\!\ov\D$, respectively.} inclusion  $\ov\D\!\lra\! \D$, $\vr$, and the natural projection $\pi_{\T}\colon\T\!\lra\! \ov\T$. Similarly to (\ref{Ggroup}), the obstruction group $\mc{G}(\ov\Gamma)$ associated to 
$\ov\Gamma$ is the quotient
$$
\mc{G}(\ov\Gamma)=\bigg(\prod_{e\in \E-\{e_1,\ldots,e_{\ell_0}\}}(\C^*)^{I_e} \times \frac{\prod_{e\in \{e_1,\ldots,e_{\ell_0}\} } (\C^*)^{I_e}}{(\C^*)^{I_m}}\bigg)/~\tn{Image}(\tn{exp}(\ov\vr_{\C})).
$$
Similarly to (\ref{ObGDfn_e}), after fixing a choice of local holomorphic coordinates $w_{\uvec{e}}$ around each nodal point $q_{\uvec{e}}\!\in\!\Si_v$, for all $v\!\in\!\ov\V$ and $\uvec{e}\!\in\!\uvec{\ov\E}_v$, 
we define 
$$
\tn{ob}_{\ov\Gamma}(\ov{f})\in \mc{G}(\ov\Gamma)
$$ 
to be class of 
$$
\aligned
\ov\eta&=
\prod_{e\in \E-\{e_1,\ldots,e_{\ell_0}\}} \prod_{i\in I_{e}} \frac{\eta_{\uvec{e},i}}{\eta_{\scz\ucev{e},i}} \times 
\prod_{i\in I_m} [\eta_{\uvec{e}_1,i},\ldots, \eta_{\uvec{e}_{\ell_0},i}] \\
& \in \prod_{e\in \E-\{e_1,\ldots,e_{\ell_0}\}}(\C^*)^{I_e} \times \big((\C^*)^{\ell_{0}}/\C^*\big)^{I_{m}}
\endaligned
$$
in $ \mc{G}(\ov\Gamma)$.
Here, we use the fact that, for each $i\!\in\!I_{m}$, the class
$$
[\eta_{\uvec{e}_1,i},\ldots, \eta_{\uvec{e}_{\ell_0},i}]\in (\C^*)^{\ell_{0}}/\C^*
$$
is independent of the identification $\cN_XD_i|_{p}=\C$.
For the same reason as in \cite[Lmm~3.7]{FT1}, the class $\tn{ob}_{\ov\Gamma}(\ov{f})\!=\![\ov\eta]$ of $\ov\eta$ in $\mc{G}(\ov\Gamma)$ is independent of the choice of local holomorphic coordinates $w_{\uvec{e}}$ and $\{\ze_v\}_{v\in \ov\V}$ up to rescaling.

\bLm{Indp1_lmm}
The obstruction groups $\mc{G}(\Gamma)$ and $\mc{G}(\ov\Gamma)$ are naturally isomorphic. There exists a holomorphic map $o_{v_0}\colon \cM_{0,k_{v_0}+\ell_{v_0}}\lra \mc{G}(\ov\Gamma)$ such that, under the natural isomorphism $\mc{G}(\Gamma)\!\cong\!\mc{G}(\ov\Gamma)$, we have
$$
\tn{ob}_{\ov\Gamma}(\ov{f})~o_{v_0}(C_{v_0})^{-1}=\tn{ob}_{\Gamma}(f).
$$
\eLm
\noindent
In other words, $\tn{ob}_{\Gamma}(f)=1$ if and only if 
$$
\tn{ob}_{\ov\Gamma}(\ov{f})= o_{v_0}(C_{v_0}).
$$

\bPf
The commutative  diagram
$$
\xymatrix{
\Z^{I_{v_0}}  \ar[rr]^{\iota_{\D}} \ar[d]^{{\cong}} &&  \D \ar[rr]^{\pi_{\D}}  \ar[d]^{\vr} && \ov\D \ar[d]^{\ov\vr}\\
 \Z^{I_m}  \ar[rr]^{\iota_{\T}} &&  \T \ar[rr]^{\pi_{\T}}  && \ov\T
}
$$
has exact rows. We get the long exact sequence,
$$
0\lra \tn{ker}(\vr) \lra \tn{ker}(\ov\vr) \lra 0\lra  \tn{coker}(\vr)\lra \tn{coker}(\ov\vr)\lra 0.
$$
The last map gives us the isomorphism $\varphi_{\Gamma,\ov\Gamma}\colon\mc{G}(\Gamma)\!\lra\! \mc{G}(\ov\Gamma)$.\\

\noindent
By definition, $\varphi_{\Gamma,\ov\Gamma}\circ \tn{ob}_{\Gamma}(f)$ is the class $[\eta]$ of 
$$
\aligned
&\eta=\prod_{e\in \E-\{e_1,\ldots,e_{\ell_0}\}} \prod_{i\in I_{e}} \frac{\eta_{\uvec{e},i}}{\eta_{\scz\ucev{e},i}} \times 
\prod_{i\in I_m} \bigg[\frac{\eta_{\uvec{e}_1,i}}{\eta_{\scz\ucev{e}_1,i}},\ldots, \frac{\eta_{\uvec{e}_{\ell_0},i}}{\eta_{\scz\ucev{e}_{\ell_0},i}}\bigg] \\
& \in \prod_{e\in \E-\{e_1,\ldots,e_{\ell_0}\}}(\C^*)^{I_e} \times \big((\C^*)^{\ell_{0}}/\C^*\big)^{I_{m}}
\endaligned
$$
in $\mc{G}(\ov\Gamma)$.  We have 
$$
\eta=\ov\eta \cdot g^{-1},
$$
where 
$$
\aligned
&g=\prod_{e\in \E-\{e_1,\ldots,e_{\ell_0}\}} \prod_{i\in I_{e}}~1 \times 
\prod_{i\in I_m} \bigg[\eta_{\scz\ucev{e}_1,i},\ldots, \eta_{\scz\ucev{e}_{\ell_0},i}\bigg] \\
& \in \prod_{e\in \E-\{e_1,\ldots,e_{\ell_0}\}}(\C^*)^{I_e} \times \big((\C^*)^{\ell_{0}}/\C^*\big)^{I_{m}}.
\endaligned
$$
Up to rescaling, the section $\ze_{v_0}$ only depends on the location of the special points $z_{v_0}\cup q_{v_0}$ (and the pre-determined multiplicities at those points).
The map
$$
\ze_{v_0}\lra \prod_{i\in I_m} \big[\eta_{\scz\ucev{e}_1,i},\ldots, \eta_{\scz\ucev{e}_{\ell_0},i}\big]\in \big((\C^*)^{\ell_{0}}/\C^*\big)^{I_{m}}
$$
only depends on the $(\C^*)^{I_{m}}$-equivalence class of $\ze_{v_0}$ and, via the inclusion 
$$
\big((\C^*)^{\ell_{0}}/\C^*\big)^{I_{m}}\hookrightarrow \prod_{e\in \E-\{e_1,\ldots,e_{\ell_0}\}}(\C^*)^{I_e} \times \big((\C^*)^{\ell_{0}}/\C^*\big)^{I_{m}}
$$ 
descends to a well-defined map
$$
o_{v_0}\colon \cM_{0,k_{v_0}+\ell_{v_0}}\lra \mc{G}(\ov\Gamma).
$$
\ePf

\noindent
\noindent
The projection map 
\bEq{FB1_e}
\pi_{\Gamma,\ov\Gamma}\colon \cM^{\tn{plog}}_{g,\mfs}(X,D,A,\nu)_{\Gamma}\lra \cM^{\tn{plog}}_{g,\mfs}(X,D,A,\nu)_{\ov\Gamma}, \qquad f\lra \ov{f},
\eEq
is a fiber bundle with fibers $\cM_{0,g_{v_0}+\ell_{v_0}}$. 
By Lemma~\ref{Indp1_lmm},  there exists a map
$$
\tn{ob}_{\ov\Gamma}\colon \cM^{\tn{plog}}_{g,\mfs}(X,D,A,\nu)_{\ov\Gamma}\lra \mc{G}_{\ov\Gamma}
$$
such that 
$$
\pi_{\Gamma,\ov\Gamma}\big(\cM_{g,\mfs}(X,D,A,\nu)_{\Gamma}\big)= \tn{ob}_{\ov\Gamma}^{-1}\big( \tn{Im}(o_{v_0})\big).
$$
In general, the map $o_{v_0}$ can be non-trivial; see Example~\ref{GhostDep_ex} or \ref{BadEx1} below. There are, however, situations where the map is constant and $\tn{ob}_{\Gamma}$ factors through the fiberation (\ref{FB1_e}); see Example~\ref{GoodEx1}. There are also examples, such as Example~\ref{GoodEx2} below, where the map $o_{v_0}$ is constant but the constant value is not $1$.

\begin{example}\label{GoodEx2}
Let $X\!=\!\P^3$, $D\!=\!D_1\!\cup\! D_2\!\cup\!D_3$ be the transverse union of three coordinate hyperplanes, 
$$
g\!=\!0,\quad k\!=\!3,\quad A\!=\![3],\quad  \mfs\!=\!(s_1,s_2,s_3)\!=\!((3,0,0),(0,3,0),(0,0,3)).
$$
Let $\Gamma$ be the dual graph in Figure~\ref{star_fg} with the set of vertices $\V\!=\{v_1,v_2,v_3,v_0\}$ and the set of edges $\E\!=\!\{e_1,e_2,e_3\}$ such that $e_i$ connects $v_0$ and $v_i$, for all $i\!=\!1,2,3$. Choose the orientations $\uvec{e}_i$ to end at $v_0$,  for all $i\!=\!1,2,3$. We have 
$$
\aligned
&I_{v_0}\!=\![N]\!=\!\{1,2,3\},\quad s_{\uvec{e}_1}=(-2,1,1)\!\in\!\Z^3,\quad s_{\uvec{e}_2}=(1,-2,1)\!\in\!\Z^3,\quad s_{\uvec{e}_3}=(1,1,-2)\!\in\!\Z^3, \\
 &I_{v_i}=\{i\}, \qquad A_{v_i}=[1]\!\in\!H_2(X_i,\Z), \quad  \forall i\!=\!1,2,3.
 \endaligned
$$
A pre-log curve with this dual graph is made of a line $\ell_i\!=\!\tn{Im}(u_{v_i})$ in $D_i\cong \P^2$ passing though the point $D_{123}$, for each $i\!\in\![3]$, and a log tuple 
$$
\big(u_{v_0}, [\ze_{v_0,i}]_{i\in [3]}, (\Si_{v_0},\mfj_{v_0})\cong \P^1, q_{v_0}=\{ q_{\scz\ucev{e}_i}\}_{i\in [3]}\big),
$$ 
where $u_{v_0}$ is the constant map onto $D_{123}$, and $\ze_{v_0,i}$, for each $i\!\in\![3]$, is a meromorphic function with poles/zeros of the prescribed order at $q_{v_0}$. Each $\ell_i$ also carries a meromorphic section $\ze_{v_i}$ of $\cO(1)|_{\ell_i}$ with a pole of order $2$ at $q_{\uvec{e}_i}$, and a  zero of order $3$ at the marked point $z_i$.
Both kernel and cokernel of 
$$
\vr\colon \D=\Z^\E \oplus \bigoplus_{i=0}^3 \Z^{I_{v_i}}\lra \T=\bigoplus_{i=1}^3 \Z^{I_{e_i}}
$$ 
are $1$-dimensional. The homomorphism 
\bEq{IsoG_e}
 (\C^*)^{3} \times (\C^*)^{3} \times (\C^*)^{3}\lra \C^*,\qquad \prod_{i=1}^3(x_{i1},x_{i2},x_{i3})\lra \frac{x_{12}x_{23}x_{31}}{x_{13}x_{32}x_{21}}.
\eEq
descends to an isomorphism $\mc{G}\lra \C^*$.
 \\

\begin{figure}
\begin{pspicture}(1,-2.5)(11,1.3)
\psset{unit=.3cm}
\psline[linewidth=0.07](20,-2)(28,-2)
\psline[linewidth=0.07](20,-2)(20,6)
\psline[linewidth=0.07](20,-2)(14.34,-7.66)

\psline[linewidth=0.15](20,-2)(22,2)
\psline[linewidth=0.15](20,-2)(16,1)
\psline[linewidth=0.15](20,-2)(22,-4)

\pscircle*(20,-2){.15}
\rput(24.5,2.5){$D_1$}
\rput(15.5,0){$D_2$}
\rput(22,-6){$D_3$}

\pscircle*(40,-2){.3}\rput(39.7,-1){\small{$v_0$}}
\psline(40,-2)(44,2)\pscircle*(44,2){.3}\rput(45,2.5){\small{$v_1$}}\psline(44,2)(44,4)
\psline[linewidth=0.15]{->}(44,2)(42.5,0.5)\rput(47,.5){\scz{$s_{\uvec{e}_1}=(-2,1,1)$}}

\psline(40,-2)(35,-2)\pscircle*(35,-2){.3}\rput(34,-2.2){\small{$v_2$}}\psline(35,-2)(35,-4)
\psline[linewidth=0.15]{->}(35,-2)(37,-2)\rput(33.5,-1){\scz{$s_{\uvec{e}_2}=(1,-2,1)$}}

\psline(40,-2)(40,-7)\pscircle*(40,-7){.3}\rput(40,-8){\small{$v_3$}}
\psline[linewidth=0.15]{->}(40,-7)(40,-5)\rput(44,-6){\scz{$s_{\uvec{e}_3}=(1,1,-2)$}}
\psline(40,-7)(38,-7)

\end{pspicture}
\caption{Dual graph $\Gamma$ and the image of a map belonging to  $\cM_{0,3}(\P^3,D,[3])_\Gamma$ in $\P^3$.}
\label{star_fg}
\end{figure}
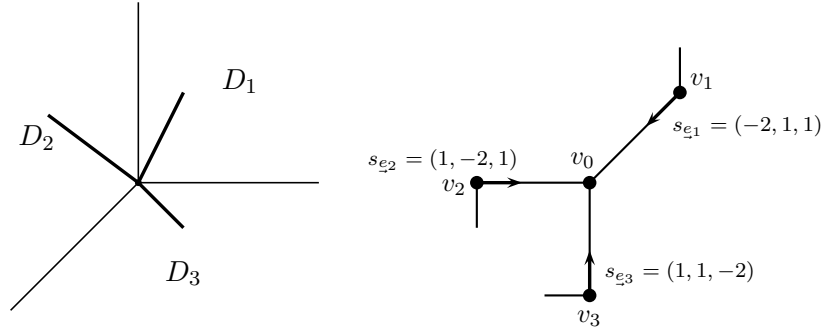

\noindent
In the pre-log space $\cM^{\tn{plog}}_{0,\mfs}(\P^3,D,[3])_\Gamma$, the three lines $\ell_1,\ell_2,\ell_3$ are allowed to be any line passing through the point $D_{123}$ with non-trivial slopes in $\C^*$. For each $i\!\in\![3]$, the line $\ell_i$ is the completion of the image of a map of the form
$$
\C \lra \C^3, \quad w\lra (w_{ij})_{j=1,2,3}\subset \C^3,~~w_{ii}=0,~w_{ij}=a_{ij}w_i,~a_{ij}\!\in\!\C^*,\quad \forall j\!\in\![3]-i.
$$
Here $D_i$ corresponds to the subspace $(x_i=0)\subset \C^3$. For $i=1,2,3$, we have $\ze_{v_i}(w_i)=\al_{ii} w_i^{-2}$.
Putting $(q_{\scz\ucev{e}_1},q_{\scz\ucev{e}_2},q_{\scz\ucev{e}_3})=(0,1,\infty)$, we get
$$
\ze_{v_0,1}= \al_1 \frac{z^2}{z-1},\quad \ze_{v_0,2}= \al_2 \frac{(z-1)^2}{z},\quad \ze_{v_0,3}= \al_3 \frac{1}{z(z-1)}.
$$
We conclude that 
$$
\tn{ob}_{\Gamma}(f)=-\frac{a_{12}a_{23}a_{31}}{a_{13}a_{32}a_{21}};
$$
i.e., $f$ is a log map if and only if the product of the slopes of three lines (in the cyclic order) is $-1$.
On the other hand, 
$$
\tn{ob}_{\ov\Gamma}(\ov{f})=\frac{a_{12}a_{23}a_{31}}{a_{13}a_{32}a_{21}}.
$$ 
We conclude that the image of $o_{v_0}$ in $\mc{G}(\ov\Gamma)\!\cong\!\C^*$ is $-1$.\qed
\end{example}

\noindent
For each expansion $\wt{\Gamma}$ of $\Gamma$ obtained by expending $v_0$ into a tree of ghost bubbles $\Gamma'$ (or in other words, expanding $\Si_{v_0}$ to $C'$), the pre-log space
$$\cM^{\tn{plog}}_{g,\mfs}(X,D,A,\nu)_{\wt\Gamma}
$$ 
is obtained by taking the closure of fibers of (\ref{FB1_e}) and letting $C_{v_0}$ to converge to a nodal domain
\bEq{NCGovG_e}
C'=\bigcup_{v\in \V'} C_v \in \partial \ov\cM_{0,g_{v_0}+\ell_{v_0}}.
\eEq
Let 
$$
\cM_{\Gamma'}\!\subset\!\partial \ov\cM_{0,g_{v_0}+\ell_{v_0}}
$$ 
denote the stratum of nodal configurations (\ref{NCGovG_e}).
Taking union over all such $\wt{\Gamma}$, we get the fiber bundle
$$
\bigcup_{\wt{\Gamma}}\cM^{\tn{plog}}_{g,\mfs}(X,D,A,\nu)_{\wt\Gamma}\lra \cM^{\tn{plog}}_{g,\mfs}(X,D,A,\nu)_{\ov\Gamma}
$$
with compact fibers $\ov\cM_{0,g_{v_0}+\ell_{v_0}}$.  Similarly to Lemma~\ref{Indp1_lmm}, for each $\wt\Gamma$, we get a map 
$$
o_{\Gamma'}\colon \cM_{\Gamma'}\lra \mc{G}(\wt\Gamma)
$$ 
such that 
$$
\tn{ob}_{\wt\Gamma}(\wt{f})=\varphi_{\ov\Gamma,\wt\Gamma}\circ \tn{ob}_{\ov\Gamma}(\ov{f})\cdot  o_{\Gamma'}(C')^{-1}, \qquad  \varphi_{\ov\Gamma,\wt\Gamma}=\varphi_{\Gamma,\wt\Gamma}\circ \varphi_{\Gamma,\ov\Gamma}^{-1}\colon \mc{G}(\ov\Gamma)\lra \mc{G}(\wt\Gamma).
$$
In other words, for the fiber bundle 
$$
\pi_{\wt\Gamma,\ov\Gamma}\colon \cM^{\tn{plog}}_{g,\mfs}(X,D,A,\nu)_{\wt\Gamma}\lra \cM^{\tn{plog}}_{g,\mfs}(X,D,A,\nu)_{\ov\Gamma}, \qquad \wt{f}\lra \ov{f},
$$
we have 
$$
\pi_{\wt\Gamma,\ov\Gamma}\big(\cM_{g,\mfs}(X,D,A,\nu)_{\wt\Gamma}\big)= (\varphi_{\ov\Gamma,\wt\Gamma}\circ \tn{ob}_{\ov\Gamma})^{-1}\big( \tn{Im}(o_{\Gamma'})\big).
$$
Note that as $\Gamma'$ gets bigger, both $\cM_{\Gamma'}$ and $\mc{G}(\wt\Gamma)$ get smaller. 
The Deligne-Mumford convergence in (\ref{DMConv_e}) is compatible with the projection $\mc{G}(\Gamma)\lra \mc{G}(\wt\Gamma)$ in the sense that
$$
\xymatrix{
\cM_{0,k_{v_0}+\ell_{v_0}}  \ar[d]^{\tn{DM-convergence}} \ar[rr]^{o_{v_0}} && \mc{G}(\Gamma)\cong \mc{G}(\ov\Gamma) \ar[d]^{\varphi_{\Gamma,\wt\Gamma}\cong \varphi_{\ov\Gamma,\wt\Gamma}} \\
\cM_{\Gamma'} \ar[rr]^{o_{\Gamma'}} &&  \mc{G}(\wt\Gamma) \;;
}
$$
i.e., if a family of smooth marked curves $\{C_t\}_{t=1}^\infty$ converges to $C'\!\in\!\cM_{\Gamma'}$, then $\pi_{\Gamma,\wt\Gamma}(o_{v_0}(C_t))$ converges to $o_{\Gamma'}(C')$. In Example~\ref{BadEx1}, as $\al$ converges to $0,1,\infty$, we get an expansion $\wt\Gamma$ of $\Gamma$ that replaces $\Si_{v_0}$ by a nodal sphere with $2$ components. Then $\cM_{\Gamma'}$ is just a point and $\mc{G}(\wt\Gamma)$ is the trivial group.

\begin{example}\label{BadEx1}
Let $X\!=\!\P^2$, $D\!=\!D_1\!\cup\! D_2$ be the transverse union of two coordinate lines, 
$$
g\!=\!0,\quad k\!=\!2,\quad A\!=\![4],\quad\tn{and}\quad \mfs\!=\!(s_1,s_2)\!=\!((5,4),(0,1)).
$$ 
Let $\Gamma$ be the configuration in Figure~\ref{2+1line_fig}. It is the nodal configuration obtained by 2 lines passing through the point $D_{12}$, a double cover  $u_{v_1}\colon \Si_{v_1}\lra D_2$ of the line $D_2$ ramified at $D_{12}$ containing the second marked point $z_2$. The three of them are connected by a ghost bubble $\Si_{v_0}$ mapped to $D_{12}$  carrying the marked point $z_1$.
\begin{figure}
\begin{pspicture}(-3,-.5)(11,3)
\psset{unit=.3cm}
\psline(0,0)(8,0)
\psline(0,0)(0,8)\rput(-3.3,6){\small{image of $\Si_{v_1}$}}
\psline(0,0)(6,6)\rput(7.7,6.7){\small{line $3$}}
\psline(0,0)(5,8)\rput(6,8.7){\small{line $2$}}

\rput(9,9){\small{$X$}}
\rput(9,0){\small{$D_1$}}\rput(0,9){\small{$D_2$}}

\rput(23,9){\small{$\Si$}}
\pscircle(20,2){3}\rput(20,2){\tiny{ghost bubble}}
\pscircle(19.13,6.92){2} \rput(19.13,6.4){\small{$\Si_{v_1}$}}
\pscircle(24.92,1.13){2}\rput(24.92,1.13){\small{$\Si_{v_3}$}}
\pscircle(23.53,5.53){2}\rput(23.53,5.53){\small{$\Si_{v_2}$}}
\pscircle*(18.3,0.3){.25}\rput(19,0.8){\small{$z_1$}}
\pscircle*(18.5,8){.25}\rput(19.4,7.8){\small{$z_2$}}
\pscircle*(19.48,4.954){.25}
\pscircle*(22.12,4.12){.25}
\pscircle*(22.954,1.48){.25}
\psline[linewidth=.15]{->}(15,3)(12,3)

\end{pspicture}
\caption{A nodal configuration in the boundary of $\ov\cM^{\log}_{0,((3,3))}(\P^2,D,[3])$.}
\label{2+1line_fig}
\end{figure}
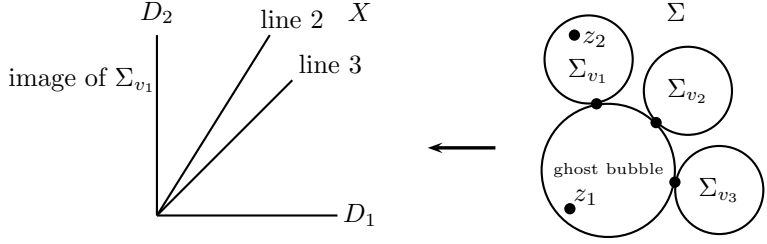
Let $\V=\{v_0,v_1,v_2,v_3\}$ where $v_1$ corresponds to the double cover of $D_2$, $v_2,v_3$ correspond to the lines $2$ and $3$, respectively, and $v_0$ corresponds to the ghost bubble.  Let $\E\!=\!\{e_1,e_2,e_3\}$ where $e_i$ corresponds to the node connecting the domain of the $i$-th component to the ghost bubble. We can choose the orientation $\uvec{e}_i$ to be the one ending at $v_0$. We have 
$$
\aligned
&I_{v_0}=\{1,2\},\quad I_{v_1}=\{2\},\qquad I_{v_i}=\eset \quad\forall~i=2,3,\qquad I_{e_i}=\{1,2\} \quad\forall~i=1,2,3,\\
&s_{\uvec{e}_1}=(2,1),\quad s_{\uvec{e}_2}=(1,1),\quad s_{\uvec{e}_3}=(1,1).
\endaligned
$$
The map 
$$
\vr\colon \Z^{\E} \oplus \Z^{I_{v_0}} \oplus \Z^{I_{v_1}}\lra \bigoplus_{i=1}^3 \Z^{I_{e_i}}
$$
has $1$-dimensional kernel cokernel. The group homomorphism
$$
(\C^*)^{I_{e_1}}\times (\C^*)^{I_{e_2}} \times (\C^*)^{I_{e_2}} \lra \C^*, \qquad ((x_1,y_1),(x_2,y_2),(x_3,y_3))\lra \frac{x_2}{y_2}\frac{y_3}{x_3}.
$$
descends to an isomorphism $\mc{G}\cong \C^*$.
We take the marked point $z_1$ to be $\infty$ and the nodal points $q_{\scz\ucev{e}_1}$, $q_{\scz\ucev{e}_2}$, and $q_{\scz\ucev{e}_3}$ to be $\al$, $0$, and $1$, respectively. Thus, $\al$ parametrizes the configuration space $\cM_{0,4}$ of the four special points on $\Si_{v_0}$. The meromorphic sections (functions) $\ze_{v_0,1}$ and  $\ze_{v_0,2}$ are given by 
$$
\ze_{v_0,i}(z)=\frac{\beta_i\,}{ z^{s_{\uvec{e}_1,i}}\,(z-1)^{s_{\uvec{e}_2,i}}\,(z-\al)^{s_{\uvec{e}_3,i}}}=
\begin{cases}
\beta_1\,\big( z\,(z-1)\,(z-\al)^2\big)^{-1} \qquad &\tn{if}\quad i=1,\\
\\
\beta_2\,\big( z\,(z-1)\,(z-\al)\big)^{-1} \qquad &\tn{if}\quad i=2.
\end{cases}
$$
where $\beta_i\!\in\!\C^*$.
Then
$$
\aligned
\tn{ob}_{\Gamma}(f)=\frac{m_2}{m_3}\frac{\al-1}{\al}\in \C^*.
\endaligned
$$
where $m_2$ and $m_3$ are slopes of the lines $2$ and $3$, respectively. Therefore, we get 
$$
\aligned
&o_{v_0}\colon \cM_{0,4}\lra \C^*, \qquad \al\lra \frac{\al}{\al-1}, \\
&\tn{and}\qquad \tn{ob}_{\ov\Gamma}(\ov{f})=\frac{m_2}{m_3}.
\endaligned
$$
For any $\ov{f}$ with $m_2/m_3\neq 1$, there exists a unique $\al$ such that $f\!\in\!\cM_{0,\mfs}(X,D,A)_{\Gamma}$. The three special values $m_2/m_3=0,1,\infty$ correspond to the limiting situations where $\al$ converges to $0,\infty,1$, respectively, and we get a tree $\Gamma'$ of two ghost bubbles instead of $\Si_{v_0}$.
\qed
\end{example}

\vskip.3in
\noindent
Department of Mathematics,\\
\noindent
The University of Iowa, Iowa City, IA 52242\\
\textit{mohammad-tehrani@iowa.edu}

%------------------------------------------------------------------------------------------------------

\end{document}